\documentclass[12pt]{book}
\usepackage[all]{xy}
\usepackage{amsmath}
\usepackage{amscd}
\usepackage{epsfig}
\usepackage{pstricks}
\usepackage[latin1]{inputenc}

\input{amssym}

\newcommand\qed{\hfill$\sqcap\kern-8.0pt\hbox{$\sqcup$}$}
\newcommand\NN{\bf N}
\newcommand\CC{{\bf C}}

\newcommand\RR{{\bf R}} 
\newcommand\sRR{{\sl \hbox{I\kern-.2em\hbox{R}}}}
\newcommand{\PP}{{\bf P}^k}

\newcommand{\pp}{{\bf P}^1}

\newcommand\ZZ{\bf Z}
\newcommand\proof{\noindent{\em{Proof}.\ }}

\newcommand{\Bif}{{ B}{\textrm{\scriptsize if}}}
\newcommand{\Mis}{{M}{\textrm{\scriptsize is}}}
\newcommand{\SMis}{{M}{\textrm{\scriptsize iss}}}
\newcommand{\Hyp}{{H}{\textrm{\scriptsize yp}}}
\newcommand{\Shi}{{S}{\textrm{\scriptsize hi}}}

\newcommand{\Bc}{T_{\textrm{\scriptsize bif}}}
\newcommand{\Bm}{\mu_{\textrm{\scriptsize bif }}}

\newcommand{\car}{\heartsuit}

\newcommand{\Crf}{{\cal C}_f}
\newcommand{\Crla}{{\cal C}_\la}
\newcommand{\Julf}{{\cal J}_f}
\newcommand{\Julla}{{\cal J}_{\la}}
\newcommand{\Jullo}{{\cal J}_{\la_0}}

\newcommand{\Fatf}{{\cal F}_f}

\newtheorem{theo}{Theorem}[section]
\newtheorem{prop}[theo]{Proposition}
\newtheorem{lem}[theo]{Lemma}
\newtheorem{rem}[theo]{Remark}
\newtheorem{cor}[theo]{Corollary}

\newtheorem{ass}[theo]{Assumptions}
\newtheorem{conj}[theo]{Fatou's conjecture}
\newtheorem{exa}[theo]{Example}

\newtheorem{defi}[theo]{Definition}

\newcommand{\aew}{\textrm{-a.e.}\;}
\newcommand{\Per}{\textrm{Per}}
\newcommand{\Rat}{\textrm{Rat}_d}

\newcommand{\la}{\lambda}

\numberwithin{equation}{section} 

\begin{document}
\begin{center}
\end{center}
\vspace{1cm}

\begin{center}
{\Huge  Bifurcation currents in holomorphic families of rational maps}
\end{center}

\vspace{2cm}
\begin{center}

 {\large \bf Fran\c{c}ois Berteloot}
\end{center}

\vspace{4cm}

\begin{center}
{\large C.I.M.E Course}\\

{\large  Cetraro\;\; July 2011}\\

\end{center}

\newpage
$\;$
\newpage
\vspace*{5cm}
\begin{flushleft}
\begin{minipage}{10cm}

\textsl{ Je d\'edie ce texte \`a mes parents ainsi qu'\`a la m\'emoire de mon ami 
Giovanni Bassanelli.}

 \end{minipage}
\end{flushleft}

\newpage

\bigskip

\begin{center}
{\bf Introduction}
\end{center}
\bigskip
\bigskip
In these lectures we will study  bifurcations within holomorphic families of polynomials or rational maps by mean of ergodic and pluripotential theoretic
 tools. \\

 A  family of rational maps $\left( f_\la\right)_{\la\in M}$, whose parameter space $M$ is a complex manifold, is called a \emph{holomorphic family} if the map $(\la,z)\mapsto f_\la(z)$
is holomorphic on $M\times \pp$ and if the degree of $f_\la$ is constant on $M$.
 The simplest example is the quadratic polynomial family $(z^2 +\la)_{\la \in\CC}$. The space of all rational maps of the same degree
may also be considered as such a family.\\
The interest for bifurcations within holomorphic families of rational maps started in the eighties with the seminal works of Ma\~{n}\'e-Sad-Sullivan \cite{MSS},  Lyubich \cite{Ly2},  and Douady-Hubbard \cite{DoHuN1} and \cite{DoHuN2}. At the end of this decade, McMullen used Ma\~{n}\'e-Sad-Sullivan ideas and Thurston's theory in his fundamental work on  iterative root-finding algorithms \cite{McMullenRF}.\\
In any holomorphic family, the \emph{stability locus} is the maximal open subset of the parameter space on which the Julia set moves continuously with the parameter.
Its complement is called the \emph{bifurcation locus}. In the quadratic polynomial family,
the bifurcation locus is nothing but the boundary of the Mandelbrot set.\\
Ma\~{n}\'e, Sad and Sullivan have shown that the stability locus is dense and  their results also enlighted the (still open) question of the density of hyperbolic parameters.
McMullen proved that any stable algebraic stable family of rational maps is either trivial or affine (i.e. consists of Latt\`es examples), his classification of generally convergent algorithms  follows from this central result.\\

Potential theory has been introduced in the dynamical study of polynomials by Brolin in 1965. These tools and,  more precisely the pluripotential theory developped
after the fundamental works of Bedford-Taylor,  turned out to be extremely powerful to study holomorphic dynamical systems depending on several complex variables.
In this context, the compactness properties of closed positive currents somehow supply to the lack of suitable normality criterions for holomorphic mappings.  We refer to the lecture notes 
\cite{Sibony}, \cite{DS2} by Sibony and Dinh-Sibony for these aspects. As we shall see,  
 these potential-theoretic tools are well adapted for studying bifurcations in one-dimensional holomorphic dynamical systems. The underlying reason for this unity of methods
 is certainly that  holomorphic families are actually holomorphic dynamical systems of the form $M\times\pp \ni (\la,z)\mapsto (\la,f_\la(z))\in M\times\pp$.\\

The use of potential theory in the study of parameter spaces actually started rather soon. Indeed, as revealed by the work of Douady-Hubbard \cite{DoHu} and Sibony \cite{S} around 1980, the Green function of $(z^2 +\la)$ evaluated at the critical value plays a crucial role in the study of the Mandelbrot set. The relation between this quantity and the Lyapunov exponent was also known and Przytycki explicitely raised the problem to ''understand the connections between Lyapunov exponent characteristic and potential theory for rational mappings" (\cite{Prz}). Let us also mention that Ma\~{n}\'e proved that the Lyapunov exponent of a  rational map and the Hausdorff dimension of its maximal entropy measure are closely related: their product is equal to the logarithm of the degree of the map  (\cite{Mane}).\\
Around 2000,  a decisive achievement was made by DeMarco. She generalized Przytycki's formula and proved that, in any holomorphic family of rational maps, the bifurcation locus 
is the support of a $(1,1)$ closed positive current admiting the Lyapunov exponent function as a global potential  (\cite{DeM1},
\cite{DeM2}). This current, denoted $\Bc$,  is now called the \emph{bifurcation current}.\\
In the recent years, several authors (\cite{BB1}, \cite{Ph}, \cite{DF}, \cite{dujardin2},  \cite{BuEp}, \cite{BB2}, \cite{BB3}, \cite{Gau}, \cite{BG})
have investigated the geometry of the bifurcation locus using the current $\Bc$ and its exterior powers $\Bc^k:=\Bc\wedge\Bc\wedge\cdot\cdot\cdot\wedge\Bc$.
We shall present here most of the results obtained in these papers. They go from laminarity statements for certain regions of the bifucation locus to
Hausdorff dimension estimates and include precise density (or equidistribution) properties relative to various classes of specific parameters.
Let us stress that a common feature of these papers is to use the supports of the $\Bc^k$ to enlight a certain stratification of the bifurcation locus.\\
We have not discussed bifurcation theory for families of endomorphisms of higher dimensional complex projective spaces but have mentionned, among the techniques presented in these notes, those which also work in this more general context. These aspects have been considered by Bassanelli, Dupont, Molino and the author (\cite{BB1}, \cite{BDM}), by Pham \cite{Ph} and also appear in Dinh-Sibony's survey \cite{DS2}. Let us mention that Deroin and Dujardin \cite{DeDu} have recently extended these methods  for studying bifurcations in the context of Kleinian groups.\\

Let us now describe the contents of these notes more precisely. Chapter 1 mainly deals with the construction of the Green measure (the maximal entropy measure) for a fixed rational map.
We show that this measure enjoys good ergodic properties and prove a very important approximation formula for its Lyapunov exponent. This formula, which roughly says that the Lyapunov exponent can be computed on cycles of big period, will be crucial for investigating the structure of the bifurcation current. The classical results presented in this chapter are essentially due to Lyubich. The approximation formula (and its generalization to higher dimensional rational maps) is due to Dupont, Molino and the author.\\

Chapter 2 has a double purpose. One is to present some concrete holomorphic families on which we shall often work, the other is to introduce the
hypersurfaces $Per_n(w)$ whose distribution turns out to shape the bifurcation locus. We study polynomial families and prove an important result
due to Branner-Hubbard about the compactness of the connectedness locus. As a by-product, we obtain some informations which will later be decisive to control the global behaviour
of the bifurcation current in such families. We also study the moduli space $Mod_2$ of degree two rational maps, that is the space of holomorphic conjugacy classes of degree two rational maps. To some extent, this space can be considered as a holomorphic family and presents some common features with the family of cubic polynomials.
Following Milnor, we show that $Mod_2$ can be identified in a dynamically natural way to $\CC^2$.\\
 In a fixed holomorphic family, $Per_n(w)$ is defined as the set of parameters
for which the corresponding map has a cycle of multiplier $w$ and exact period $n$. We show that the $Per_n(w)$ are complex hypersurfaces (co-dimension one complex analytic subsets) in the parameter space and give some precise global defining functions for them. This will be an important tool to study the equidistribution of the bifurcation current by the hypersurfaces $Per_n(w)$.\\

The core of Chapter 3 is devoted to the definition of the bifurcation current $\Bc$. It starts with a brief survey of Ma\~{n}\'e-Sad-Sullivan work and a motivated introduction of
$\Bc$. We then formally define the bifurcation current, we also introduce the activity currents of marked critical points and relate them to $\Bc$. Most of the chapter is devoted to  DeMarco's fundamental results. We prove that $\Bc$ is supported by the bifucation locus and admits both the Lyapunov exponent function and the sum of the Green function evaluated on critical points as a global potential. We end the chapter with a proof of DeMarco's formula which precisely relates the Lyapunov exponent with the Green function
evaluated on the critical points. Our presentation of these results  only relies on computations with closed positive currents. It is due to Bassanelli and the author who also generalized it to higher dimensional rational maps. \\

In Chapter 4, we study how the asymptotic distribution of dynamically defined hypersurfaces is governed by the bifurcation current. Most of the results presented here are valid in polynomial
families and in $Mod_2$ but some are actually true in any holomorphic family. We are mainly interested by the hypersurfaces $Per_n(w)$ (for $\vert w\vert \le 1$) and the hypersurfaces $Per(c,k,n)$ defined by the pre-periodicity of a critical point $c$ (i.e. $f^n(c)=f^k(c)$). We prove that, conveniently weighted, these hypersurfaces equidistribute the bifurcation current. The results concerning $Per(c,k,n)$ are due to Dujardin and Favre while those concerning $Per_n(w)$ are due to Bassanelli and the author, we have tried here to unify their presentation.
These results yield a precise description of the laminated structure of the bifurcation locus in large regions of $Mod_2$. More precise laminarity results, due to Dujardin, are  presented at the end of the chapter.\\

In Chapter 5, we investigate the higher exterior powers $\Bc^k$ of the bifurcation current. Although we prove some general results which indicate that the supports of 
$\Bc^k$ induce a dynamically meaningful stratification of the bifurcation locus, we concentrate our attention on the highest power, a measure, denoted $\Bm$ and called the \emph{bifurcation measure}. Our main goal is to show that the support of $\Bm$ is the seat of the strongest bifurcations. To this purpose we first prove that this support is simultaneously approximated by extremely stable parameters (hyperbolic) and highly unstable ones (Shishikura, Misiurewicz). The approximation by hyperbolic or Shishikura parameters is due to Bassanelli and the author, the approximation by Misiurewicz parameters is due to Dujardin and Favre.   We then characterize the support of $\Bm$ as being the closure of the sets of Shishikura or Misiurewicz parameters. The proof of this result, which is due to Buff and Epstein, requires to introduce new ideas based on transversality techniques. These methods have recently been extended by Gauthier for proving that the support of the bifurcation measure, in the moduli space of degree $d$ rational maps, has full Hausdorff
dimension. We end this last chapter with a presentation of Gauthier's proof.\\

We have tried to give a synthetic and self-contained presentation of the subject. In most cases, we have given complete and detailed proofs and, sometimes have substantially simplify those available in the literature. Although basics on ergodic theory are discussed in the first chapter, we have not included any element of pluripotential theory. For this we refer the reader to the appendix of Sibony's \cite{Sibony} and Dinh-Sibony's \cite{DS2} lecture notes or to the book of Demailly \cite{Demailly}. We finally would like to recommend the recent survey of Dujardin \cite{DuS} on bifurcation currents and equidistribution in parameter space.\\

{\bf Acknowledgements} It is my pleasure to thank my colleagues Charles Favre, Thomas Gauthier and Nessim Sibony for their useful comments on the first draft of these notes.
I also would like to thank the anonymous referee for having carefully read the manuscript and having helped me to improve it

\tableofcontents

\chapter{Rational maps as ergodic dynamical systems}

Among ergodic properties of rational maps we present those which will be used in our study. A particularly important result of this first chapter is an approximation
formula for the Lyapunov exponent of a rational map with respect to its maximal entropy measure.

\section{Potential theoretic aspects}
\subsection{The Fatou-Julia dichotomy}

A rational function 
$f$ is a holomorphic map of the Riemann sphere to itself and may be represented as the ratio of two polynomials
\begin{center}
$f=\frac{a_0+a_1z+a_2z^{2}+\cdot\cdot\cdot +a_d z^{d}}{{b_0+b_1 z+b_2z^{2}+\cdot\cdot\cdot + b_d z^{d}}}$
\end{center}
where at least one of the coefficents $a_d$ and $b_d$ is not zero. The number $d$ is the algebraic degree of $f$.
In the sequel we shall more likely speak of {\it rational map}. 
Such a map may also be considered as a holomorphic ramified self-cover of the Riemann sphere whose topological degree
is equal to $d$. Among these maps, polynomials are exactly those for which $\infty$ is totally invariant:
$f ^{-1}\{\infty\}=f\{\infty\}={\infty}$.
\\

It may also be convenient to identify the Riemann sphere with the one-dimensional complex projective space $\pp$
that is the quotient of $\CC^{2}\setminus\{0\}$ by the action $z\mapsto u\cdot z$ of $\CC^{*}$.
Let us recall that the Fubini-Study form $\omega$ on $\pp$ satisfies $\pi^{*}(\omega)= dd^{c} \ln \Vert\;\Vert$ where 
the norm is the euclidean one on $\CC^{2}$.\\
In this setting, the map $f$ can be seen as induced on $\pp$ by a non-degenerate and $d$-homogenous map of $\CC^{2}$
\begin{center}
$F(z_1,z_2):=\big(a_0 z_2^{d}+a_1z_1z_2^{d-1}+\cdot\cdot\cdot +a_d z_1^{d},
b_0 z_2^{d}+b_1z_1z_2^{d-1}+\cdot\cdot\cdot +b_d z_1^{d}\big)$
\end{center}
through the canonical projection $\pi: \CC^{2}\setminus\{0\}\to \pp$.

The homogeneous map $F$ is called a lift of $f$; all other lifts are proportional to $F$.\\

A point around which a rational map $f$ does not induce a local biholomorphism is called \emph {critical}.
The image of such a point is called a \emph {critical value}. A degree $d$ rational map has exactly $(2d-2)$ critical points
counted with multiplicity.
The \emph{critical set} of $f$ is the collection of all critical points and is denoted $\Crf$.\\

As for any self-map, we may study the dynamics of rational ones, that is trying to understand the behaviour of the sequence of iterates

\begin{center}
$f^{n}:=f\circ\cdot\cdot\cdot\circ f$.
\end{center}

The Fatou-Julia dynamical dichotomy consists in a splitting of $\pp$ into two disjoint subsets on which the dynamics
of $f$ is radically different.
 The \emph{Julia set}  of a rational map $f$ is the subset of $\pp$ on which the dynamics of $f$ may drastically change under a small perturbation
of initial conditions while the \emph{Fatou set}  is the complement of the Julia set.

\begin{defi} The Julia set $\Julf$ and the Fatou set $\Fatf$ of a rational map $f$ are respectively defined by:
\begin{center}
$\Julf:=\{z\in\pp\;/\;(f^{n})_n\;\textrm{is not equicontinuous near}\; z\}$\\
$\Fatf:=\pp\setminus\Julf.$
\end{center}
\end{defi}

Both the Julia and the Fatou set are totally invariant: $\Julf=f\big(\Julf\big)=f^{-1}\big(\Julf\big)$
and $\Fatf=f\big(\Fatf\big)=f^{-1}\big(\Fatf\big)$. In particular, $f$ induces two distinct dynamical systems
on $\Julf$ and $\Fatf$.\\
The dynamical system $f:\Julf\to\Julf$ is chaotical. However,
as it results from the  Sullivan non-wandering theorem and the Fatou-Cremer classification,
 the dynamics of a rational map is, in some sense, totally predictible on its Fatou set.\\
 
 The periodic orbits are called \emph{cycles} and play a very important role in the understanding of the global dynamical behaviour of a map $f$.
 
 \begin{defi}
 A $n$-cycle is a set of $n$ distinct points $z_0,z_1,\cdot\cdot\cdot,z_{n-1}$ such that $f(z_i)=z_{i+1}$ for $0\le i\le n-2$
 and $f(z_{n-1})=z_0$. One says that $n$ is the exact period of the cycle.
 \end{defi}
 
 Each point $z_i$ is fixed by $f^{n}$. The \emph{multiplier} of the cycle is the derivative of $f^{n}$
 at some point $z_i$ of the cycle and computed in a local chart: $\big( \chi \circ f^n\circ \chi^{-1}\big)'(\chi(z_i))$. It is easy to see that this number depends only on the cycle and neither on the point $z_i$ or the chart $\chi$. By abuse we shall denote it $(f^{n})'(z_i)$.\\
 The local dynamic of $f$ near a cycle is governed by the multiplier $m$. This leads to the following
 
 \begin {defi}
 The multiplier of a $n$-cycle is a complex number $m$ which is equal to the derivative of $f^{n}$ computed in any local chart at any point of the cycle.
 \begin{itemize}
 \item[]When $\vert m\vert >1$ the cycle is called \emph{repelling}
\item[] when $\vert m\vert <1$ the cycle is called \emph{attracting}
\item[] when $\vert m\vert =1$ the cycle is called \emph{neutral}.
 \end{itemize}
 \end{defi}
 
 Repelling cycles belongs to the Julia set and attracting ones to the Fatou set. For neutral cycles
 this depends in a very delicate way on the diophantine properties of the argument of $m$.\\
 
 The first fundamental result about Julia sets is the following.
 
 \begin{theo}\label{theoFJ}
 Repelling cycles are dense in the Julia set.
 \end{theo}
 
It is possible to give an elementary proof of that result using the Brody-Zalcman renormalization technique (see \cite{BM}). We shall see later that repelling cycles
actually equidistribute a measure whose support is exactly the Julia set.
 
\subsection{The Green measure of a rational map}

Our goal  is to endow the dynamical system $f: \Julf\to\Julf$ with an ergodic structure
 capturing most of its chaotical nature. This is done by exhibiting an invariant measure $\mu_f$ on $\Julf$ which
 is of constant Jacobian. Such a  measure was first constructed by Lyubich \cite{Lyu}. For our purpose it will 
 be extremely important to use a potential-theoretic approach which goes back to Brolin \cite{Bro}
 for the case of polynomials. We follow here the presentation given by Dinh and Sibony in their survey
 \cite{DS2}  which also covers \emph {mutatis mutandis} the construction of Green currents for holomorphic
 endomorphisms of ${\bf P}^{k}$. \\

The following Lemma is the key of the construction. It relies on the fundamental fact that $$d^{-1}f^{\star} \omega =\omega +dd^{c}v$$
for some smooth function $v$ on $\pp$. This follows from a standard cohomology argument or may be seen concretely by setting 
$v:=d^{-1}\ln \frac{\Vert F(z)\Vert}{\Vert z\Vert^{d}}$ for some lift $F$ of $f$.\\

\begin{lem}\label{lemexg}
Up to some additive constant, there exists a unique continuous function $g$ on $\pp$ such that
$d^{-n}f^{n\star}\nu \to dd^{c}g +\omega$ for any positive measure $\nu$ which is given by 
$\nu=\omega +dd^{c}u$ where $u$ is continuous.
\end{lem}

\proof 
Let us set $g_n:= v+\cdot\cdot\cdot +d^{-n+1}v\circ f^{n-1}$. One sees by induction that 
$d^{-n}f^{n\star}\nu=\omega + dd^{c}g_n + dd^{c} (d^{-n}u\circ f^{n})$. As the sequence $(g_n)_n$
is clearly uniformly converging, the conclusion follows by setting $g:=\lim_n g_n$.\qed\\

It might be useful to see how the function $g$ can be obtained by using lifts. 

\begin{lem}\label{PropG_F}
Let $f$ be a degree $d$ rational map. For any lift $F$ of $f$, the sequence $d^{-n}\ln \Vert F^{n}(z)\Vert$ converges uniformly on compact subsets of 
$\CC^{2}\setminus\{0\}$ to a function $G_F$ which satisfies the following invariance and homogeneity properties:
\begin{itemize}
\item[i)] $G_F\circ F=d G_F$
\item[ii)] $G_F(tz)=G_F(z)+ln \vert t\vert$,\;$\forall t\in \CC$.
\end{itemize}  
Moreover,  $G_F -\ln \Vert\;\Vert=g\circ \pi$ where $g$ is a function given by Lemma \ref{lemexg}.

\end{lem}

\proof Let us set $G_n(z):=d^{-n} \ln \Vert F^n (z)\Vert$. As $F$ is homogeneous and non-degenerate there exists a constant $M>1$ such that
\begin{center}
$\frac{1}{M} \Vert z\Vert^d\le \Vert F(z)\Vert\le M \Vert z\Vert^d$.
\end{center}
Thus $\frac{1}{M} \Vert F^n (z)\Vert^d\le \Vert F^{n+1} (z)\Vert\le M \Vert F^n (z)\Vert^d$ which, taking logarithms and dividing by 
$d^{n+1}$ yields $\vert G_{n+1}(z)-G_n(z)\vert \le \frac{\ln M}{d^{n+1}}$. This shows that $G_n$ is uniformly converging to $G_F$.
The properties i) and ii) follows immediately from the definition of $G_F$.\\

According to the proof of Lemma \ref{lemexg},  
$g= \lim_n \big(v+\cdot\cdot\cdot +d^{-n+1}v\circ f^{n-1}\big)$ where a possible choice of $v$ is
$v\circ \pi=d^{-1}\ln \frac{\Vert F(z)\Vert}{\Vert z\Vert^{d}}$. To get the last assertion, it suffices to observe that
$d^{-k}v\circ f^{k}\circ \pi=G_{k+1}-G_k$.\qed\\

The two above lemmas lead us to coin the following

\begin{defi}\label{DefiGreen}
Let $F$ be a lift of a degree $d$ rational map $f$. The Green function $G_F$ of $F$ on $\CC^2$ is defined by
\begin{center}
 $G_F:=\lim_n d^{-n}\ln \Vert F^{n}(z)\Vert$.
\end{center}
The Green function $g_F$ of $f$ on $\pp$ is defined by
\begin{center}
$G_F -\ln \Vert\;\Vert=g_F\circ \pi$.
\end{center}
We will often commit the abuse to denote $g_f$ any function which is equal to $g_F$ up to some additive constant and to call it \emph{the} Green function of $f$.
\end{defi}

The function $G_F$ is $p.s.h$ on $\CC^2$ with a unique pole at the origin.\\

It is worth emphasize that both $g_F$ and $G_F$  are uniform limits of
smooth functions. 
In particular, these functions are continuous.
One may actually prove more  (see \cite{DS2} Proposition 1.2.3 or \cite{BB1} Proposition 1.2):

\begin{prop}\label{PropGHC}
The Green functions $G_F (z)$ and $g_F(z)$ are H\"{o}lder continuous in $F$ and $z$.
\end{prop}

We are now ready to define the Green measure $\mu_f$ and verify its first properties.

\begin{theo} Let $f$ be a degree $d\ge 2$ rational map and $g$ be the Green function of $f$. Let $\mu_f:=\omega +dd^{c}g$.
Then $\mu_f$ is a $f$-invariant probability measure whose support is equal to $\Julf$. Moreover $\mu_f$ has constant Jacobian:
$f^{\star}\mu_f= d\mu_f$.
\end{theo}

\proof That $\mu_f$ is a probability measure follows immediately from Stokes theorem since $\int_{\pp} \omega=1$.\\

We shall use Lemma \ref{lemexg} for showing that  $f^{\star}\mu_f= d\mu_f$. By construction $v+d^{-1}g\circ f=\lim_n v+d^{-1}g_n\circ f=\lim_n g_{n+1}=g$ and thus $d^{-1}f^{\star}\mu_f =\omega +dd^{c}v +d^{-1}dd^{c}(g\circ f)= \omega+
dd^{c}(v+d^{-1}g\circ f)=\mu_f$.\\

The invariance property $f_{\star}\mu_f =\mu_f$ follows immediately from $f^{\star}\mu_f= d\mu_f$ by using the fact that $f_{\star}f^{\star}=d\;Id$.\\

Let us show that the support of $\mu_f$ is equal to $\Julf$. If $U\subset \Fatf$ is  open then $f^{n\star}\omega$ is uniformly bounded on $U$ and therefore
$\mu_f(U)=\lim\int_U d^{-n}f^{n\star}\omega=0$, this shows that $Supp\;\mu_f\subset \Julf$. Conversely the identity
$f^{\star}\mu_f =d \mu_f$ implies that $\left(Supp\;\mu_f\right)^{c}$ is invariant by $f$ which, by Picard-Montel's theorem, implies that
$\left(Supp\;\mu_f\right)^{c}\subset \Fatf$. One should first observe that, since $\mu_f$ has continuous potentials, $Supp\;\mu_f$ certainly contains more than three points.
\qed\\
 
It is sometimes useful to use the Green function $G_F$ for defining local potentials of $\mu_f$.

\begin{prop}\label{locpot}
Let $f$ be a rational map and $F$ be a lift. For any section $\sigma$ of the canonical projection $\pi$ defined on some open subset $U$ of $\pp$, the function $G_F\circ \sigma$ is a potential for $\mu_f$ on $U$. 
\end{prop}

\proof On $U$ one has
$dd^{c} G_F\circ \sigma=dd^{c} g_F + dd^{c} \ln \Vert\sigma\Vert = dd^{c} g_F + \sigma^{*} dd^{c} \ln \Vert\cdot\Vert=
dd^{c} g_F + (\pi\circ\sigma)^{*} \omega=
dd^{c} g_F+\omega=\mu_f$.\qed\\

Let us underline that, by construction, the measure $\mu_f$ has continuous local potentials and in particular cannot give mass to points.\\

 The measure $\mu_f$ can also be obtained as the image by the canonical projection $\pi: \CC^{2}\setminus\{0\}\to \pp$ of a Monge-Amp\`ere
measure associated to the positive part $G_F^{+}:=\max\left(G_F,0\right)$ of the Green function $G_F$.

 \begin{prop}\label{Propmu_F}
 Let $F$ be a lift of a degree $d$ rational map $f$ and $G_F$ be the Green function of $F$.
 The measure $\mu_F:=dd^{c} G_F^{+}\wedge dd^{c}G_F^{+}$ is supported on the compact set
 $\{G_F=0\}$ and satisfies $F^{\star}\mu_F=d^{2}\mu_F$ and $\pi_{\star}\mu_F=\mu_f$.
 \end{prop}
 
  This construction will be used only once in this text and we therefore skip its proof. 
Observe that the support of $\mu_F$ is contained in the boundary of the compact set $K_F:=\{G_F\le 0\}$ which is precisely the set of points $z$ with bounded forward orbits by $F$.\\
 
The case of polynomials presents very interesting features, in particular the Green measure coincides with the harmonic measure of the
filled-in Julia set.

\begin{prop}\label{Greenpoly}
Let $P$ be a degree $d\ge 2$ polynomial on $\CC$. The Green function $g_{\footnotesize{\CC},P}$ of $P$ is the subharmonic function defined on the complex plane by
\begin{center}
$g_{\footnotesize{\CC},P}:=\lim_n d^{-n} \ln ^{+}\vert P^{n}\vert$.
\end{center}
The Green measure of $P$ is compactly supported in the complex plane and admits $g_{\footnotesize{\CC},P}$ as a global potential there.
\end{prop}

\proof 
We may take $F:=\big(z_2^{d} P(\frac{z_1}{z_2}),z_2^{d}\big)$ as a lift  of $P$.
Then $F^{n}:=\big(z_2^{d^n}P^n(\frac{z_1}{z_2}),z_2^{d^n}\big)$ and 
\begin{eqnarray*}
G_F(z_1,1)=\frac{1}{2} \lim_n d^{-n} \ln \big( 1+\vert P^{n} (z_1)\vert^{2} \big)=
\lim_n d^{-n} \ln ^{+}\vert P^{n}(z_1)\vert.
\end{eqnarray*}
The conclusion then follows from Proposition \ref{locpot}. \qed\\
 
\section{Ergodic aspects}
\subsection{Mixing, equidistribution towards the Green measure}

\begin{defi} Let $(X,f,\mu)$ be a dynamical system. One says that the measure $\mu$ is mixing if and only if 
\begin{center}
$\lim_n \int_X (\varphi\circ f^{n})\;\psi\;\mu=\int_X \varphi\;\mu\; \int_X \psi\;\mu$
\end{center}
for any test functions $\varphi$ and $\psi$.
\end{defi}
This means that the events $\{f^{n}(x)\in A\}$ and $\{x\in B\}$ are asymptotically independants for any pair of Borel sets
$A,B$.\\

As we shall see, the constant Jacobian property implies that Green measures are mixing and therefore ergodic.

\begin{theo}
 The Green measure $\mu_f$ of any degree $d$ rational map $f$ is mixing.
\end{theo}

\proof
Let us set $c_{\varphi}:=\int\varphi\;\mu_f$ and $c_{\psi}:=\int\psi\;\mu_f$
where $\varphi$ and $\psi$ are two test functions. We may assume that $c_{\varphi}=1$.\\

Since $\mu_f$ and $\varphi\mu_f$ are two probability measures, there exists a smooth function $u_{\varphi}$ on $\pp$
such that:
\begin{eqnarray}\label{cohomes}
\varphi\mu_f=\mu_f +\Delta u_{\varphi}.
\end{eqnarray}

On the other hand, by the constant Jacobian property $f^{*}\mu_f=d \mu_f$ we have: 

\begin{eqnarray}\label{CJ}
d^{-n}f^{n*}\big((\varphi -c_{\varphi})\mu_f\big) =\big(\varphi\circ f^{n}  -c_{\varphi}\big)\mu_f
\end{eqnarray}

Now, combining \ref{cohomes} and \ref{CJ} we get:

\begin{eqnarray*}
\int (\varphi\circ f^{n})\psi\;\mu_f - \big(\int \varphi\;\mu_f\big)\big( \int \psi\;\mu_f\big)=
\int (\varphi\circ f^{n})\psi\,\mu_f - c_{\varphi} c_{\psi}=\\
\int \psi\big(\varphi\circ f^{n}  -c_{\varphi}\big)\;\mu_f= \int \psi\; d^{-n}f^{n*}\big((\varphi -c_{\varphi}) \mu_f\big)
=\int \big(d^{-n}f^{n}_*\psi\big) \;( \varphi -1)\;\mu_f =\\
\int \big(d^{-n}f^{n}_*\psi\big) \Delta u_{\varphi}
=\int \psi\; d^{-n}f^{n*}(\Delta u_{\varphi}) =\int \psi\;d^{-n}\Delta(u_{\varphi}\circ f^{n})=
d^{-n}\int (u_{\varphi}\circ f^{n})\;\Delta \psi.
\end{eqnarray*}

As $\int (u_{\varphi}\circ f^{n})\;\Delta \psi$ is bounded, this leads to the desired conclusion.\qed\\

It is not hard to show that a mixing measure is also ergodic.

\begin{defi} Let $(X,f,\mu)$ be a dynamical system. One says that the measure $\mu$ is ergodic if and only if 
all integrable $f$-invariant functions are constants. 
\end{defi}

In particular this allows to use the classical Birkhoff ergodic theorem
which says that time-averages along typical orbits coincide with the spatial-average:

\begin{theo}\label{theobirk}
Let $\left(X,f,\mu\right)$ be an ergodic dynamical system and $\varphi\in L^1(\mu)$. Then
\begin{center}
$ \lim_n \frac{1}{n} \sum_{k=0}^{n-1} \varphi(f^k(x)) = \int_{X}  \varphi\;\mu$
\end{center}
for $\mu$ almost every $x$ in $X$.
\end{theo}

The measure-theoretic counterpart of Fatou-Julia theorem \ref{theoFJ} is the following equidistribution result which has been first proved by Lyubich \cite{Lyu}. The content of subsection \ref{ssLyapMult} will provide another proof which exploits the mixing property.

\begin{theo} \label{theoLyu} Let $f$ be a rational map of degree $d$. Let $R_n^{\star}$ denote the set of all $n$-periodic repelling points
of $f$. Then 
$d^{-n}\sum_{R_n^{\star}}\delta_z$ is weakly converging to $\mu_f$.
\end{theo}

Let us finally mention another classical equidistribution result.
We refer to \cite{DS2} for a potential theoretic proof.

\begin{theo}
Let $f$ be a degree $d$ rational map. Then
\begin{center}
$\lim_n d^{-n}\sum_{\{f^{n}(z)=a\}} \delta_z = \mu_f$
\end{center}
for any $a\in \pp$ which is not exceptional for $f$.
\end{theo}

We recall that a map $f$ has \emph{no} exceptional point unless $f$ is a polynomial or a map of the form  $z^{\pm d}$.  A map of the form $z^{\pm d}$ has two exceptional points : $0$ and $\infty$. For polynomials other than $z^d$, the only exceptional point is $\infty$.\\

Although this will not be used in this text, we mention that the Green measure $\mu_f$ of any degree $d$ rational map
$f$ is the unique measure of maximal entropy for $f$. This means that the entropy of $\mu_f$ is maximal and, according to the variational principle, equals $\ln d$ which is the value of the topological entropy of $f$.

\subsection{Natural extension and iterated inverse branches}\label{ssNatExt}
To any ergodic dynamical system, it is possible to associate a new system  which is invertible and contains all the information of the original one. It is basically obtained by considering the set of all complete orbits on which is acting a shift.
 This general construction is the so-called \emph {natural extension} of a dynamical system;
here is a formal definition.\\

\begin{defi} The natural extension of a dynamical system $\left(X,f,\mu\right)$ is the dynamical system 
$\left(\widehat X,\hat f,\hat \mu\right)$ where
\begin{center}
$\widehat{X}:=\{\hat x:=(x_n)_{n\in\ZZ}\;/\;x_n\in {X},\;f(x_n)=x_{n+1}\}$\\
$\hat f(\hat x):=(x_{n+1})_{n\in\ZZ}$\\
$\hat\mu \{(x_n)\;\textrm{s.t.}\; x_0\in B\}=\mu (B)$.
\end{center}
The canonical projection $\pi_0: \hat X\to X$ is given by $\pi_0(\hat x)=x_0$. One sets $\tau$ for $(\hat f)^{-1}$.
\end{defi}

Let us stress that $\pi_0  \circ \hat f=f\circ \pi_0$
and $(\pi_0)_* (\hat {\mu})=\mu$. The measure $\hat \mu$ inherits most of the ergodic properties of $\mu$.\\

\begin{prop}\label{propmixnatex}
The measure $\hat\mu$ is ergodic (resp. mixing) if and only if 
$\mu$ is ergodic (resp. mixing).
\end{prop}

We refer the reader to the chapter 10 of \cite{CFS}  for this construction and its properties.\\

 A powerful way to control the behaviour of inverse branches along typical orbits of the system $\left(\Julf,f,\mu_f\right)$ is to apply
 standard ergodic theory to its natural extension. This is what we shall do now. The first point is to observe that one may work
with orbits avoiding the critical set of $f$. To this purpose one
 considers
\begin{center}
${\widehat X}_{reg}=\{\hat x\in\hat{\Julf}\;/\;x_n\notin \Crf\;;\;\forall n\in \ZZ\}$.
\end{center}
As  $\hat\mu_f$ is $\hat f$-invariant and $\mu_f$ does not give mass to points, one sees that
$\hat \mu_f({\widehat X}_{reg})=1$.

\begin{defi}\label{defiInv}
Let $\hat x\in{\widehat X}_{reg}$ and $p\in\ZZ$. The injective map induced by $f$ on some neighbourhood of $x_p$
is denoted 
$f_{x_p}$.
The inverse of $f_{x_p}$
is defined on some neighbourhood of $x_{p+1}$ and is denoted $f^{-1}_{x_p}$. We then 
set
\begin{center}
$f^{-n}_{\hat x}:=f^{-1}_{x_{-n}}\circ\cdot\cdot\cdot\circ f^{-1}_{x_{-1}}$. 
\end{center}
The map $f^{-n}_{\hat x}$ is called \emph{iterated inverse branch  of $f$ along $\hat x$ and of depth $n$}.
\end{defi}

It will be good to keep in mind that $f^{-n}_{\hat x}(x_0)=x_{-n}$ and that $f^{-1}_{\hat x}=f^{-1}_{x_{-1}}$.\\
 
 The following Proposition yields a control of the disc on which $f^{-1}_{ x_{-k-1}}$ is defined.
 
\begin{prop}\label{propDomInv}
For any sufficently small and strictly positive $\epsilon$,
there exists a function $\alpha_\epsilon:{\widehat X}_{reg} \to ]0,1[$ such that
\begin{center}
$\alpha_{\epsilon}(\tau({\hat x}))\ge e^{-\epsilon} \alpha_{\epsilon}({\hat x})$ and
\end{center}
\begin{center}
 $f^{-1}_{ x_{-k-1}}$ is defined on $ D(x_{-k},\alpha_{\epsilon}(\tau^k(\hat x)))$
 \end{center}
  for $\hat\mu_f\aew\hat x \in{\widehat X}_{reg}$ and every $k\in\ZZ$.\\
\end{prop}

The function $\alpha_{\epsilon}$ is a so-called \emph{slow function}. The interest of such a function relies on the fact that its decreasing  might be negligeable with respect to other datas. For instance, in some circumstances, Proposition \ref{propDomInv}
will tell us that the local inverses 
$f^{-1}_{ x_{-k}}$ are defined on discs whose radii may essentially be considered as constant along the orbit $\hat x$.\\

\proof We need the following quantitative version of the inverse mapping theorem (see \cite{BrDu} lemme 2).\\

\begin{lem}\label{FactBD} Let $\rho(x):=\vert f'(x)\vert$ , $r(x):=\rho(x)^2$. There exists $\epsilon_0>0$ and, for $\epsilon\in ]0,\epsilon_0]$, 
$0<C_1(\epsilon),C_2(\epsilon)$ such that for every $x\in \Julf$:
\begin{itemize}
\item[1-] $f$ {\it is one-to-one on} $D\left(x,C_1(\epsilon)\rho(x)\right)$,
\item[2-]  $D(f(x),C_2(\epsilon) r(x))\subset f\left[D\left(x,C_1(\epsilon)\rho(x)\right)\right]$,
\item[3-] $\textrm{Lip}\;f^{-1}_x \le e^{\frac{\epsilon}{3}} \rho(x)^{-1}$ {\it on} $D(f(x),C_2(\epsilon) r(x))$.
\end{itemize}
\end{lem}

Let us set $\beta_{\epsilon}(\hat x):= \textrm{Min}\;(1,C_2(\epsilon) r(x_{-1}))$. According to the two first assertions of the above Lemma, $
f^{-1}_{ x_{-1}}=f^{-1}_{\hat x}$ is defined on $D(x_0,\beta_{\epsilon}(\hat x))$ and, similarly, $
f^{-1}_{ x_{-k-1}}=f^{-1}_{\tau^k(\hat x)}$ is defined on $D(x_{-k},\beta_{\epsilon}(\tau^k(\hat x)))$.
All we need is to find a function $\alpha_{\epsilon}$ such that $0<\alpha_{\epsilon}<\beta_{\epsilon}$ and 
$\alpha_{\epsilon}(\tau({\hat x}))\ge e^{-\epsilon} \alpha_{\epsilon}({\hat x})$.\\
As $\mu_f$ admits continuous local potentials, the function  $\ln \beta_{\epsilon}$ is $\hat\mu_f$-integrable. Then, by Birkhoff ergodic theorem \ref{theobirk},
$ \int_{\widehat X} \ln \beta_{\epsilon}\; \hat\mu_f=\lim_{\vert n\vert\to +\infty} \frac{1}{\vert n\vert}\sum_{k=1}^{n} \ln \beta_{\epsilon}(\tau^k(\hat x))$ and, in particular
\begin{center}
$\lim_{\vert n\vert \to +\infty} \frac{1}{\vert n\vert} \ln \beta_{\epsilon}(\tau^n(\hat x))=0\;\;\textrm{for}\; \hat \mu_f\aew  \hat x \in\widehat X$.
\end{center}
In other words, for $\hat\mu_f\aew\hat x\in {\widehat X}_{reg}$ there exists $n_0(\epsilon,\hat x)\in\NN$ such that 
$\beta_{\epsilon}(\tau^n (\hat x))\ge e^{-\vert n\vert \epsilon}$ for $\vert n\vert \ge n_0(\epsilon,\hat x) $. Setting then 
$V_{\epsilon}:= \inf_{\vert n\vert \le n_0(\epsilon,\hat x)} \big(\beta_{\epsilon}(\tau^n (\hat x))e^{\vert n\vert \epsilon}\big)$ we obtain
 a measurable function $V_{\epsilon}:{\widehat X}_{reg} \to ]0,1]$ such that:
$\beta_{\epsilon}(\tau^n(\hat x)) \ge e^{-\vert n\vert\epsilon} V_{\epsilon} (\hat x)$ for $\hat\mu_f\aew$ $\hat x \in{\widehat X}_{reg}$ and every $n\in\ZZ$.
It suffices to take $\alpha_{\epsilon}(\hat x):=\textrm{Inf}_{n\in \ZZ} \{\beta_{\epsilon}(\tau^n(\hat x))  e^{\vert n\vert\epsilon}\} $.\qed\\

\section{The Lyapunov exponent}

\subsection{Definition, formulas and some properties}\label{sslyap}

Let us consider the ergodic dynamical system $\left(\Julf,f,\mu_f\right)$ which has been constructed in the last section. As the measure $\mu_f$ has continuous local potentials,
the function $\ln\vert f'\vert$ belongs to $L^1(\mu_f)$ for any choice of a metric $\vert\;\vert$ on $\pp$. We may therefore apply the Birkhoff ergodic theorem \ref{theobirk} to get:
\begin{eqnarray}\label{LyapLim}
\lim_n \frac{1}{n} \ln\vert (f^n)'(z)\vert
 = \lim_n \frac{1}{n} \sum_{k=0}^{n-1} \ln \vert f'(f^k(z))\vert = \int_{\pp} \ln \vert f'\vert\;\mu_f,\;\;\mu_f\aew
\end{eqnarray}

This identity shows that the integral $\int_{\pp} \ln \vert f'\vert\;\mu_f$ does not depend on the choice of the metric $\vert\;\vert$ and leads to the following definition.
\begin{defi}\label{defilyap}
The Lyapunov exponent of the ergodic dynamical system $\left(\Julf,f,\mu_f\right)$ is the number
\begin{center}
$L(f)=\int_{\pp} \ln \vert f'\vert\;\mu_f$.
\end{center}
For simplicity we shall say that $L(f)$ is  the Lyapunov exponent of $f$.
\end{defi}

As the identity \ref{LyapLim} shows, the Lyapunov exponent $L(f)$ is the exponential rate of growth of $\vert (f^n)'(z)\vert$ for a typical $z\in \Julf$.\\

\begin{rem}
Using the invariance property $f_{*}\mu_f=\mu_f$ one immediately sees that $L(f^n)=nL(f)$.
\end{rem}

We shall need an expression of $L(f)$ which uses the formalism of line bundles. In order to prove it, we first
compare the Lyapunov exponents of $f$ with the sum of Lyapunov exponents of one of its lifts $F$.

\begin{prop}\label{propJon}
Let $F$ be a lift of some rational map $f$of degree $d$. Then the sum of Lyapunov exponents of $F$ with respect to $\mu_F$ is
given by $L(F):=\int \ln \vert \det F'\vert\;\mu_F$ and is equal to 
$L(f)+\ln d$.
\end{prop}

\proof
let $F$ be a polynomial lift of $f$. We shall compute the Lyapunov exponent using the spherical metric $\vert\;\vert_{s}$. Exploiting the fact that $f^*\omega=\vert f'\vert_{s}^2 \omega$, it is not difficult to check that 
\begin{center}
$|f'(\xi)|_s =\displaystyle\frac{1}{d}\frac{\Vert z\Vert^2}{\Vert F(z)\Vert^2}\left|\det F'(z)\right|$
\end{center}\label{lmL}
for any $z$ such that $\pi(z)=\xi$. We thus have
\begin{eqnarray*}
\frac{1}{n}\ln|(f^n)'(\xi)|_s +\ln d=\frac{1}{n}\ln\frac{\Vert z\Vert^2}{\Vert F^n(z)\Vert^2}+\frac{1}{n}\ln\left|\det (F^n)'(z)\right|.
\end{eqnarray*}

Then the conclusion follows by Birkhoff theorem \ref{theobirk} since $\Vert F^n(z)\Vert$ stays away from $0$ and $+\infty$ when $z$ is in the support of $\mu_F$ and 
$\pi_{\star} \mu_F=\mu_f$.\qed\\

 For any integer $D$ the line bundle ${\cal O}_{\pp}(D)$ over $\pp$ is the quotient of $(\CC^2\setminus\{0\})\times\CC$ by the action 
 of  $\CC^*$ defined by  
$(z,x)\mapsto (uz,u^Dx)$. We denote by  $[z,x]$ the elements of this quotient.\\

The canonical metric on ${\cal O}_{\pp}(D)$ may be written
\begin{center}
$\|[z,x]\|_0=e^{-D\ln\|z\|}|x|$. 
\end{center}
The homogeneity property of  $G_F$ allows us to define another metric on  ${\cal O}_{\pp}(D)$ by setting
\begin{center}
$\|[z,x]\|_{G_F}=e^{-D G_F(z)}|x|$.
\end{center}

Let us underline that, according to Definition \ref{DefiGreen}, $\Vert\cdot\Vert_{G_F}=e^{-Dg_F}\Vert\cdot\Vert_0$.\\

The expression given in the following Lemma will allow us to perform integration by parts and get some fundamental formulas.
This will turn out to be extremely useful when we shall relate the Lyapunov exponent with bifurcations in Chapter \ref{chapBC}.

\begin{lem}\label{lemLyapBund}
Let $f$ be a rational map of degree $d\ge 2$ and $F$ be one of its lifts. Let $D:=2(d-1)$
and $Jac_F$ be the holomorphic section of ${\cal O}_{\pp}(D)$ induced  by $\det F'$.
Then
\begin{center} 
 $L(f)+\ln d=\displaystyle\int_{\pp}\ln\|Jac_F\|_{G_F}\;\mu_f$.
\end{center}
\end{lem}

\proof
The section $Jac_F$ is defined by $Jac_F(\pi (z)):=[z,\det F']$ for any $z\in\CC^2\setminus\{0\}$. 
Using Proposition \ref{propJon}, the fact that $G_F$ vanishes on the support of $\mu_F$ and $\pi_*\mu_F=\mu_f$ we get
\begin{eqnarray*}
%\begin{center}
L(f)+\ln d=\int_{\tiny\CC^2}\ln|\det F'|\;\mu_F=\int_{\tiny \{G_F=0\}}\ln|\det F'|\;\mu_F=\\
\int_{\{G_F=0\}}\ln\left(e^{-D G_F(z)}|\det F'|\right)\;\mu_F=
\int_{\tiny\CC^2}\ln\|Jac_F\circ\pi\|_{G_F}\;\mu_F=\\
\int_{\pp}\ln\|Jac_F\|_{G_F}\;\mu_f.
\end{eqnarray*}
%\end{center}
\qed\\

It is an important and not  obvious fact that $L(f)$ is stricly positive. It actually follows from
the Margulis-Ruelle inequality that $L(f)\ge \frac{1}{2}\ln d$ where $d$ is the degree of $f$.
We will present later a simple argument which shows that this bound is equal to $\ln d$ for polynomials (see Theorem \ref{PrzFor}).
Zdunik and Mayer (\cite{Zd}, \cite{Ma}) have proved that the bound $\frac{1}{2}\ln d$ is taken if and only if the map $f$ is a
Latt\`es example. Let us summarize these results in the following statement.

\begin{theo}\label{theoCarLat} The Lyapunov exponent of a degree $d$ rational map is always greater than $\frac{1}{2}\ln d$
and the equality occurs if and only if the map is a Latt\`es example.
\end{theo}

We recall that a Latt\`es map is, by definition, induced on the Riemann sphere from an expanding map on a complex torus by mean of some elliptic function.
We refer to the survey paper of Milnor \cite{Mi1} for a detailed discussion of these maps.\\

A remarkable consequence of the positivity of $L(f)$ is that the iterated inverse branch  $f^{-n}_{\hat x}$ (see definition \ref{defiInv}) are approximately $e^{-nL}$-Lipschiptz
and are defined on a disc whose size
only depends  on $\hat x$.

\begin{prop}\label{PropBrInv}
There exists $\epsilon_0>0$ and, for $\epsilon\in ]0,\epsilon_0]$, two  measurable functions $\eta_\epsilon:{\widehat X}_{reg}\to ]0,1]$ and  $S_\epsilon:{\widehat X}_{reg} \to ]1,+\infty]$ such that
the maps
$f^{-n}_{\hat x}$ are defined on $D(x_0,\eta_\epsilon(\hat x))$ and
$\textrm{Lip}\;f^{-n}_{\hat x}\le S_\epsilon (\hat x) e^{-n(L-\epsilon)}$
for $\hat\mu_f\aew\hat x\in {\widehat X}_{reg}$ and for every $n\in \NN$. 
\end{prop}

\proof 
 We may assume that $0<\epsilon_0<\frac{L}{3}$.
 Since $f^{-n}_{\hat x}=f^{-1}_{x_{-n}}\circ\cdot\cdot\cdot\circ f^{-1}_{x_{-1}}$, the third assertion of Lemma \ref{FactBD} yields
$\ln \textrm{Lip}\;f^{-n}_{\hat x} \le n\frac{\epsilon}{3} -\sum_{k=1}^n \ln \rho(x_{-k}) $. By Birkhoff ergodic theorem we thus have
\begin{center}
$\limsup \frac{1}{n}\ln \textrm{Lip}\;f^{-n}_{\hat x}  \le -L+\frac{\epsilon}{3}$\;  for\; $\hat \mu_f\aew\hat x \in\widehat X$.
\end{center}
Then there exists
$n_0(\hat x)$ such that $\textrm{Lip}\;f^{-n}_{\hat x} \le e^{-n(L-\epsilon)}$ for $n\ge n_0(\hat x)$ and
it suffices to set $S_{\epsilon}:=\max_{0\le n\le n_0(\hat x)} \left(e^{n(L-\epsilon)}\textrm{Lip}\;f^{-n}_{\hat x}\right)$
to get the estimate 
\begin{center}
$\textrm{Lip}\;f^{-n}_{\hat x}\le S_\epsilon (\hat x) e^{-n(L-\epsilon)}$\;
for every \;$n\in \NN$ \;and\; $\hat\mu_f\aew\hat x\in {\widehat X}_{reg}$.
\end{center}

We now set $\eta_{\epsilon}:=\frac{\alpha_{\epsilon}}{S_{\epsilon}}$ where $\alpha_{\epsilon}$ is the given by Proposition \ref{propDomInv}. Let us check by induction on $n\in\NN$ that
$ f^{-n}_{ \hat x}$ is defined on $D(x_0,\eta_{\epsilon}(\hat x))$ for $\hat\mu_f\aew\hat x \in\widehat X$ and every $n\in\NN$.
Here we will use the fact that the function $\alpha_{\epsilon}$ is slow: $\alpha_{\epsilon}(\tau({\hat x}))\ge e^{-\epsilon} \alpha_{\epsilon}({\hat x})$.\\
Assume that $f^{-n}_{\hat x}$ is defined on $D\big(x_0,\eta_{\epsilon}(\hat x)\big)$. Then,
by our estimate on $\textrm{Lip}\;f^{-n}_{\hat x}$, we have 
\begin{center}
$f^{-n}_{\hat x}\big(D\big(x_0,\eta_{\epsilon}(\hat x)\big)\big)
\subset D\big(x_{-n}, e^{-n(L-\epsilon)} \alpha_{\epsilon}(\hat x)\big)$.
\end{center}
 On the other hand, by Proposition \ref{propDomInv}, 
the branch $f_{x_{-n-1}}^{-1}$ is defined on the disc $D\big(x_{-n}, \alpha_{\epsilon}(\tau^{n+1}(\hat x))\big)$ which, as $\alpha_{\epsilon}$
is slow, contains $D\big(x_{-n}, e^{-(n+1)\epsilon}\alpha_{\epsilon}(\hat x)\big)$. Now, since $0<\epsilon_0 <\frac{L}{3}$
one has $e^{-(n+1)\epsilon}  \ge e^{-n(L-\epsilon)}$ and thus $f^{-(n+1)}_{\hat x}=f_{x_{-n-1}}^{-1}\circ f^{-n}_{\hat x}$ is defined on $D\big(x_0,\eta_{\epsilon}(\hat x)\big)$.
\qed\\

\subsection{Lyapunov exponent  and multipliers of repelling cycles}\label{ssLyapMult}
 
The following approximation property will play an important role in our study of bifurcation currents. We would like to mention that Deroin and Dujardin have 
recently used similar ideas to study the bifurcation in the context of Kleinian groups (see \cite{DeDu}).

\begin{theo}\label{theoapprox}
 Let $f:\pp\to\pp$ be a rational map of degree $d\ge 2$ and $L$ the Lyapunov exponent of $f$ with respect to its Green measure.
Then:
\begin{center}
$L=\lim_n d^{-n}\sum_{p\in R^{*}_n} \frac{1}{n} \ln \vert (f^n)'(p)\vert$
\end{center}
where $R^{*}_n:=\{p\in \pp \;/\; \;p\; \textrm{has exact period}\; n \;  \textrm{and}\;\vert (f^{n })' (p)\vert > 1 \}$.
\end{theo}

Observe that the Lyapunov exponent $\lim_{k} \frac{1}{k}\ln  \vert (f^k)'(p)\vert$ of $f$ along the orbit of a point $p$
is precisely equal to $\frac{1}{n} \ln \vert (f^n)'(p)\vert$ when $p$ is $n$ periodic. The above Theorem thus shows that the Lyapunov exponent $L$ of $f$ is the limit, when 
$n\to+\infty$, of the averages of Lyapunov exponents of repelling $n$-cycles.\\

To establish the above Theorem, we will prove that the repelling cycles equidistribute the Green measure $\mu_f$ in  a somewhat constructive way and control the multipliers of the cycles which
appear.  For this purpose, we follow the approach used by Briend-Duval \cite{BrDu} 
where the positivity of the Lyapunov exponent plays a crucial role and which also works for endomorphisms of $\PP$.
This strategy actually yields to a version of Theorem \ref{theoapprox} for endomorphisms of $\PP$; this has been done in 
\cite{BDM}. Okuyama has given a different proof of Theorem \ref{theoapprox} in \cite{Oku,Oku2}, his proof actually does not use the positivity of the Lyapunov exponent.
The proof we present here is that of \cite{Be} with a few more details.\\

\proof 
For the simplicity of notations we consider polynomials and therefore work on ${\bf C}$ with the euclidean
metric. We shall denote $D(x,r)$ the open disc centered at $x\in {\bf C}$ and radius $r>0$. From now on, $f$ is a degree $d\ge 2$ polynomial whose Julia set is denoted 
${ J}$ and whose Green measure is denoted $\mu$. \\

We shall use the natural extension (see subsection \ref{ssNatExt}) and exploit the positivity of $L$ through Proposition \ref{PropBrInv}.
Let us add a few notations to those already introduced in Propositions \ref{propDomInv} and \ref{PropBrInv}.
Let $0<\epsilon_0$ be given by Proposition \ref{propDomInv}.

For $0<\epsilon\le \epsilon_0$ and $n,N\in\NN$ we set:
\begin{center}
${\widehat X}^{\epsilon}_N:=\{\hat x \in \widehat X \;/\; \eta_{\epsilon} (\hat x)\ge \frac{1}{N}\;\textrm{and}\; S_{\epsilon} (\hat x)\le N\}$\\

${\hat\nu}^{\epsilon}_N:=1_{{\widehat X}^{\epsilon}_N} \hat \mu$\\

$\nu^{\epsilon}_N:=\pi_{0\star}\hat \nu^{\epsilon}_N.$
\end{center}

For $0<\epsilon\le L$ and $n,N\in\NN$ we set:
\begin{center}
$R^{\epsilon }_n:=\{p\in {\bf C} \;/\; f^n (p)=p\;\textrm{and}\;\vert (f^{n })' (p)\vert \ge e^{n(L-\epsilon)}\}$\\
$\mu^{\epsilon}_n:=d^{-n}\sum_{R^{\epsilon }_n} \delta_p$
\end{center}
\begin{center}
$R_n:=R_n^L=\{p\in {\bf C}\;/\; f^n (p)=p\;\textrm{and}\; \vert (f^{n})' (p)\vert \ge 1\}$\\
$\mu_n:=\mu^{L}_n=d^{-n}\sum_{R_n} \delta_p.$
\end{center}

%We have to show that  $\mu^{L}_n \to \mu$ and $L=\lim_n\frac{d^{-n}}{n}\sum_{R_n^*} \ln\vert (f^{n})'(p)\vert$. \\

The following Lemma reduces the problem to some estimates on Radon-Nikodym derivatives.

\begin{lem}\label{lemRND}
If any weak limit $\sigma$ of $\big(\mu^{\epsilon}_n\big)_n$ for $\epsilon \in ]0,\epsilon_0[$ satisfies
$\frac{d\sigma}{d\nu^{\epsilon'}_N}\ge 1$ for some $\epsilon'>0$ and every $N\in \NN$ then
$\mu^{L}_n \to \mu$ and $L=\lim_n\frac{d^{-n}}{n}\sum_{R_n^*} \ln\vert (f^{n})'(p)\vert$. 
\end{lem}

\proof We start by showing that
$\mu^{\epsilon}_n \to \mu$ for any $\epsilon \in\; ]0,L]$. Let $\sigma$ be a weak limit of $\big(\mu^{\epsilon}_n\big)_n$.
Since all the $\mu^{\epsilon}_n$ are probability measures, it suffices to show that $\sigma=\mu$.\\

Assume first that $0<\epsilon< \epsilon_0$.
By assumption $\frac{d\sigma}{d\nu^{\epsilon'}_N}\ge 1$ and therefore $\sigma\ge \nu^{\epsilon'}_N$ for every $N\in\NN$.
Letting $N\to +\infty$ one gets $\sigma\ge \mu$.
This actually implies that $\sigma=\mu$ since
\begin{center}
$\sigma({J}) \le \limsup_n \mu^{\epsilon}_n({J})\le \lim_n \frac{d^n+1}{d^n} =1= \mu(J)$.
\end{center}
We have shown that $\mu^{\epsilon}_n\to \mu$ for
 $0<\epsilon < \epsilon_0$. Let us now assume that $\frac{\epsilon_0}{2}=:\epsilon_1\le \epsilon$.\\
As $\mu^{\epsilon}_n\ge \mu^{\epsilon_1}_n $ and $\mu^{\epsilon_1}_n\to \mu$, one gets
$\sigma \ge \mu$. Just as before this implies that $\sigma=\mu$.\\

We now want to show that $L=\lim_n\frac{d^{-n}}{n}\sum_{R_n^*} \ln\vert (f^{n})'(p)\vert$.
Let us set $\varphi_n(p):=\frac{1}{n} \ln\vert (f^{n})'(p)\vert$. For  $M>0$ one has
\begin{eqnarray*}
\mu^{\epsilon}_n({J}) (L-\epsilon)\le d^{-n} \sum_{R^{\epsilon}_n} \varphi_n(p) \le d^{-n} \sum_{R_n} \varphi_n(p)= \int_{J} \ln \vert f'\vert \mu_n\le\\
\le \int_{J} Max \big(\ln \vert f'\vert, -M\big)  \mu_n
\end{eqnarray*}
since $\mu^{\epsilon}_n\to \mu$ and $\mu_n =\mu^{L }_n\to\mu$ we get
\begin{eqnarray*}
 (L-\epsilon)\le \liminf d^{-n} \sum_{R_n} \varphi_n(p)\le \limsup d^{-n} \sum_{R_n} \varphi_n(p) \le\\
 \int_{J} Max \big(\ln \vert f'\vert, -M\big)  \mu.
\end{eqnarray*}
To obtain  $\lim d^{-n} \sum_{R_n} \varphi_n(p)=L$ it suffices to make first $M\to +\infty$ and then $\epsilon\to 0$.\\
Since there are less than $2n d^\frac{n}{2}$ periodic points whose period strictly divides $n$, one may replace $R_n$ by
$R^{*}_n:=\{p\in \pp \;/\; \;p\; \textrm{has exact period}\; n \;  \textrm{and}\;\vert (f^{n })' (p)\vert \ge 1 \}$.\qed\\

Let us now finish the proof of Theorem \ref{theoapprox}. We assume here that $0<\epsilon < \frac{\epsilon_0}{2}$. Let $\hat a\in {\widehat X}^{\epsilon}_N$ and $a:=\pi (\hat a)$. For every $r>0$ we denote by
$D_r$ the \emph{closed} disc centered at $a$ of radius $r$. 
According to Lemma \ref{lemRND}, it suffices to show that \emph{any weak limit} $ \sigma$  of $ \big(\mu^{2\epsilon}_n\big)_n$ satisfies
\begin{eqnarray}\label{estiRN}
 \sigma({D}_{r'})\ge 
 \nu^{\epsilon}_N({D}_{r'}),\;\; \textrm{for any integer}\; N\;\textrm{and all} \;0<r'<\frac{1}{N}.
\end{eqnarray}

Let us pick $r'<r<\frac{1}{N}$. We set $\widehat D_r:=\pi^{-1} (D_r)$ and :
\begin{eqnarray*}
\widehat C_n:=\{\hat x\in \widehat D_r \cap \widehat{X}_{reg\;N}^{\epsilon}\;/\; f^{-n}_{\hat x} (D_r)\cap D_{r'} \ne \emptyset \}.
\end{eqnarray*}
Let also consider the collection $S_n$ of sets of the form $f^{-n}_{\hat x} (D_r)$ where $\hat x$ runs in $\widehat C_n$. As $f^{-n}_{\hat x}$ is an inverse branch on $D_r$ of the ramified cover $f^n$,
one sees that the sets of the collection $S_n$ are mutually disjoint.\\

Let us momentarily admit the two following estimates:
\begin{eqnarray} \label{Fact1}
 d^{-n} \big(\textrm{Card}\;  S_n\big) \le \mu^{2\epsilon}_n (D_r) \;\textrm{ for}\; n\;\textrm{ big enough}
\end{eqnarray}
\begin{eqnarray}\label{Fact2}
d^{-n} \big(\textrm{Card}\; S_n\big) \;\mu(D_r) \ge \hat \mu\big(\hat f^{-n}(\widehat D_r\cap \widehat{X}_{reg\;N}^{\epsilon}) \cap \widehat D_{r'}\big).
\end{eqnarray}

Combining \ref{Fact1} and \ref{Fact2} yields:
\begin{eqnarray*}
\hat \mu\big(\hat f^{-n}(\widehat D_r\cap \widehat{X}_{reg\;N}^{\epsilon}) \cap \widehat D_{r'}\big)\le \mu (D_r) \mu^{2\epsilon}_n (D_r)
\end{eqnarray*}
which, by the mixing property of $\hat\mu$, implies
\begin{eqnarray*}
\nu^{\epsilon}_N(D_r) \mu(D_{r'}) =\hat \mu (\widehat D_r\cap \widehat{X}_{reg\;N}^{\epsilon}) \hat\mu(\widehat D_{r'}) \le \mu (D_r) \sigma (D_r)
\end{eqnarray*}
since $\mu(D_{r'})>0$,  one gets \ref{estiRN} by making $r\to r'$ .\\

Let us now prove the estimate \ref{Fact1}. We have to show that $D_r$ contains at least $\big(\textrm{Card}\;S_n\big)$ elements of $R^{2\epsilon}_n$ when $n$ is big enough. Here we shall use Proposition \ref{PropBrInv}.
 For every $\hat x\in \widehat C_n \subset \widehat{X}_{reg\;N}^{\epsilon}$ one has $\eta_{\epsilon}(\hat x) \ge \frac{1}{N}$ and
$S_{\epsilon}(\hat x) \le N$ and thus the map $f^{-n}_{\hat x}$ is defined on $D_r$ ($r<\frac{1}{N}$) and $\textrm{Diam}\;f^{-n}_{\hat x}(D_r)\le 2r\;\textrm{Lip}\;f^{-n}_{\hat x}
\le 2r S_{\epsilon}(\hat x) e^{-n(L-\epsilon)} \le  2r N e^{-n(L-\epsilon)}$.\\
As moreover $f^{-n}_{\hat x}(D_r)$ meets $D_{r'}$,
there exists $n_0$, which depends only on $\epsilon$, $r$ and $r'$, such that $f^{-n}_{\hat x}(D_r) \subset D_{r}$ for every $\hat x\in\widehat C_n$ and $n\ge n_0$.
Thus, by Brouwer theorem, $f^{-n}_{\hat x}$ has a fixed point $p_n \in f^{-n}_{\hat x}(D_r)$ for every $\hat x\in\widehat C_n$ and $n\ge n_0$.
Since the elements of $S_n$ are mutually disjoint sets, we have produced $\big(\textrm{Card}\; S_n\big)$ fixed points of $f^n$ in $D_r$ for 
$n\ge n_0$. It remains to check that these fixed points  belong to $R^{2\epsilon}_n$. This actually follows immediately from the estimates on $\textrm{Lip}\;f^{-n}_{\hat x}$. Indeed:
\begin{center}
$\vert (f^{n})'(p_n)\vert =\vert (f^{-n}_{\hat x})' (p_n)\vert^{-1} \ge \big(\textrm{Lip}\;f^{-n}_{\hat x} \big)^{-1} \ge N^{-1} e^{n(L-\epsilon)} \ge e^{n(L-2\epsilon)}$
\end{center}
 for $n$ big enough.\\

Finally we prove the estimate \ref{Fact2}.
Let us first observe that
\begin{eqnarray}\label{observ}
 \pi\big(\hat f^{-n}(\widehat D_r\cap \widehat{X}_{reg\;N}^{\epsilon}\big) \cap \widehat D_{r'}\big)\subset\bigcup_{\hat x \in\widehat C_n} f^{-n}_{\hat x} (D_r).
\end{eqnarray}
This can be easily seen : if $\hat u \in \hat f^{-n}(\widehat D_r\cap \widehat{X}_{reg\;N}^{\epsilon})\cap \widehat D_{r'}$ then $u_0=\pi(\hat u) \in D_{r'}\cap  f^{-n}_{\hat x} (D_r)$ 
where $\hat x:=\hat f^n (\hat u) \in \widehat D_r\cap \widehat{X}_{reg\;N}^{\epsilon}.$\\
By the constant Jacobian property we have $\mu\big(f^{-n}_{\hat x} (D_r)\big) = d^{-n} \mu(D_r)$ and, since the sets $f^{-n}_{\hat x} (D_r)$ of the collection $S_n$
are mutually disjoint, we obtain
\begin{eqnarray}\label{equal}
\mu\big(\bigcup_{\hat x \in\widehat C_n} f^{-n}_{\hat x} (D_r)\big) =\big(\textrm{Card}\; S_n\big)\;d^{-n} \mu(D_r).
\end{eqnarray}
 Combining \ref{observ} with \ref{equal} yields \ref{Fact2}:
\begin{eqnarray*}
\big(\textrm{Card}\; S_n\big)\;d^{-n} \mu(D_r)\ge
\mu\big[ \pi\big(\hat f^{-n}(\widehat D_r\cap \widehat{X}_{reg\;N}^{\epsilon}\big) \cap \widehat D_{r'}\big) \big] =\\\hat\mu \big[\pi^{-1}\circ \pi
\big(\hat f^{-n}(\widehat D_r\cap \widehat{X}_{reg\;N}^{\epsilon}\big) \cap \widehat D_{r'}\big)\big]
\ge\hat\mu 
\big(\hat f^{-n}(\widehat D_r\cap \widehat{X}_{reg\;N}^{\epsilon}\big) \cap \widehat D_{r'}\big). 
\end{eqnarray*}
\qed

\chapter{Holomorphic families}

We introduce here the main spaces in which we shall work in the next chapters and present some of their structural properties.

\section{Generalities}
\subsection{Holomorphic families and the space $Rat_d$}

Let us start with a formal definition.

\begin{defi}Let $M$ be a complex manifold.
A holomorphic map 
\begin{center}
$f:M\times\pp\to\pp$
\end{center}
 such that all rational maps $f_\la:=f(\la,\cdot):\pp \to\pp$ have the same degree $d\ge 2$ is called \emph {holomorphic family} of degree $d$ rational maps parametrized by $M$. For short, any such family will be denoted $\left(f_{\la}\right)_{\la\in M}$.
\end{defi}

Any degree $d$ rational map $f:=\frac{a_d z^d+\cdot\cdot\cdot +a_1 z+a_0}{b_d z^d+\cdot\cdot\cdot +b_1 z+b_0}$ is totally defined by the point 
$[a_d:\cdot\cdot\cdot:a_0:b_d:\cdot\cdot\cdot:b_0]$ in the projective space ${\bf P}^{2d+1}$. This allows to identify the space 
$Rat_d$ of 
degree $d$ rational maps with a Zariski dense open subset of ${\bf P}^{2d+1}$.\\
We can be more precise by looking at the space of  homogeneous polynomial maps of $\CC^2$ which is identified to $\CC^{2d+2}$
by the correspondance 
\begin{center}$(a_d,\cdot\cdot\cdot,a_0,b_d,\cdot\cdot\cdot,b_0)\mapsto \left(\sum_{i=1}^d a_i z_1^i z_2^{d-i},\sum_{i=1}^d b_i z_1^i z_2^{d-i}\right)$.
\end{center}
Indeed, $Rat_d$ is precisely the image by the canonical projection $\pi:\CC^{2d+2}\to{\bf P}^{2d+1}$ of the subspace $H_d$ of $\CC^{2d+2}$ consisting of
non-degenerate polynomials. As $H_d$ is the complement in $\CC^{2d+2}$ of the projective variety defined by the vanishing of the resultant 
$Res(\left(\sum_{i=1}^d a_i z_1^i z_2^{d-i},\sum_{i=1}^d b_i z_1^i z_2^{d-i}\right)$, one sees that $Rat_d={\bf P}^{2d+1}\setminus \Sigma_d$
where $\Sigma_d$ is an (irreducible) algebraic hypersurface of ${\bf P}^{2d+1}$.\\

From now on, we will always consider $Rat_d$ as a quasi-projective manifold. 
We may therefore also see any holomorphic family of degree $d$ rational maps with parameter space $M$
as a holomorphic map $f$ from $M$ to $Rat_d$. In particular we may take for $M$ any submanifold of $Rat_d$; this is especially
interesting when $M$ is dynamically defined as are, for instance, the hypersurfaces $Per_n(w)$ which will be defined in the next section.\\

The simplest example of holomorphic family is the family of quadratic polynomials. Up to affine conjugation, any degree $2$ polynomial is of the form $z^2+a$.
To understand quadratic polynomials it is therefore sufficent to work with the family $\left(z^2+a\right)_{a\in {\tiny \CC}}$.\\

In most cases, when considering a holomorphic family $\left(f_{\la}\right)_{\la\in M}$,  we shall make the two following mild assumptions. 

\begin{ass}\label{Mass}
Let  $\left(f_{\la}\right)_{\la\in M}$ be any holomorphic family of degree $d$ rational maps.
\begin{itemize}
\item[A1] The \emph {marked critical points assumption} means that
the critical set $C_{\la}$ of $f_\la$ is given by $2d-2$ graphs: $C_{\la}=\cup_{1}^{2d-2}\{c_i(\la)\}$ where the maps $M\ni\la\mapsto c_i(\la)\in \pp$ are holomorphic.
\item[A2] The \emph {no persistent neutral cycles assumption} means that
if $f_{\la_0}$
has a neutral cycle then this cycle becomes attracting or repelling under a suitable small perturbation of $\la_0$.
\end{itemize}
\end{ass}

Using a ramified cover of the parameter space $M$ one may actually always make shure that the assumption $A2$ is satisfied.\\

The group of M\"obius transformations, which is isomorphic to 
$PSL(2,{\bf C})$, acts by conjugation on the space $Rat_d$ of degree $d$ rational maps.
The dynamical properties of two conjugated rational maps are clearly equivalent and it is therefore natural to
work within the quotient of $Rat_d$ resulting
from this action.\\

The moduli space $Mod_d$ is, by definition, the quotient of $Rat_d$ under the action of $PSL(2,{\bf C})$ by conjugation.
 We shall denote as follows the canonical
projection:

\begin{eqnarray*}
\Pi: Rat_d &\longrightarrow& Mod_d\\
f&\longmapsto&\bar f 
\end{eqnarray*}

We shall usually commit the abuse of language which consists
in considering an element of $Mod_d$ as a rational map.
For instance, "$\bar f$ has a $n$-cycle of multiplier $w$" means that
every element of $\bar f$
posseses such a cycle. We shall also sometimes write $f$ instead of $\bar f$.\\

Although the action of $PSL(2,{\bf C})$ is not free, it may be proven that 
$Mod_d$ is a normal quasi-projective variety \cite{Sil}.

\begin{rem}\label{helpwhen}
The following property is helpful for working in $Mod_d$.
Every element $f$ of $Rat_d$ belongs to a local submanifold $T_f$ whose dimension equals $2d-2$
and which is transversal to the orbit of $f$ under the action of $PSL(2,{\bf C})$. Moreover, $T_f$ is invariant
under the action of the stabilizer $Aut(f)$ of $f$ which is a finite subgroup of $PSL(2,{\bf C})$. Finally,
$\Pi\big(T_f\big)$ is a neighborhood of $\bar f$ in $Mod_d$ and $\Pi$ induces a biholomorphism
between $T_f/Aut(f)$ and $\Pi\big(T_f\big)$.
\end{rem}

\subsection{The space of degree $d$ polynomials}\label{ssPd} 

As for quadratic polynomials, there exists a nice parametrization of the space of degree $d$ polynomials.\\

Let ${\cal P}_d$ be the space of polynomials of degree $d\ge 2$ with $d-1$ marked critical points up to conjugacy by affine transformations. Although this space has a natural structure of affine variety of dimension $d-1$,
we may actually work with a specific parametrization of ${\cal P}_d$ which we shall now present.

For every $(c,a):=(c_1,c_2,\cdot\cdot\cdot,c_{d-2},a)\in{\bf C}^{d-1}$ we denote by $P_{c,a}$ the polynomial of degree $d$ whose critical points are $(0,c_1,\cdot\cdot\cdot,c_{d-2})$
and such that $P_{c,a}(0)=a^d$. This polynomial is explicitely given by:

$$P_{c,a}:=\frac{1}{d} z^d+\sum_{2}^{d-1}\frac{(-1)^{d-j}}{j} \sigma_{d-j} (c) z^j +a^d$$

where $\sigma_i(c)$ is the symmetric polynomial of degree $i$ in $(c_1,\cdot\cdot\cdot,c_{d-2})$. 
For convenience we shall set $c_0:=0$.\\

Thus, when considering degree $d$ polynomials, instead of working in  ${\cal P}_d$ we may consider the holomorphic family 
\begin{center}
$\big(P_{c,a}\big)_{(c,a)\in {\tiny {\bf C}^{d-1}}}$
\end{center}
whose parameter space $M$ is simply ${\bf C}^{d-1}$.
Using this parametrization, one may exhibit a finite ramified cover $\pi: {\bf C}^{d-1}\to {\cal P}_d$ (see \cite{DF} Proposition 5.1).\\

It will be crucial 
to consider the projective compactification ${\bf P}^{d-1}$ of  ${\bf C}^{d-1}=M$. This is why we wanted the expression of $P_{c,a}$ to be homogeneous in $(c,a)$ and have used the parameter $a^d$ instead of $a$. In this context, we shall denote by ${\bf P}_{\infty}$
 the projective space at infinity : ${\bf P}_{\infty}:=\{[c:a:0]\; ;(c,a)\in{\bf C}^{d-1}\setminus\{0\}\}$.

\subsection{The moduli space of degree two rational maps}\label{secMod}

In his paper \cite {Mi2}, Milnor has given a particularly nice description of  $Mod_2$
which we will now present. The reader may also consult the fourth chapter of
book of Silverman \cite{Sil}.\\

Any $f\in Rat_2$ has $3$ fixed points (counted with multiplicities) whose multipliers may be denoted 
$\mu_1,\mu_2,\mu_3$. Let us observe that $\mu_i=1$ if and only if one of the fixed point is not simple. The symmetric functions 
$$\sigma_1:=\mu_1+\mu_2+\mu_3,\;\;\;\sigma_2:=\mu_1\mu_2+\mu_1\mu_3+\mu_2\mu_3,
\;\;\;\sigma_3:=\mu_1\mu_2\mu_3$$
are clearly well defined on $Mod_2$. It follows from the holomorphic fixed point formula $\sum \frac{1}{1-\mu_i}=1$ (see \cite{Milnor4},  Lecture 12), applied in the generic case of three distinct fixed points,
that 

\begin{eqnarray}\label{index}
\sigma_3-\sigma_1 +2=0.
\end{eqnarray}

The above identity is the crucial point. It shows that the set of three multipliers $\{\mu_1,\mu_2,\mu_3\}$ is entirely determined by $(\sigma_1,\sigma_2)$.
More precisely,  $\{\mu_1,\mu_2,\mu_3\}$ is the set of roots of the polynomial $X^3 -\sigma_1 X^2 +\sigma_2 X -(\sigma_1-2)$.\\

Let us mention some further useful facts.

\begin{lem}\label{nextlem} The multipliers $\mu_1,\mu_2,\mu_3$ satisfy the following identities
\begin{eqnarray*}
(\mu_1-1)^2=(\mu_1\mu_2-1)(\mu_1\mu_3-1)\\
(\mu_2-1)^2=(\mu_2\mu_1-1)(\mu_2\mu_3 -1)\\
(\mu_3-1)^2=(\mu_3\mu_1-1)(\mu_3\mu_2 -1)
\end{eqnarray*}
In particular, if two fixed points are distinct and have multipliers $\mu_i$,$\mu_j$ then $\mu_i\mu_j\ne 1$ and the multiplier $\mu_k$ of the remaining fixed point is given by
$\mu_k= \frac{2-\mu_i-\mu_j}{1-\mu_i\mu_j}$.
\end{lem}
\proof The three  identities
are  deduced from $(X-1)^2-(XY-1)(XZ-1)=X(X+Y+Z-2-XYZ)$ taking $\mu_1\mu_2\mu_3-(\mu_1 +\mu_2 +\mu_3) +2=0$ (see \ref{index}) into account.\\

If two fixed points are distinct then one of them must be simple and thus $\mu_i\ne1$. Then $\mu_i\mu_j\ne 1$ follows from one of the previous identities.
The expression of $\mu_k$ is then obtained using  $\mu_1\mu_2\mu_3-(\mu_1 +\mu_2 +\mu_3) +2=0$.\qed\\

 Milnor has actually shown that $(\sigma_1,\sigma_2)$ induces a good 
parametrization of $Mod_2$ (\cite{Mi2}).

\begin{theo}\label{theoMiMod2}
The map $M:Mod_2 \to {\bf C}^2$ defined by $\bar f \mapsto (\sigma_1,\sigma_2)$
is one-to-one and onto.
\end{theo}  

\proof
Let us denote by $m$ the map which associates to any $\sigma:=(\sigma_1,\sigma_2)\in {\bf C}^2$ the set $\{\mu_1,\mu_2,\mu_3\}$ of roots of the polynomial $X^3 -\sigma_1 X^2 +\sigma_2 X -(\sigma_1-2)$. It is clear that $m(\sigma)\ne \{1,1,1\}$ if and only if $\sigma\ne (3,3)$.\\

Let $\sigma\in {\bf C}^2\setminus\{(3,3)\}$. Let us show that  $M(\bar f)=\sigma$ uniquely determines $\bar f\in Mod_2$. Set $\{\mu_1,\mu_2,\mu_3\}=m(\sigma)$.
At least one of the $\mu_i$ does not equal $1$ and thus the map $f$ has at least two distinct fixed points whose multipliers are, say, $\mu_1\ne 1$ and $\mu_2$. After conjugation we may assume that these fixed points are $0$ and $\infty$ and that 
\begin{center}
$f\sim z\frac{\alpha z + \beta}{\delta z +1}$ with $\alpha\ne 0$ and $\alpha-\beta\gamma\ne 0$.
\end{center}
After conjugating by $z\mapsto \alpha^{-1} z$ one gets 
\begin{center}
$f\sim z\frac{ z + \beta'}{\delta' z +1}$ with  $1-\beta'\gamma'\ne 0$
\end{center}
but, since $\beta'$ and $\delta'$ are obviously respectively equal to $\mu_1$ and $\mu_2$ one actually has $f\sim z\frac{ z + \mu_1}{\mu_2 z +1}$ with $\mu_1\mu_2\ne 1$. 
We have shown that $f\sim g\sim z\frac{ z + \mu_1}{\mu_2 z +1}$ as soon as $M(\bar f)=M(\bar g)$.\\

Let us now assume that $M(\bar f)=(3,3)$. In that case the multipliers of the fixed points of $f$ are all equal to $1$ and $f$ has actually only one fixed point. After conjugation, we may assume that this fixed point is $\infty$ and that $f^{-1}\{\infty\}=\{0,\infty\}$.
Then $f \sim\frac{p(z)}{z}$ but, since $f(z)-z=\frac{p(z)-z^2}{z}$ does not vanish on the complex plane, we must have $p(z)=z^2+c$ and $f \sim z+\frac{c}{z}$. After conjugating
by $z\mapsto \sqrt{c} z$ we get $f \sim z+\frac{1}{z}$.\\

We have shown that $M$ is one-to-one. It clearly follows from the above computations that $M$ is also onto.\qed\\

\begin{rem}\label{NorFor} The proof of the above Theorem shows that the expression $ z\frac{ z + \mu_i}{\mu_j z +1}$ can be used as a normal form where $\mu_i\mu_j\ne 1$.
\end{rem}

It will be extremely useful to consider the projective compactification of $Mod_2$ obtained through the above Theorem:

\begin{eqnarray*}
Mod_2 \ni \bar f \longmapsto [\sigma_1:\sigma_2:1] \in {\bf P}^2
\end{eqnarray*}
whose corresponding line at infinity
will be denoted by ${\cal L}$
\begin{eqnarray*}
{\cal L}:=\{[\sigma_1:\sigma_2:0];\;(\sigma_1,\sigma_2)\in{\bf C}^2\setminus\{0\}\}.
\end{eqnarray*}

It is important to stress that this compactification is actually natural in the sense 
that the "behaviour near ${\cal L}$'' captures a lot of dynamically meaningful 
informations. This will be discussed in subsection \ref{Onthegeom}. Let us now simply mention that the line at infinity $\cal L$ may somehow be parametrized by the limits of 
$\frac{\sigma_2}{\sigma_3}$.

\begin{prop}\label{paramL}
If $[\sigma_1:\sigma_2:1]$ converges to $P\in {\cal L}$ then the point $P$ is of the form $[1:\mu+\frac{1}{\mu}:0]$ where $\mu+\frac{1}{\mu}$ is a limit of 
$\frac{\sigma_2}{\sigma_3}=\frac{1}{\mu_1}+\frac{1}{\mu_2}+\frac{1}{\mu_3}$ in $\pp$. Moreover, $\mu\ne 0,\infty$ if and only if two of the multipliers 
$\mu_j$ stay bounded as $[\sigma_1:\sigma_2:1]$ approaches $P$.
\end{prop}

\proof Since  $(\sigma_1,\sigma_2)$ tends to $\infty$ then at least one of the $\mu_i$'s tends to $\infty$ too. Let us assume that $\mu_3\to \infty$.
We first look at the case where both $\mu_1$ and $\mu_2$ stay bounded. The third identity in Lemma \ref{nextlem} shows that $\mu_1\mu_2\to 1$. Thus $\sigma_3\to\infty$
and $\frac{\sigma_1}{\sigma_3}=\frac{\mu_1+\mu_2}{(\mu_1\mu_2)\mu_3}+ \frac{1}{\mu_1\mu_2} \to 1$. As 
$[\sigma_1:\sigma_2:1]=[\frac{\sigma_1}{\sigma_3}:\frac{\sigma_2}{\sigma_3}:\frac{1}{\sigma_3}]$ converges to $P\in {\cal L}$ one sees that $\frac{\sigma_2}{\sigma_3}$ must converge in $\pp$ and that $P$ is of the required form.\\
Let us now consider the case where $\mu_2 \to \infty$. The first identity of Lemma \ref{nextlem} implies that $\mu_1\to 0$. Thus $\sigma_2\sim \mu_2\mu_3$ and therefore 
$\frac{\sigma_2}{\sigma_3} \to \infty$ and $\frac{\sigma_1}{\sigma_2} \to 1$.  In that case $P=[0:1:0]$ which is again  the announced form.\qed\\ 

\section{The connectedness locus in polynomial families}\label{conloc}

\subsection{Connected and disconnected Julia sets of polynomials}\label{condiscon}

Among rational functions, polynomials are characterized by the fact that $\infty$ is a totally invariant
critical point. For any polynomial $P$ the super-attractive fixed point $\infty$ determines a basin of attraction
\begin{center}
${\cal B}_P(\infty):=\{z\in \CC\;\textrm{s.t.}\; \lim_n P^{n}(z)=\infty\}$.
\end{center}

This basin is always connected and its boundary is precisely the Julia set $J_P$ of $P$. The complement of 
${\cal B}_P(\infty)$ is called the filled-in Julia set of $P$.\\

As we already saw, another nice feature of polynomials is the possibility to define a Green function $g_{\tiny{\CC},P}$  
The Green function $g_{\tiny{\CC},P}$ is a subharmonic function on the complex plane which vanishes exactly on the filled-in Julia set of 
$P$ (see Proposition \ref{Greenpoly}).\\

Any degree $d$ polynomial $P$ is locally conjugated at infinity with the polynomial $z^{d}$. This means that there exists 
a local change of coordinates $\varphi_P$ (which is called B\"{o}ttcher function) such that $\varphi_P \circ P=(\varphi_P)^{d}$
on a neighbourhood of $\infty$. It is important to stress the following relation between the B\"{o}ttcher and Green functions:
\begin{center}
$\ln \vert \varphi_P \vert =g_{\tiny{\CC},P}$.
\end{center}

The only obstruction to the extension of the B\"{o}ttcher function $\varphi_P$ to the full basin ${\cal B}_P(\infty)$ is the presence of other critical points than $\infty$ in ${\cal B}_P(\infty)$.
This leads to the following important result:
\begin{theo}
For any polynomial $P$ of degree $d\ge 2$ the following conditions are equivalent:
\begin{itemize}
\item[i)] ${\cal B}_P(\infty)$ is simply connected
\item[ii)] $J_P$ is connected
\item[iii)] ${\cal C}_P \cap {\cal B}_P(\infty)=\{\infty\}$
\item[iv)] $P$ is conformally conjugated to $z^{d}$ on ${\cal B}_P(\infty)$.\\
\end{itemize}
\end{theo}

The above Theorem gives a nice characterization of polynomials having a connected Julia set. Let us apply it to the 
quadratic family $\left(z^2+a\right)_{a\in {\tiny \CC}}$. The Julia set $J_a$ of $P_a:=z^{2}+a$ is connected if and only if
the orbit of the critical point $0$ is bounded. In other words, the set of parameters $a$ for which 
$J_a$ is connected is the famous Mandelbrot set.
\begin{defi}{\bf Mandelbrot set.}\label{defimandel} Let $P_a$ denote the quadratic polynomial $z^{2}+a$. The Mandelbrot set ${\cal M}$ is defined by
\begin{center}
${\cal M}:=\{a\in \CC\;\textrm{s.t.}\; \sup_n\vert P_a^{n}(0)\vert <\infty\}$.
\end{center}
\end{defi}
The  Mandelbrot set is thus the connectedness locus of the quadratic family. It is not difficult to show that ${\cal M}$ is compact. The compacity of the connectedness locus in the polynomial families of degree $d\ge 3$ is a much more
delicate question which has been solved by Branner and Hubbard \cite{BH}. We shall treat it in the two next subsections
and also present a somewhat more precise result which will turn out to be very useful later.

\subsection{Polynomials with a bounded critical orbit}\label{ssPBCO}

We work here with the parametrization
$\big(P_{c,a}\big)_{(c,a)\in {\tiny {\bf C}^{d-1}}}$ of ${\cal P}_d$ and will use the projective compactification 
${\bf P}^{d-1}$ introduced in the subsection \ref{ssPd}.\\ 
We aim to show that  the subset of parameters $(c,a)$ for which the polynomial $P_{c,a}$ has at least one bounded critical orbit can only cluster on certain hypersurfaces of
${\bf P}_{\infty}$. The ideas here are essentially those used by Branner and Hubbard for proving the compactness of the connectedness locus (see \cite{BH} Chapter 1, section 3)
but we also borrow from the paper
(\cite{DF}) of Dujardin and Favre.\\

We shall use the following
\begin{defi} The notations are those introduced in subsection \ref{ssPd}.
For every $0\le i\le d-2$, the hypersurface $\Gamma_i$ of ${\bf P}_{\infty}$ is defined by: 
$$\Gamma_i:=\{[c:a:0]/\; \alpha_i(c,a)=0\}$$
where $\alpha_i$ is the homogeneous polynomial given by:
$$\alpha_i (c,a):=P_{c,a}(c_i)=\frac{1}{d}c_i^d+\sum_{j=2}^{d-1}\frac{(-1)^{d-j}}{j} \sigma_{d-j} (c) c_i^j +a^d.$$
We  denote by ${\cal B}_i$ the set of parameters $(c,a)$ for which the critical point $c_i$ of $P_{c,a}$ has a bounded forward orbit (recall that $c_0=0$):
\begin{center}
${\cal B}_i:=\{(c,a)\in {\bf C}^{d-1}\;\textrm{s.t.}\; \sup_n\vert P_{c,a}^n(c_i)\vert <\infty\}$.
\end{center}
\end{defi}

A crucial observation about the intersections of hypersurfaces $\Gamma_i$ is given by the next Lemma. 

\begin{lem}\label{leminter}
The intersection $\Gamma_0\cap\Gamma_1\cap\cdot\cdot\cdot\cap\Gamma_{d-2}$ is empty and
$\Gamma_{i_1}\cap\cdot\cdot\cdot\cap\Gamma_{i_k}$ has codimension $k$ in ${\bf P}_{\infty}$ if $0\le i_1<\cdot\cdot <i_k \le d-2$.
\end{lem}

\proof A simple degree argument shows that $P_{c,a}(0)=P_{c,a}(c_1)=\cdot\cdot\cdot=P_{c,a}(c_{d-2})=0$ implies that $c_1=\cdot\cdot\cdot=c_{d-2}=a=0$. Thus $\Gamma_0\cap\Gamma_1\cap\cdot\cdot\cdot\cap\Gamma_{d-2}=\emptyset$.
Then the conclusion follows from
Bezout's theorem.\qed\\

Since the connectedness locus coincides with $\cap_{0\le i\le d-2} {\cal B}_i$,
the announced result can be stated as follows.

\begin{theo}\label{controlinfty}
For every $0\le i\le d-2$, the cluster set of ${\cal B}_i$ in ${\bf P}_{\infty}$ is contained in $\Gamma_i$. In particular, the connectedness locus is compact in 
${\bf C}^{d-1}$.
\end{theo}

As the Green function $g_{c,a}$ of the polynomial $P_{c,a}$ is defined by
\begin{center} 
$g_{c,a}(z):=\lim_nd^{-n}\ln^+ \vert P_{c,a}^n (z)\vert $
\end{center}
one sees that
\begin{center}
${\cal B}_i=\{(c,a)\in {\bf C}^{d-1}\;\textrm{s.t.}\; g_{c,a}(c_i)=0\}$. 
\end{center}
This is why the proof of Theorem \ref{controlinfty} will rely on estimates on the Green functions 
and, more precisely, on the following result.

\begin{prop} \label{estimgreen} Let  
$g_{c,a} $ be the Green function of  $P_{c,a}$ and $G$ be 
 the function  defined on ${\bf C}^{d-1}$ by:
$G(c,a):=\max \{g_{c,a}(c_k);\;0\le k\le d-2\}$.
Let $\delta:= \frac{\sum_{k=0}^{d-2}c_k}{d-1}$. Then the following estimates occur:
\begin{itemize}
\item[1)] $G(c,a)\le\ln \max\{\vert a\vert,\vert c_k\vert\} +O(1)$ if $\max\{\vert a\vert,\vert c_k\vert\} \ge 1$
\item[2)] $\max \{g_{c,a}(z),G(c,a)\}\ge \ln \vert z-\delta\vert -\ln 4$.
\end{itemize}
\end{prop}

Let us first see how Theorem \ref{controlinfty} may be deduced from Proposition \ref{estimgreen}.\\

\noindent {\it Proof of Theorem \ref{controlinfty}.}  
Let $\Vert(c,a)\Vert_{\infty}:=\max\{\vert a\vert,\vert c_k\vert\}$. We simply have to check that  $\alpha_i\big(\frac{(c,a)}{\Vert (c,a)\Vert_{\infty}}\big)$
tends to $0$ when $g_{c,a}(c_i)$ stays equal to $0$ and $\Vert(c,a)\Vert_{\infty}$ tends to $+\infty$.
As $P_{c,a}(c_i)=\alpha_i(c,a)$ and $g_{c,a}(c_i)=0$, the estimates given by Proposition \ref{estimgreen} yield:
\begin{eqnarray*}
\ln\Vert (c,a)\Vert_{\infty} +O(1)\ge\max\big(dg_{c,a}(c_i),G(c,a)\big)=\max\{g_{c,a}\circ P_{c,a} (c_i),G(c,a)\}\ge \\
\ge\ln\frac{1}{4}\vert\alpha_i(c,a)-\delta\vert 
\end{eqnarray*}
 since $\alpha_i$ is $d$-homogeneous we then get:
$$(1-d)\ln \Vert (c,a)\Vert_{\infty} +O(1) \ge \ln\frac{1}{4}\vert\alpha_i\big(\frac{(c,a)}{\Vert (c,a)\Vert_{\infty}}\big)-\frac{\delta}{\Vert (c,a)\Vert_{\infty}^d}\vert$$
and the conclusion follows since $\frac{\delta}{\Vert (c,a)\Vert_{\infty}^d}$ tends to $0$ when $\Vert(c,a)\Vert_{\infty}$ tends to $+\infty$.\qed\\

Let us end this subsection by giving a\\
 
\noindent{\it Proof of Proposition \ref{estimgreen}}.
The first estimate is a standard consequence of the uniform growth of $P_{c,a}$ at infinity. Let us however prove it with care.
We will set $A:= \Vert (c,a)\Vert_{\infty}$ and $M_A(z):=\max\{A,\vert z\vert\}$. From
\begin{center}
$\vert P_{c,a}(z)\vert \le \frac{1}{d}\vert z\vert ^{d}\big(1+d\max\{\frac{\vert\sigma_{d-j}(c)\vert}{j\vert z\vert ^{d-j}},\frac{\vert a\vert^{d}}{\vert z\vert^{d}}\}\big)$
\end{center}
we get $\vert P_{c,a}(z)\vert \le C_d \vert z\vert^{d}$ for $\vert z\vert \ge A$ where the constant $C_d$ only depends on $d$. We may assume that $C_d\ge 1$.
By the maximum modulus principle this yields
\begin{center}
$\vert P_{c,a}(z)\vert \le C_d M_A(z)^{d}$.
\end{center}
It is easy to check that
\begin{center}
$M_A(Cz)\le CM_A(z)\;\;\textrm{if}\;C\ge 1$\\
$M_A\big(M_A(z)^{N}\big) =M_A(z)^{N}\;\;\textrm{if}\;A\ge 1$.
\end{center}
From now on we shall assume that $A\ge 1$.
By induction one gets 
\begin{center}
$\vert P_{c,a}^{n}(z)\vert \le C_d^{1+d+\cdot\cdot\cdot+d^{n-1}} M_A(z)^{d^{n}}$
\end{center}
which implies
\begin{center}
$g_{c,a}(z)\le \frac{\ln C_d}{d-1}+\ln \max\{\Vert (c,a)\Vert_{\infty},\vert z\vert\}\;\;\textrm{if}\; \Vert (c,a)\Vert_{\infty}\ge 1$
\end{center}
and in particular
\begin{center}
$G(c,a)\le \frac{\ln C_d}{d-1}+\ln\Vert (c,a)\Vert_{\infty}\;\;\textrm{if}\; \Vert (c,a)\Vert_{\infty}\ge 1$.
\end{center}

The second estimate is really more subtle. It exploits the fact that the Green function $g_{\tiny{\CC},P}$ coincides with the log-modulus of the B\"{o}ttcher coordinate function $\varphi_P$
and relies on a sharp control of the distorsions of this holomorphic function.\\

The B\"{o}ttcher coordinate function $\varphi_{c,a}:\{g_{c,a}>G(c,a)\}\to \CC$ is a univalent function  such that $\varphi_{c,a} \circ P_{c,a}=\varphi_{c,a}^{d}$. It is easy to check that $\ln \vert \varphi_{c,a}\vert =g_{c,a}$ where it makes sense and that $\varphi_{c,a}(z) = z-\delta +O(\frac{1}{z})$ where $\delta:= \frac{\sigma_1(c)}{d-1}=\frac{\sum c_k}{d-1}$.\\

One thus sees that $\varphi_{c,a}:\{g_{c,a}>G(c,a)\}\to \CC\setminus \overline{D}\big(0,e^{G(c,a)}\big)$ is a univalent map whose inverse $\psi_{c,a}$ satisfies
$\psi_{c,a}(z) = z+\delta +O(\frac{1}{z})$ at infinity.
We shall now apply the following result, which is a version of the Koebe $\frac{1}{4}$-theorem (see \cite{BH}, Corollary 3.3),
to $\psi_{c,a}$.
\begin{theo}
If $F:\widehat{\CC}\setminus \overline{D}_r \to \widehat {\CC}$ is holomorphic and injective and 
\begin{center}
$F(z)=z+\sum_{n=1}^{\infty} \frac{a_n}{z^{n}},\;z\in \CC\setminus \overline{D}_r$
\end{center}
then $  \CC\setminus \overline{D}_{2r} \subset F\big(\CC\setminus \overline{D}_r\big).$
\end{theo}
Pick $z\in \CC$ and set $r:=2\max\{e^{g_{c,a}(z)},e^{G(c,a)}\}$. Then $z\notin \psi_{c,a}\big(\CC\setminus \overline{D}_r\big)$ since otherwise
we would have $e^{g_{c,a}}(z)=\vert \varphi_{c,a} (z)\vert >r \ge 2 e^{g_{c,a}}(z)$.\\
Thus, according to the above distorsion theorem, $z\notin \CC\setminus \overline{D}\big(\delta,2r\big)$.
In other words $\vert z-\delta\vert \le 2r=4\max\{e^{g_{c,a}(z)},e^{G(c,a)}\}$ and the desired estimate follows by taking logarithms.\qed

\section{The hypersurfaces $\Per_n(w)$}\label{secPer}

We will consider here some dynamically
defined subsets of the parameter space which will play a central role in our study. 

\subsection{Defining the $\Per_n(w)$ using dynatomic polynomials}\label{ssdefper}

For any holomorphic family of rational maps,
the following result describes precisely the set of maps
having a cycle of given period and multiplier.

\begin{theo}\label{theopoly}
Let $f:M\times\pp\to\pp$ be a holomorphic family of degree $d\ge 2$ rational maps. Then for every integer $n\in\NN^*$ there exists a
holomorphic function $p_n$ on $M\times {\bf C} $ which is polynomial on ${\bf C}$ and such that:
\begin{itemize}
\item[1-] for any $w \in {\bf C} \setminus \{1\}$, the function
$p_n(\la,w)$ vanishes if and only if $f_{\la}$ has a cycle of exact period $n$ and multiplier $w$
\item[2-] $p_n(\la,1)=0$ if and only if $f_{\la}$ has a cycle of exact period $n$ and multiplier $1$ or a cycle of exact period $m$ whose multiplier is a primitive $r^{th}$ root of unity
with $r\ge 2$ and $n=mr$
\item[3-] for every $\la\in M$, the degree $N_d(n)$ of $p_n(\la,\cdot)$ satisfies $ d^{-n} N_d(n) \sim\frac{1}{n}$.
\end{itemize}
\end{theo}

This leads to the following
\begin{defi}
Under the assumptions and notations of Theorem \ref{theopoly}, one sets
\begin{center}
$\Per_n(w):=\{\la\in M/\;p_n(\la,w)=0\}$
\end{center}
for any integer $n$ and any complex number $w$.
\end{defi}

According to Theorem \ref{theopoly}, $Per_n(w)$ is (at least when $w\ne 1$) the set of parameters $\la$ for which $f_{\la}$ has a cycle of exact period $n$ and multiplier $w$. Moreover, $\Per_n(w)$ is an hypersurface in the parameter space $M$ or coincides with $M$.
We also stress that the estimate on the degree
$N_d(n)$ of $p_n(\la,\cdot)$ will be important in some of our applications. \\

We now start to explain the construction of the functions $p_n$. It clearly suffices to treat the case of the family $\Rat$ and then set $p_n(\la,w):=p_n(f_{\la},w)$
for any holomorphic family $M\ni\la\mapsto f_{\la}\in \Rat$. Our presentation is essentially based on the fourth chapter of the book
of Silverman \cite{Sil} and also borrows to the paper \cite{Mi2} of Milnor.\\

We will consider polynomial families; to deal with the general case one may adapt the proof by using lifts to $\CC^2$. 
According to the discussion
we had in subsection \ref{ssPd},  any degree $d$ polynomial $\varphi$ will be identified to a point in $\CC^{d-1}$.\\

The key point is to associate to any integer $n$ and any polynomial $\varphi$ of degree $d\ge 2$ a  polynomial $\Phi_{\varphi,n}^*$  whose roots, for a generic $\varphi$,  are precisely
the periodic points of $\varphi$  with exact period $n$.
 Such polynomials are called \emph {dynatomic} since they generalize cyclotomic ones, they are defined as follows.

\begin{defi}\label{defdyna}
For a degree $d$ polynomial $\varphi$ and an integer $n$ one sets
\begin{center}
$\Phi_{\varphi,n}(z):=\varphi^n(z)-z$.
\end{center}
 The associated \emph{dynatomic polynomials} are then defined by setting
\begin{center}
$\Phi_{\varphi,n}^*(z):=\prod_{k\vert n}\big(\Phi_{\varphi,k}(z)\big)^{\mu(\frac{n}{k})}$
\end{center}
where $\mu:\NN^*\to \{-1,0,1\}$ is the classical M\"{o}bius function.
\end{defi}

It is clear that $\Phi_{\varphi,n}$ is a polynomial whose roots are all periodic points of $\varphi$ with exact period \emph {dividing} $n$,
and that $\Phi_{\varphi,n}^*$ is a fraction whose roots and poles belong to the same set. Actually 
$\Phi_{\varphi,n}^*$ is
still a polynomial but this is not at all obvious! To prove it, one will systematically exploits the
fact that the sum $\sum_{k\vert n}\mu(\frac{k}{n})$ vanishes if $n>1$ and computes the valuation of  $\Phi_{\varphi,n}^*$ at
any $m$-periodic point of $\varphi$ for $m\vert n$. We then obtain a precise description of the roots of $\Phi_{\varphi,n}^*$:

\begin{theo}\label{theodyna}
Let $\varphi$ be a polynomial of degree $d\ge 2$. Then $\Phi_{\varphi,n}^*$ is a polynomial whose roots are the periodic points of $\varphi$ with exact period $m$
dividing $n$ and multiplier $w$ satisfying $w^r=1$ when $2\le r:=\frac{n}{m}$. The degree $\nu_d(n)$ of $\Phi_{\varphi,n}^*$ is equivalent to $d^{n}$.\\
\end{theo}

The proof of the above Theorem will be given in the next subsection, for the moment we admit it and prove Theorem \ref{theopoly}.
Let us note that the degree $\nu_d(n)$ of $\Phi_{\varphi,n}^*$ is given by $\nu_d(n)=\sum_{k\vert n} d^k \mu(\frac{n}{k})$ and is clearly equivalent to $d^{n}$ since 
$\vert \mu\vert \in \{0,1\}$ and $\mu(1)=1$. We may also observe that $\left(\varphi^n\right)'(z)=1$ for any root $z$ of $\Phi_{\varphi,n}^*$  whose period strictly divides $n$.\\

The construction of
 $p_n(\varphi,w)$ requires to understand the structure of the zero set $\textrm{Per}_n$  of  $(\varphi,z)\mapsto\Phi_{\varphi,n}^*(z)$. 
We recall that $\varphi$ is seen as a point in $\CC^{d-1}$.
Here are the informations we need.

\begin{prop}\label{Per_n}
The set $\textrm{Per}_n:=\{\left(\varphi,z\right)\;/\; \Phi_{\varphi,n}^*(z)=0\}$ is an algebraic subset of $\CC^{d-1}\times \CC$.
The roots of $\Phi_{\varphi,n}^*$ are simple and have exact period $n$
when $\varphi\in \CC^{d-1}\setminus X_n$ for some proper algebraic subset $X_n$ of $\CC^{d-1}$.
\end{prop}

\proof One sees on definition \ref{defdyna} that $\Phi_{\varphi,n}^*(z)$ is rational in $\varphi$. On the other hand, 
 $\Phi_{\varphi,n}^*(z)$ is locally bounded as it follows from the description of its roots given by Theorem \ref{theodyna}. Thus $\Phi_{\varphi,n}^*(z)$ is actually polynomial in $\varphi$.\\
Let us set $\Delta(\varphi):=\prod_{i\ne j} \left(\alpha_i(\varphi)-\alpha_j(\varphi)\right)$ where the $\alpha_i$ are the roots of $\Phi_{\varphi,n}^*$ counted with multiplicity.
This is a well defined function which vanishes exactly when $\Phi_{\varphi,n}^*$ has a multiple root. This function is holomorphic outside its zero set
and therefore everywhere by Rado's theorem.
Then $\{\Delta=0\}$ is an analytic subset of $\CC^{d-1}$ which is proper since $\varphi_0:=z^d \notin \{\Delta=0\}$.\\
Let $Y_n$ be the projection of $Per_n\cap\{(\varphi^n)'(z)=1\}$ onto $\CC^{d-1}$. By Remmert mapping theorem $Y_n$ is an analytic subset of $\CC^{d-1}$. Using $\varphi_0$ again one sees
that $Y_n\ne \CC^{d-1}$. Since $\left(\varphi^n\right)'(z)=1$ when $z$ belongs to a cycle of $\varphi$  whose period strictly divides $n$, one may take $X_n:=\{\Delta=0\}\cup Y_n$.\qed\\

Let $\textrm{Z}(\Phi_{\varphi,n}^*)$ be the set of roots of $\Phi_{\varphi,n}^*$ taken with multiplicity.
If $z\in \textrm{Z}(\Phi_{\varphi,n}^*)$ has exact period $m$ with $n=mr$,  we  denote by $w_n(z)$ the $r$-th power of the multiplier of $z$ (that is 
$(\varphi^n)'(z)$).
As  Theorem \ref{theodyna} tells us :
\begin{center}
a point $z$ is periodic of exact period $n$ and $w_n(z)\ne 1$ if and only if\\
$z\in \textrm{Z}(\Phi_{\varphi,n}^*)$ and $w_n(z)\ne 1$.
\end{center}

Let us now consider the sets
\begin{center}
$\Lambda_n^*(\varphi):=\{w_n(z);\;z\in \textrm{Z}(\Phi_{\varphi,n}^*)\}$
\end{center}
 and let us denote by $\sigma_i ^{*(n)}(\varphi)$, $1\le i\le \nu_d(n)$,  the associated symmetric functions.
The symmetric functions  $\sigma_i ^{*(n)}$ are globally defined and continuous on $\CC^{d-1}$ and, according to Proposition \ref{Per_n}, are holomorphic outside $X_n$. 
These functions are therefore holomorphic on $\CC^{d-1}$.
We set 
\begin{center}
$q_n(\varphi,w):=\prod_{i=0}^ {\nu_d(n)} \sigma_i ^{*(n)}(\varphi) (-w)^ {\nu_d(n)-i}$.
\end{center}
By construction $q_n(\varphi,w)=0$ if and only if $w\in\Lambda_n^*(\varphi)$ and $q_n$ is holomorphic in $(\varphi,w)$ and polynomial in $w$. 
As Proposition \ref{Per_n} shows, the elements of $\textrm{Z}(\Phi_{\varphi,n}^*)$ are cycles of exact period $n$ and therefore
each element of $\Lambda_n^*(\varphi)$ is repeated $n$ times when $\varphi\notin X_n$. This means that there exists a polynomial $p_n(\varphi, \cdot)$
such that $q_n(\varphi, \cdot)=\left(p_n(\varphi, \cdot)\right)^n$ when $\varphi\notin X_n$. 
As $p_n(\varphi,w)$ is holomorphic where it does not vanish, one sees that $p_n$ extends holomorphically to all $\CC^{d-1}\times \CC$. In other words,
$p_n$ may be defined by 
\begin{center}
$\left(p_n(\varphi,w)\right)^n:=q_n(\varphi,w)=\prod_{i=0}^ {\nu_d(n)} \sigma_i ^{*(n)}(\varphi) (-w)^ {\nu_d(n)-i}$.
\end{center}

The degree $N_d(n)$ of $p_n(\la,\cdot)$ is equal to $\frac{1}{n}\nu_d(n)=\frac{1}{n}\sum_{k\vert n}\mu(\frac{n}{k})d^k$. In particular $  N_d(n) \sim \frac{d^{n}}{n}$.
\qed
\\

\subsection{The construction of dynatomic polynomials}

We aim here to  prove Theorem \ref{theodyna}. For this purpose, let us recall that the M\"{o}bius function $\mu:\NN^*\to \{-1,0,1\}$ enjoys the following fundamental property:
\begin{eqnarray}\label{fondamu}
\sum_{k\vert n}\mu\left(\frac{n}{k}\right)=0\;\textrm{for any}\;n\in \NN^*.
\end{eqnarray}

Let us also adopt a few more notations. The valuation of $\Phi_{\varphi,n} \left(\;\textrm{resp.}\;\Phi_{\varphi,n}^*\right)$ at some point $z$ will be denoted
$a_z(\varphi,n) \left(\;\textrm{resp.}\;a_z^*(\varphi,n)\right)$. The set of $m$-periodic points of $\varphi$ will be denoted $Per(\varphi,m)$.\\

The following Lemma summarizes elementary facts.

\begin{lem}\label{lemval}
Let $\psi$ be a polynomial and $z\in Per(\psi,1)$. Let $\lambda:=\psi'(z)$, then for $q\ge 2$ one has:
\begin{itemize}
\item[i)] $\lambda^q\ne 1\Rightarrow a_z(\psi,q)=a_z(\psi,1)=1$
\item[ii)] $\lambda \ne 1\;\textrm{and}\;\lambda^q = 1\Rightarrow a_z(\psi,q)>a_z(\psi,1)=1$
\item[iii)] $\lambda= 1\Rightarrow a_z(\psi,q)=a_z(\psi,1)\ge 2$.
\end{itemize} 
\end{lem}

\proof  We may assume $z=0$ and set $\psi=\lambda X+\alpha X^e +o(X^e)$. Then $\Phi_{\psi,q}$ is equal to  $\left(\lambda^q -1\right) X+o(X)$
  and to $q\alpha X^e +o(X^e)$ if $\lambda=1$. The assertions $i)$ to $iii)$ then follow immediately.\qed\\

We have to compute $a_z^*(\varphi,n)$ for $z\in Per_n(\varphi,m)$ and $m\vert n$. We denote by $\lambda$ the multiplier of 
$z$ (i.e. $\la=\left(\varphi^m\right)'(z)$) and set $N:=\frac{n}{m}$. When $\lambda$ is a root of unity we denote by $r$ its order. Clearly, 
 Theorem \ref{theodyna} will be proved if we establish the following three facts:

\begin{itemize}
\item[F1 :] $N=1\Rightarrow  a_z^*(\varphi,n)>0$
\item[F2 :] $N\ge 2\;\textrm{and}\; \la^N\ne 1\;\textrm{or}\; \la=1\Rightarrow a_z^*(\varphi,n)=0$
\item[F3 :] $N\ge 2,\;\la^N= 1\;\textrm{and}\; \la\ne1\Rightarrow a_z^*(\varphi,n)\ge 0\;\textrm{and}
\;a_z^*(\varphi,n)>0\;\textrm{iff}\; r=N.$
\end{itemize}

 Besides definitions, the following computation is based on the obvious facts that $a_z(\varphi,k)=0$ when $k\ne qm$ for some $q\in \NN^*$
 and that $\varphi^{qm}=(\varphi^{q})^{m}$:
\begin{eqnarray*}\label{fst}
a_z^*(\varphi,n)=\sum_{k\vert n} \mu\left(\frac{n}{k}\right)a_z(\varphi,k)=\sum_{qm\vert n} \mu\left(\frac{n}{qm}\right)a_z(\varphi,qm)=\\
\sum_{q\vert N} \mu\left(\frac{N}{q}\right)a_z(\varphi,qm)=
\sum_{q\vert N} \mu\left(\frac{N}{q}\right)a_z(\varphi^m,q).
%=a_z^*(\varphi^m,N).
\end{eqnarray*}

We now proceed Fact by Fact, always starting with the above identity.\\

F1) $a_z^*(\varphi,n)=\mu(1)a_z(\varphi^n,1)>0$.\\

F2) Since $\la^q\ne 1$ when $q\vert N$ although $\la=1$, the assertions i) and iii) of Lemma \ref{lemval} tell us that $a_z(\varphi^m,q)=a_z(\varphi^m,1)$.
Then $a_z^*(\varphi,n)=
\left(\sum_{q\vert N} \mu\left(\frac{N}{q}\right)\right)a_z(\varphi^m,1)$ which, according to \ref{fondamu}, equals $0$.\\

F3) When $r$ is not dividing $q$ then $\la^q\ne 1$ and, by the assertion i) of Lemma \ref{lemval} we have $a_z(\varphi^m,q)=a_z(\varphi^m,1)$.
We may therefore write
\begin{eqnarray*}
a_z^*(\varphi,n)=
\left(\sum_{q\vert N} \mu\left(\frac{N}{q}\right)\right)a_z(\varphi^m,1) +
\sum_{q\vert N,\;r\vert q} \mu\left(\frac{N}{q}\right)\left[ a_z(\varphi^m,q) -a_z(\varphi^m,1) \right].
\end{eqnarray*}
By \ref{fondamu}, the first term in the above expression vanishes and we get
\begin{eqnarray*}
a_z^*(\varphi,n)=
\sum_{k\vert \frac{N}{r}} \mu\left(\frac{N/r}{k}\right)\left[ a_z(\varphi^m,rk) -a_z(\varphi^m,1) \right]=\\
\sum_{k\vert \frac{N}{r}} \mu\left(\frac{N/r}{k}\right)\left[ a_z(\varphi^{mr},k) -a_z(\varphi^m,1) \right]=
a_z^*(\varphi^{mr},\frac{N}{r}) - a_z(\varphi^m,1)\sum_{k\vert \frac{N}{r}} \mu\left(\frac{N/r}{k}\right).
\end{eqnarray*}

We finally consider two subcases.\\

If $N\ne r$ then, by \ref{fondamu}, $\sum_{k\vert \frac{N}{r}} \mu\left(\frac{N/r}{k}\right)=0$ and thus
$a_z^*(\varphi,n)=a_z^*(\varphi^{mr},\frac{N}{r})$. Since $z$ is a fixed point of $\varphi^{mr}$ whose multiplier
equals $\la^r=1$, Fact F2 shows that $a_z^*(\varphi^{mr},\frac{N}{r})=0$.\\

If $N=r$ we get 
\begin{center}
$a_z^*(\varphi,n)=a_z^*(\varphi^{mr},1) - a_z(\varphi^m,1)=a_z(\varphi^{mr},1) - a_z(\varphi^m,1)
=a_z(\varphi^{m},r) - a_z(\varphi^m,1).$
\end{center} Since $z$ is a fixed point of $\varphi^m$ whose multiplier $\la$ satisfies $\la^r=\la^N=1$ and $\la\ne 1$,
the assertion ii) of Lemma \ref{lemval} shows that this quantity is strictly positive.\qed\\

\begin{rem}\label{RemAlg}
It follows easily from the above construction that the growth of $p_n(\la,w)$ is polynomial in $\la$
when $\la\in {\CC}^{d-1}$. This shows that $p_n(\la,w)$ is actually a polynomial function on ${\CC}^{d-1}\times \CC$.
\end{rem}

\subsection{On the geometry of the $\Per_n(w)$  in particular families}\label{Onthegeom}

In this subsection we will describe some geometric properties of the hypersurfaces $\Per_n(w)$ in specific families of rational maps. We are mainly interested in the behaviour at infinity in the projective compactifications of the polynomial families and the moduli space of degree two rational maps.\\

We start with the polynomial family of degree $d$ polynomials (see subsection \ref{ssPd}). 
The following result is a consequence of our investigations of subsection \ref{ssPBCO}.
We recall that, according to Remark \ref{RemAlg}, the sets $\Per_m(\eta)$ may be seen as 
algebraic subsets of the projective space 
${\bf P}^{d-1}$.

\begin{prop}\label{interPer}
If $1\le k\le d-1$, $m_1<m_2<\cdot\cdot\cdot< m_k$ and $\sup_{1\le i\le k}\vert \eta_i\vert <1$ then 
$\Per_{m_1}(\eta_1)\cap\cdot\cdot\cdot\cap \Per_{m_k}(\eta_k)$ is an algebraic subset of codimension $k$ whose intersection with ${\bf C}^{d-1}$ is not empty. 
\end{prop}

\proof
By Bezout's theorem,  $\Per_{m_1}(\eta_1)\cap\cdot\cdot\cdot\cap \Per_{m_k}(\eta_k)$ is  a non-empty algebraic subset of
${\bf P}^{d-1}$ whose dimension is bigger than $(d-1-k)$.\\
Any cycle of attracting basins capture a critical orbit. Therefore, Theorem \ref{controlinfty} implies that the intersection of 
${\bf P}_{\infty}$ with $\Per_{m_1}(\eta_1)\cap\cdot\cdot\cdot\cap \Per_{m_k}(\eta_k)$ is contained
in some $\Gamma_{i_1}\cap\cdot\cdot\cdot\cap\Gamma_{i_k}$ since the $m_i$ are mutually distinct and the $\vert\eta_i\vert$ strictly smaller than $1$. Then, according to Lemma \ref{leminter}, ${\bf P}_{\infty}\cap\Per_{m_1}(\eta_1)\cap\cdot\cdot\cdot\cap \Per_{m_k}(\eta_k)$ has codimension $k$ in ${\bf P}_{\infty}$. The conclusion now follows from obvious dimension considerations. \qed\\

We now consider the space $Mod_2$.  As it has been discussed in subsection \ref{secMod}, this space can be identified to ${\bf C}^2$ and has a natural projective compactification which is given by 
$$Mod_2 \ni \bar f \longmapsto [\sigma_1:\sigma_2:1]\in {\bf P}^2$$
where $\sigma_1$,$\sigma_2$,$\sigma_3$ are the symetric functions of the  three multipliers $\{\mu_1,\mu_2,\mu_3\}$ of the fixed points 
 (we recall that $\sigma_3=\sigma_1-2$, see \ref{index}).\\
The line at infinity
${\cal L}=\{[\sigma_1:\sigma_2:0];\;(\sigma_1,\sigma_2)\in{\bf C}^2\setminus\{0\}\}$ enjoys an interesting dynamical parametrization (see Proposition \ref{paramL}).\\

 Using this identification, the defining functions  $p_n(\la,w)$ of $Per_n(w)$ are polynomials on 
${\bf C}^2\times {\bf C}$. 
Any $Per_n(w)$ may be seen as a curve in ${\bf P}^2$.
Understanding the behaviour of these curves at infinity is crucial to investigate the structure of the bifurcation locus in $Mod_2$. The following facts have been proved by Milnor 
(\cite{Mi2}). 

\begin{prop}\label{Per/L}
\begin{itemize}
\item[1)] For all $w\in{\bf C}$ the curve $Per_1(w)$ is actually a line  whose equation in ${\bf C}^2$ is
$(w^2+1)\la_1-w\la_2-(w^3+2)=0$ and whose
point at infinity is $[w:w^2+1:0]$. 
In particular, $Per_1(0)=\{\la_1=2\}$ is the line of 
quadratic polynomials, its point at infinity is $[0:1:0]$.
\item[2)] For $n>1$ and $w\in {\bf C}$ the points at infinity of the curves
$Per_n(w)$ are of the form $[u:u^2+1:0]$ with $u^q=1$ and $q\le n$.
\end{itemize}
\end{prop}

\proof 
1) This is a straightforward computation using the fact that
 the multipliers of the fixed points $\{\mu_1,\mu_2,\mu_3\}$ are the roots of the polynomial $X^3 -\sigma_1 X^2 +\sigma_2 X -(\sigma_1-2)$.\\
 
 2) Let $P\in Per_n(w)\cap {\cal L}$.  Let $\bar\sigma:=[\sigma_1:\sigma_2:1]\in Per_n(w)$ such that $\bar\sigma\to P$.
 We denote by $\mu_j$ the multipliers of the fixed points of $\bar\sigma$.
According to Proposition \ref{paramL} we may write $P=[1:\mu+\frac{1}{\mu}:0]$ where $\mu$ is a limit of 
 $\frac{1}{\mu_1}+\frac{1}{\mu_2}+\frac{1}{\mu_3}$. \\
 
 We first consider the case where $\mu\ne 0,\infty$. Then we may assume that $\mu_3\to \infty$, $\mu_2\to \mu$ and $\mu_1\to\frac{1}{\mu}$
 (recall that $\mu_1\mu_2$ must tend to $1$ as it follows from the identities of Lemma \ref{nextlem}). 
 By Remark \ref{NorFor} we may use the normal form  
 \begin{center}
 $f=f_{\bar\sigma}:= z\frac{ z + \mu_i}{\mu_j z +1}$.
 \end{center}
 Let us set $\delta:=1-\mu_1\mu_2$ and $l(z):=\mu_2 z+1$. An easy computation yields
 \begin{eqnarray}\label{f/z}
\mu_2
 \frac{f(z)}{z}=
 1-\frac{\delta}{l(z)}
 \end{eqnarray}
 and 
  \begin{eqnarray}\label{f'}
\mu_2f'(z)=1-\frac{\delta}{l(z)^2}.
 \end{eqnarray}
For $\vert \delta\vert <1$, we set:
 \begin{center}
 $D:=\{\vert l\vert <\vert \delta\vert^{\frac{2}{3}}\},\; A:=\{\vert \delta\vert^{\frac{2}{3}}\le\vert l\vert <\vert \delta\vert^{\frac{1}{3}}\},\;
 C:=\{\vert l\vert \ge\vert \delta\vert^{\frac{1}{3}}\}$.
 \end{center}
We thus have a decomposition of the Riemann sphere into three distinct regions:$$\pp=D\cup A\cup C.$$ 
 Observe that the disc $D\cup A$ degenerates to the point $\{\frac{-1}{\mu}\}$ as $\bar\sigma\to P$.\\
 
 Using the identities \ref{f/z} and \ref{f'} one gets the following estimates
 
 \begin{eqnarray}\label{D}
 \frac{f(z)}{z}=\frac{1}{\mu_2} +O\left(\vert \delta\vert^{\frac{1}{3}}\right)\;\textrm{on}\; A\cup C
 \end{eqnarray}
 \begin{eqnarray}\label{AUC}
 \vert f'(z)\vert \ge \vert \delta\vert^{\frac{-1}{3}}\;\textrm{on}\; D
 \end{eqnarray}
  \begin{eqnarray}\label{C}
f'(z)=\frac{1}{\mu_2} +O\left(\vert \delta\vert^{\frac{1}{3}}\right)\;\textrm{on}\;  C.
 \end{eqnarray}
 
 Let $O:=\{z_0,f(z_0),\cdot\cdot\cdot,f^{n-1}(z_0)\}$ be the orbit of some $n$-periodic point $z_o$ of $f=f_{\bar\sigma}$ whose mutiplier equals $w$
 (recall  that $\bar\sigma\in Per_n(w)$). We assume here that $\bar\sigma$ is very close to $P$. Then it is impossible that $O\subset D\cup C$ with $O\cap D\ne\emptyset$
 since otherwise, by the estimates \ref{D} and \ref{C}, we would have 
 $\vert w\vert =\prod_0^{n-1}\vert f'\left(f^j(z_0)\right)\vert \ge C \vert \delta\vert^{-\frac{1}{3}}$. Thus, either $O\subset A\cup C$ or $O\cap D\ne\emptyset$ and $O\cap A\ne\emptyset$.\\
 
 If $O\subset A\cup C$ then, using the estimate \ref{AUC}, we get
 \begin{eqnarray*}
1= \frac{f^n(z_0)}{z_0}=\prod_0^{n-1} \frac{f^{j+1}(z_0)}{f^j(z_0)}= \frac{1}{\mu_2^n}\left(1+O\left(\vert \delta\vert^{\frac{1}{3}}\right)\right)^n
 \end{eqnarray*}
 which implies that $\mu^n=1$ since $\mu_2\to \mu$ and $\delta\to 0$ when $\bar\sigma\to P$.\\
 
 If $O\cap D\ne\emptyset$ and $0\cap A\ne\emptyset$, we may assume that $z_0\in A$, $f^q(z_0)\in D$ and $\{z_0,\cdot\cdot\cdot,f^{q-1}(z_0)\}\subset A\cup C$.
 Then
  \begin{eqnarray*}
\frac{f^q(z_0)}{z_0}=\prod_0^{q-1} \frac{f^{j+1}(z_0)}{f^j(z_0)}= \frac{1}{\mu_2^q}\left(1+O\left(\vert \delta\vert^{\frac{1}{3}}\right)\right)^q
 \end{eqnarray*}
 
 which implies that $\mu^q=1$ since $\mu_2\to \mu$ and both $z_0$ and $f^q(z_0)$ belong to $A\cup D$ which tends to  $\{\frac{-1}{\mu}\}$ when $\bar\sigma\to P$.\\

 When $\mu=0$ or $\mu=\infty$ the proof is similar but slighly more subtle. We refer the reader to the paper \cite{Mi2} for details.\qed\\

The following Proposition, also due to Milnor (see Theorem 4.2 in \cite{Mi2}), implies that the curves $Per_n(w)=\{p_n(\cdot,w)=0\}$ have no multiplicity.

\begin{prop}\label{multi}
Let $N_2(n):=Card\left(Per_n(0)\cap Per_1(0)\right)$ be the number of hyperbolic components of period $n$ in the Mandelbrot set. Then $N_2(n)=\frac{\nu_2(n)}{2}$
where $\nu_2(n)$ is defined inductively by $\nu_2(1)=2$ and $2^n=\sum_{k\vert n}\nu_2(k)$.
Moreover, for any $w\in \Delta$ and any $\eta\in\Delta$ we have $Deg\;p_n(\cdot,w)=N_2(n)=Card\left(Per_n(w)\cap Per_1(\eta)\right).$
\end{prop}

 Epstein \cite{epstein} has obtained some far advanced generalizations of the above result and, in particular, has proved the boundedness of certain hyperbolic components of $Mod_2$:

\begin{theo}\label{thEps}
Let $H$ be a hyperbolic component of $Mod_2$ whose elements admit two distinct attracting cycles. 
If neither attractor is a fixed point then $H$ is relatively compact in $Mod_2$
\end{theo}

\chapter{The bifurcation current}\label{chapBC}

In this chapter, we consider an arbitrary holomorphic family $\left(f_{\la}\right)_{\la\in M}$ of degree $d$ rational maps  with marked critical points 
 (see \ref{Mass}).
Our first aim is to describe necessary and sufficient conditions for the Julia set $\Julla$ of $f_\la$ to move  holomorphically with the parameter $\la$. The parameters around which such a motion does exist are called stable. We will show that the set of stable parameters is dense in the parameter space, to this purpose we will relate the stability of
$\Julla$ with the stability of the dynamics on the  critical set $\Crla$ of $f_\la$; this is the essence of the Ma\~{n}\'e-Sad-Sullivan theory.
We will then exhibit a closed positive $(1,1)$-current on the parameter space whose support is precisely the complement of the stability locus; this is the bifurcation current.

\section{Stability versus bifurcation}
\subsection{Motion of repelling cycles and  Julia sets}

As Julia sets coincide with the closure of the sets of repelling cycles (see Theorems \ref{theoFJ} and \ref{theoLyu}), it is natural to investigate how $\Julla$ varies with $\la$
through the parametrizations of such cycles.\\
Assume that $f_{\la_0}$ has a repelling $n$-cycle $\{z_0,f_{\la_0}(z_0),\cdot\cdot\cdot,f_{\la_0}^{n-1}(z_0)\}$. Then, by the implicit function theorem applied to the equation 
$f_{\la}^n(z)-z=0$ at $(\la_0,z_0)$, there exists a neighbourhood $U_{0}$ of $\la_0$ in $M$ and a holomorphic map 
$$U_0\ni \la\mapsto h_{\la}(z_0)\in  \pp$$ such that $h_{\la_0}(z_0)=z_0$ and $h_{\la}(z_0)$ is a $n$-periodic
repelling point of $f_\la$ for all $\la\in U_0$. Moreover,  $h_\la(\cdot)$ can be extended to the full cycle so that $$f_\la\circ h_\la=h_\la\circ f_{\la_0}.$$
This lead us to say  that every repelling $n$-cycle of $f_{\la_0}$ \emph{moves holomorphically} on some neighbourhood of $\la_0$ and to set
 the following formal definition.

\begin{defi} 
 Let us denote by ${\cal R}_{\la,n}$ the set of repelling $n$-cycles of $f_{\la}$.
Let $\Omega$ be a neighbourhood of  $\la_0$ in $M$. One says that 
${\cal R}_{\la_0,n}$ moves holomorphically on $\Omega$ if there exists a map 
\begin{center}
$h: \Omega\times {\cal R}_{\la_0,n}\ni(\la,z)\mapsto h_{\la}(z) \in  {\cal R}_{\la_,n}$
\end{center}
which depends holomorphically on $\la$ and satisfies $h_{\la_0}=Id$,
$f_\la\circ h_\la=h_\la\circ f_{\la_0}$.
\end{defi}
 
More generally, the holomorphic motion of an arbitrary subset of the Riemann sphere is defined in the following way.

\begin{defi} 
Let $E$ be subset of the Riemann sphere and $\Omega$ be a complex manifold. Let $\la_0\in\Omega$. An holomorphic motion of $E$ over $\Omega$ and
centered at $\la_0$ is a map
\begin{center}
$h: \Omega\times E \ni(\la,z)\mapsto h_{\la}(z) \in  {\widehat \CC}$
\end{center}
which satisfies the following properties:
\begin{itemize}
\item[i)] $h_{\la_0}=Id\vert_{E}$
\item[ii)] $E\ni z\mapsto h_{\la}(z)$ is one-to-one for every $\la\in \Omega$
\item[iii)] $\Omega \ni \la \mapsto h_{\la}(z)$ is holomorphic for every $z\in E$. 
\end{itemize} 
\end{defi}

The interest of holomorphic motions relies on the fact that any holomorphic motion of a set $E$ extends to the closure of $E$. 
This is a quite simple consequence of Picard-Montel theorem.

\begin{lem}{\bf (basic $\la$-lemma)}\label{lamlem}
Let $E\subset {\widehat \CC}$ be a subset of the Riemann sphere and
$\sigma : E\times \Omega \ni(z,\la)\mapsto \sigma(z,\la)\in {\widehat \CC}$ be a holomorphic motion of $E$ over $\Omega$.
Then $\sigma$ extends to a holomorphic motion $\tilde \sigma$ of $\overline{E}$ over $\Omega$. Moreover
$\tilde\sigma$ is continuous on $\overline{E}\times \Omega$.
\end{lem}

 As $\Julla$ is the closure of the set of repelling cycles of $f_\la$,
this Lemma implies that the Julia set $\Jullo$ moves holomorphically over a neighbourhood $V_{\la_0}$ of $\la_0$ in $M$ as soon as \emph{all} repelling cycles of $f_{\la_0}$
move holomorphically on $V_{\la_0}$. Moreover, the holomorphic motion obtained in this way clearly conjugates the dynamics: $h_{\la} (\Jullo)=\Julla$ and
$f_\la\circ h_\la=h_\la\circ f_{\la_0}$ on $\Jullo$.\\

Our observations may now be gathered in the following basic Lemma.

\begin{lem}\label{lemmotcyc}
If there exists a neighbourhood $\Omega$ of $\la_0$ in the parameter space $M$ such that
${\cal R}_{\la_0,n}$ moves holomorphically on $\Omega$ for all sufficently big $n$, then there exists a
holomorphic motion $h_{\la} :\Jullo\to\Julla$ which conjugates the dynamics.
\end{lem}

Let us mention here that there exists a much stronger version of the $\la$-lemma, which is due to Slodkowski and shows that the
holomorphic motion actually extends as a quasi-conformal transformation of the full Riemann sphere (see Theorem \ref{lamslo}).
In particular, under the assumption of the above Lemma, $f_\la$ and $f_{\la_0}$ are quasi-conformally conjugated
when $\la\in\Omega$.\\

We may now define the set of stable parameters and its complement; the bifurcation locus.

\begin{defi}\label{defimhj}
Let $\left(f_{\la}\right)_{\la\in M}$ be a  holomorphic family of degree $d$ rational maps.\\
The \emph {stable set} $\cal S$ is the set of
parameters $\la_0\in M$ for which there exists a neighbourhood $\Omega$
of $\la_0$ and a holomorphic motion $h_\la$ of $\Jullo$ over $\Omega$,
centered at $\la_0$, and such that
$f_\la\circ h_{\la}=h_\la\circ f_{\la_0}$ on $\Jullo$.\\
The \emph {bifurcation locus} $\Bif$ is the complement $M\setminus {\cal S}$.
\end{defi}

By definition, $\cal S$ is an open subset of $M$ but it is however not yet clear that it is not empty. We shall actually show that $\cal S$ is dense in $M$.
To this purpose we will prove that the stability is characterized by the stability of the critical orbits. The next subsection will be devoted to this simple but
 remarkable fact.

\subsection{Stability of critical orbits}

Let us start by explaining why bifurcations are related with the instability of critical orbits.\\

As Lemma \ref{lemmotcyc} shows, a parameter $\la_0$ belongs to the bifurcation locus if
for any neighbourhood $\Omega$ of $\la_0$ in the parameter space $M$ there exists $n_0\ge 0$ for which
${\cal R}_{\la_0,n_0}$ does not move holomorphically on $\Omega$.\\
It is not very difficult to see that this forces one of the repelling $n_0$-cycles of $f_{\la_0}$, say $R_{\la_0}$, to
become neutral and then attracting for a certain value $\la_1 \in \Omega$. 
Now comes the crucial point. A classical result asserts that the basin of attraction of any attracting cycle of a rational map contains a critical point (see \cite{BM} Th\'eor\`eme II.5). Thus, one of the critical orbits $f_{\la}^{k}\left(c_i(\la)\right)$ is uniformly converging to $R_\la$ on a neighbourhood of $\la_1$. Then the sequence
$f_{\la}^{k}\left(c_i(\la)\right)$ cannot be normal on $\Omega$ since otherwise, by Hurwitz lemma, it should converge uniformly to $R_\la$ which is repelling for $\la$ close to $\la_0$. This arguments show that the bifurcation locus is contained in the set of parameters around which the post-critical set does not move continuously.\\ 

This is an extremely important observation because it will allow us to detect bifurcations by considering only the critical orbits. It leads to the following definitions.\\

\begin{defi}\label{defiact}
Let $\left(f_{\la}\right)_{\la\in M}$ be a  holomorphic family of degree $d$ rational maps. A marked critical point $c(\la)$ is said to be 
\emph{passive} at $\la_0$ if the sequence  $\big(f_{\la}^{n}\left(c(\la)\right)\big)_n$ is normal on some neighbourhood of $\la_0$.
If $c(\la)$ is not passive at $\la_0$ one says that it is \emph{active}.
The \emph{activity locus} of $c(\la)$ is the set of parameters at which $c(\la)$ is active. 
\end{defi}

The key result may now be given.

\begin{lem}\label{lembifa}
In a holomorphic family with marked critical points the bifurcation locus coincides with the union of the activity loci of the critical points.
\end{lem}

\proof According to our previous arguments, the bifurcation locus is contained in the union of the activity loci.
It remains to show that, for any marked critical point $c(\la)$, the sequence $\big(f_{\la}^{n}\left(c(\la)\right)\big)_n$
is normal on $\cal S$.
Assume that $\la_0\in {\cal S}$. As $\Jullo$ is a perfect compact set, we may find three distinct points $a_1, a_2,a_3$ on $\Jullo$ which are avoided by the orbit of $c(\la_0)$ .
Since the holomorphic motion $h_\la$ of $\Jullo$ conjugates the dynamics, the orbit of $c(\la)$ avoids $\{h_\la(a_j);\;1\le j\le3\}$ for all $\la$ in a small neighbourhood of $\la_0$. The conclusion then follows from Picard-Montel's Theorem.\\
\qed\\

The following Lemma is quite useful.

\begin{lem}\label{lemappperno}
If $\la_0$ belongs to the activity locus of some marked critical point $c(\la)$ then there exists a sequence of parameters $\la_k\to\la_0$ such that
$c(\la_k)$ belongs to some super-attracting cycle of $f_{\la_k}$ or is strictly preperiodic to some repelling cycle of $f_{\la_k}$.
\end{lem}

\proof
To simplify we will assume that all critical points are marked.
Since $c(\la)$ is active at $\la_0$ the critical point $c(\la_0)$ cannot be periodic. In particular, its pre-history cannot stay in the critical set.
Then, after maybe replacing $c(\la_0)$ by another critical point, we may assume that there exist
holomorphic maps $c_{-2}(\la),\;c_{-1}(\la)$ near $\la_0$ such that $f_\la\left(c_{-2}(\la)\right)= c_{-1}(\la)$,
$f_\la\left(c_{-1}(\la)\right)= c(\la)$ and $Card\;\{c_{-2}(\la),\;c_{-1}(\la),c(\la)\}=3$. Now, by Picard-Montel Theorem, the sequence
 $\big(f_{\la}^{n}\left(c(\la)\right)\big)_n$ cannot avoid the set $\{c_{-2}(\la),c_{-1}(\la),c(\la)\}$ on any neighbourhood of $\la_0$.\\
A similar argument, using Picard-Montel Theorem and a repelling cycle of period $n_0\ge 3$, shows that 
 $c(\la)$ becomes strictly preperiodic for $\la$ arbitrarily close to $\la_0$.
\qed\\

We are now ready to state and prove Ma\~{n}\'e-Sad-Sullivan theorem.

\begin{theo}\label{theoMSS}
Let $\left(f_{\la}\right)_{\la\in M}$ be a holomorphic family of degree $d$ rational maps with marked critical points
$\{c_1(\la),\cdot\cdot\cdot,c_{2d-2}(\la)\}$.\\
A parameter $\la_0$ is stable if one of the following equivalent conditions is satisfied.
\begin{itemize}
\item[1)] $\Jullo$ moves holomorphically around $\la_0$ (see definition \ref{defimhj})
\item[2)]  the critical points are passive at $\la_0$ (see definition \ref{defiact})
\item[3)] $f_\la$ has no unpersistent neutral cycles for $\la$ sufficently close to $\la_0$.
\end{itemize}
The set $\cal S$ of stable parameters is dense in $M$.
\end{theo}

\proof
The equivalence between 1) and 2) is given by Lemma \ref{lembifa}. Similar arguments may allow to show that 3) is an equivalent statement
(we will give an alternative proof of that later). It remains to show that $\cal S$ is dense in $M$. According to Lemma \ref{lemappperno} we may perturb $\la_0$
and assume that $f_{\la_0}$ has a superattracting cycle of period bigger than 3 which persists, as an attracting cycle, for $\la$ close enough to $\la_0$.
 If $\la_0$ is still active, Picard-Montel's theorem shows that a new perturbation
 guarantees that a critical point falls in the attracting cycle  and, therefore, becomes passive. Since the number of critical points is finite, we may make all critical points
passive after a finite number of perturbations.
\qed\\

\begin{exa} In the quadratic polynomial family, the bifurcation locus is the boundary of  the connectivity locus (or the Mandelbrot set). 
Indeed, it follows immedialey from the definition \ref{defimandel} of the Mandelbrot set that its boundary is the activity locus.
\end{exa}

\begin{exa} The situation is more complicated in the family $\big(P_{c,a}\big)_{(c,a)\in \tiny{\CC^{d-1}}}$ of degree $d$ polynomials when $d\ge 3$ (see subsection \ref{ssPd}).
Indeed, Theorem \ref{controlinfty} shows that the bifurcation (i.e. activity) locus is not bounded since it coincides with the boundary of $\cup_{0\le i\le d-2} {\cal B}_i$ (where ${\cal B}_i$ is the set of parameters for which the orbit of
the critical point $c_i$ is bounded) while the connectedness locus  $\cap_{0\le i\le d-2} {\cal B}_i$ is bounded.
\end{exa}

Although we shall not use it, we end this section by quoting a very interesting classification of the activity situations which is due to Dujardin and Favre (see \cite{DF} Theorem 4).

\begin{theo}
Let $\left(f_{\la}\right)_{\la\in M}$ be a holomorphic family of degree $d$ rational maps with a marked critical point
$c(\la)$. If $c$ is passive on some connected open subset $U$ of $M$ then exactly one of the following cases holds:
\begin{itemize}
\item[1)] $c$ is never preperiodic in $U$ and the closure of its orbits move holomorphically on $U$
\item[2)] $c$ is persitently preperiodic on $U$
\item[3)] the set of parameters for which $c$ is preperiodic is a closed subvariety in $U$. Moreover, either there exists a persistently attracting cycle
attracting $c$ throughout $U$, or $c$ lies in the interior of a linearization domain associated to a persistent irrationally neutral periodic point.
\end{itemize}
\end{theo}

It is worth emphasize that the proof of that result relies on purely local, and subtle, arguments.

\subsection{Some remarkable parameters}\label{ssRP}

\noindent$\bullet$ As we already mentionned,  the basin of any attracting cycle of a rational map contains a critical point. As a consequence,
any rational map of degree $d$ has at most $2d-2$ attracting cycles. Then, using  perturbation arguments, Fatou and Julia proved that the number of non-repelling cycles
is bounded by $6d-6$. The precise bound has been obtained by Shishikura \cite{Shi1} using quasiconformal surgery, Epstein has given a more algebraic proof based on  quadratic
differentials (see \cite{epstein2}). 

\begin{theo}\label{TheoShi}
A rational map of degree $d$ has at most $2d-2$ non-repelling cycles.
\end{theo}

In particular, any degree $d$ rational map cannot have more than $2d-2$ neutral cycles. Shishikura has also shown that the bound
$2d-2$ is sharp (we will give another proof, using bifurcation currents, in subsection \ref{densShiHyp}). Let us mention that the Julia set of  a degree $d$ map having $2d-2$ Cremer cycles coincides with the full Riemann sphere. These results motivate the following definition.

\begin{defi}\label{defiShishi}
 The set $\Shi$ of degree $d$ \emph{Shishikura rational maps} is defined by
\begin{center}
$\Shi=\{f\in Rat_d\;/\; f\;\textrm{has}\;2d-2\;\textrm{neutral cycles}\}$.
\end{center}
In a holomorphic family $\left(f_{\la}\right)_{\la\in M}$ we shall denote by $\Shi(M)$ the set of parameters $\la$
for which $f_{\la}$ is Shishikura.
\end{defi}

According to Theorem \ref{theoMSS}, one has $\Shi\subset\Bif$. In the last chapter, we will obtain some informations on the geometry of the set $\Shi$ and, in particular, reprove that
it is not empty.\\

\noindent$\bullet$ Any repelling cycle of $f\in Rat_d$ is an invariant (compact) set on which $f$ is uniformly expanding.
Some rational map may be uniformly expanding on much bigger compact sets. Such sets are called \emph{hyperbolic}
and are necessarily contained in the Julia set. A rational map which is uniformly expanding on its Julia set is said to be hyperbolic.
Let us give a precise definition.

\begin{defi}
Let $f$ be a rational map. A compact set $K$ 
of the Riemann sphere is said to be \emph{hyperbolic} for $f$ if it is invariant
($f(K)\subset K$) and if there exists $C>0$ and $M>1$ such that
\begin{center}
$\vert (f^n)'(z)\vert_{s} \ge CM^n\;;\;\;\forall z\in K,\;\forall n\ge 0$.
\end{center}
One says that $f$ is uniformly expanding on $K$ (recall that $\vert\;\vert_s$ is the spherical metric).
\end{defi}

It may happen that a critical orbit is captured by  some hyperbolic set and, in particular, by a repelling cycle. Such rational maps play 
a very important role since they allow to define transfer maps carrying  informations from the
dynamical plane to the parameter space. A particular attention will be devoted to those having all critical orbits captured by a hyperbolic set.

\begin{defi}\label{defiMis}
 The set $\Mis$ of degree $d$ \emph{Misiurewicz rational maps} is defined by
\begin{center}
$\Mis=\{f\in Rat_d\;/\;  \textrm{all critical orbits of}\;f\;\textrm{are captured by a compact hyperbolic set}\}.$
\end{center}
When the hyperbolic set is an union of repelling cycles the map is said to be strongly Misiurewicz and the set of such maps is denoted $\SMis$.
In a holomorphic family $\left(f_{\la}\right)_{\la\in M}$ we shall denote by $\Mis(M)$ (resp. $\SMis(M)$) the set of parameters $\la$
for which $f_{\la}$ is Misiurewicz (resp. strongly Misiurewicz).
\end{defi}

Within a holomorphic family  $\left(f_{\la}\right)_{\la\in M}$,
one may show that a critical point whose orbit is captured by a hyperbolic set and leaves this set under a small perturbation
is active. In particular, we have the following inclusion: $\Mis\subset \Bif$. To prove this, one first has to construct a holomorphic motion of the hyperbolic set
and then linearize along its orbits (see \cite{Gau}).
When the hyperbolic set is a cycle, the motion is given by the implicit function theorem and the linearizability is a well known fact.

\begin{lem}\label{lemAcMi}
Let $\left(f_{\la}\right)_{\la\in M}$ be a  holomorphic family of degree $d$ rational maps with a marked critical point $c(\la)$.
Assume that $f_{\la}$ has a repelling $n$-cycle ${\cal R}(\la):=\{z_\la,f_{\la}(z_\la),\cdot\cdot\cdot,f_{\la}^{n-1}(z_\la)\}$ for $\la\in U$.
If $f_{\la}^k(c(\la)\in {\cal R}(\la)$ for $\la=\la_0\in U$ but not for all $\la\in U$ then $c(\la)$ is active at $\la_0$.
\end{lem}

\proof 
We may assume that $n=1$ which means that $z_\la$ is fixed by $f_\la$. Shrinking $U$ and linearizing we get a 
a family of local biholomorphisms $\phi_{\la}$ which depends holomorphically on $\la$ and 
such that $\phi_\la(0)=z_\la$ and $f_\la \circ \phi_\la ( u)= \phi_\la (m_\la u)$ (see \cite{BM} Th\'eor\`eme II.1 and Remarque II.2). As $z_\la$ is repelling, one has $\vert m_\la\vert >1$.
Let us set $u(\la):=\phi_\la^{-1}(f_{\la}^{k}(c(\la))$,  
then $f_{\la}^{p+k}(c_\la)=f_{\la}^{p}\left(\phi_\la(u(\la))\right)=\phi_\la\left((m_\la)^p  u(\la)\right)$ which shows that $f_{\la}^{p+k}(c_\la)$ is not normal at $\la_0$
since, by assumption, $u(\la_0)=0$ but $u$ does not vanish identically on $U$.

\qed\\

\noindent $\bullet$ As we already mentionned, a rational map which is uniformly expanding on its Julia set is called hyperbolic.
The dynamical study of such maps turns out to be much easier.

\begin{defi}
 The set $\Hyp$ of degree $d$ \emph{hyperbolic rational maps} is defined by
\begin{center}
$\Hyp=\{f\in Rat_d\;/\; f\;\textrm{is uniformly expanding on its Julia set}\}$.
\end{center}
In a holomorphic family $\left(f_{\la}\right)_{\la\in M}$ we shall denote by $\Hyp(M)$ the set of parameters $\la$
for which $f_{\la}$ is hyperbolic.
\end{defi}

There are several characterizations of hyperbolicity. One may show that a rational map $f$ is hyperbolic if and only if its postcitical set does not contaminate its
Julia set: 
\begin{center}
$f$ is hyperbolic $\Leftrightarrow \overline{\cup_{n\ge0}f^n(\Crf)}\cap \Julf=\emptyset$.
\end{center}

As a consequence, in any holomorphic family  $\left(f_{\la}\right)_{\la\in M}$, hyperbolic parameters are stable : $\Hyp(M)\subset {\cal S}$.
This characterization also implies that a hyperbolic map has only attracting or repelling cycles.\\

The Fatou's conjecture asserts that the hyperbolic parameters are dense in any holomorphic family,
it is an open problem even for quadratic polynomials.
According to Ma\~{n}\'e-Sad-Sullivan Theorem it can be rephrased as follows

\begin{conj}
 $\Hyp(M)={\cal S}$ for any holomorphic family $\left(f_{\la}\right)_{\la\in M}$. 
 \end{conj}
 
 Let us also mention that Ma\~{n}\'e, Sad and Sullivan have shown that  hyperbolic and non-hyperbolic parameters
 cannot coexist in the same stable component (i.e. connected component of ${\cal S}$). A given stable component is therefore either \emph{hyperbolic}
(if  all parameters in the component are hyperbolic) or \emph{non-hyperbolic}  (when the component does not contain any 
 hyperbolic parameter). Fatou's conjecture  thus claims that non-hyperbolic components do not exist.
 \\
 
Although the bifurcation locus of the quadratic polynomial family is clearly accumulated by hyperbolic parameters, this is far from being clear in other families and seems to be an interesting question. In the last chapter we will show that parameters which, in some sense, produce the strongest bifurcations, are accumulated by hyperbolic parameters.

\section{Potential theoretic approach}

\subsection{Lyapunov exponent and bifurcations}

We wish here to explain why the current $dd^c L(\la)$, where $L(\la)$ is the Lyapunov exponent of $f_\la$,  is a natural candidate for being
a bifurcation current in any holomorphic family $\left(f_{\la}\right)_{\la\in M}$.\\

To support this idea we will first use the material discussed in subsections \ref{ssLyapMult} and \ref{ssdefper} and relate the Lyapunov exponent with the distribution of hypersurfaces $\Per_n(e^{i\theta})$ in $M$.\\

According to the Ma\~{n}\'e-Sad-Sullivan Theorem \ref{theoMSS}, the bifurcation locus is the closure of the union  of all such hypersurfaces:
\begin{center}
$\Bif=\overline{\cup_n\cup_\theta \Per_n(e^{i\theta})}$.
\end{center}

One may thus expect that the sequence of currents $\frac{1}{2\pi}\int_0^{2\pi} [\Per_n(e^{i\theta})]$, correctly weighted, converges to some current supported by $\Bif$. As we shall see, 
 the approximation formula \ref{theoapprox} actually implies that the following convergence occurs in the sense of currents:
\begin{center}
$\frac{d^{-n}}{2\pi}\int_0^{2\pi} [\Per_n(e^{i\theta})]\;d\theta \to \Bc$.
\end{center}

We will now briefly prove this,  more details and a generalization will be provided by Theorem \ref{theodistaver}.\\

Let us denote
 $p_n(\cdot,e^{i\theta})$  the canonical defining functions for  the hypersurfaces $\Per_n(e^{i\theta})$ given by Theorem \ref{theopoly}.
We have to show that the sequence of $p.s.h$ functions
\begin{center} 
$\displaystyle L_n(\la):=\frac{d^{-n}}{2\pi}\int_0^{2\pi} \ln \vert p_n(\la,e^{i\theta})\vert\;d\theta$.
\end{center}
 converges to $L$ in $L^1_{loc} (M)$.\\

Writting $$p_n(\la,e^{i\theta})=:\prod_{i=1}^{N_d(n)}\big(e^{i\theta}-w_{n,j}(\la)\big)$$

and
using the fact that $\ln^+ \vert a \vert=\frac{1}{2\pi}\int_0^{2\pi}\ln \vert a-e^{i\theta}\vert d\theta$
 yields

\begin{eqnarray}\label{Com}
L_n(\la)=\frac{1}{2\pi d^n}\int_0^{2\pi} \ln \prod_j\vert e^{i\theta}-w_{n,j}(\la)\vert d\theta=
d^{-n}\sum_j \ln^+ \vert w_{n,j}(\la)\vert.
\end{eqnarray}
Using the fact that
$d^{-n}N_d(n)\sim\frac{1}{n}$ (see Theorem \ref{theopoly}) one sees that the sequence $L_n$ is locally bounded from above and it thus remains to see that it converges pointwise to $L$.

 According to Theorem \ref{theopoly}, the set $\{w_{n,j}(\la)\; /\;w_{n,j}(\la)\ne 1\}$ coincides with the set of multipliers
of cycles of exact period $n$ (counted with multiplicity) from which the cycles of multiplier $1$ are deleted.

Since $f_\la$ has a finite number of non-repelling cycles (Fatou's theorem), all cycles appearing in \ref{Com} are repelling for $n$ big enough. The conclusion then follows 
from Theorem \ref{theoapprox}.
\qed\\

As a second evidence, we will relate $dd^c L$ with the instability of the critical dynamics. 
 Ma\~{n}\'e-Sad-Sullivan Theorem \ref{theoMSS} also revealed that the bifurcations are due to the activity of one (or more) critical point. It is not very difficult to check that the
 activity locus of a marked critical point $c(\la)$ is exactly supported by the current $dd^c g_\la\left(c(\la)\right)$ where $g_\la$ is the Green function of $f_\la$.
 A bifurcation current could therefore also be defined by $$dd^c \sum_j  g_\la\left(c(\la)\right)$$ where the sum is taken over the critical points. It turns out that this current coincides with the current $dd^c L$.\\
 
  We will justify this for a polynomial family and, to this purpose, will
 establish a fundamental formula which relates the Lyapunov exponent (see definition \ref{defilyap}) of a polynomial with its critical points. 

\begin{theo}{\bf (Przytycki's formula)}\label{PrzFor}
Let $P$ be a unitary degree $d$ polynomial, $L(P)$ its Lyapunov exponent and $g_{\tiny{\CC},P}$  its Green function. Then
\begin{center}
$L(P)=\ln d + \sum g_{\tiny{\CC},P} (c)$
\end{center}
where the sum is taken over the critical points of $P$ counted with multiplicity.
\end{theo}

\proof
Let us write $c_1,c_2,\cdot\cdot\cdot,c_{d-1}$ the critical points of $P$. Then
\begin{eqnarray*}
L(P)=\int_{\CC} \ln \vert P'\vert\;\mu_P=\int_{\CC} \ln \vert d\prod_{j=1}^{d-1} (z-c_j)\vert\;\mu_P=\ln d+\sum_{j=1}^{d-1}\int_{\CC} \ln  \vert z-c_j\vert\;\mu_P.
\end{eqnarray*}
Now, since $\mu_P=dd^c g_{\tiny{\CC},P}$, the formula immediately follows by an integration by parts.
\qed\\

\begin{rem}
As the Green function $g_{\tiny{\CC},P}$  of a polynomial $P$ is positive (see Proposition \ref{Greenpoly}) the above formula implies that $L(P)\ge \ln d$.
\end{rem}

 Most of the remaining of this chapter will be devoted to extend these ideas to  families of rational maps.
 
\subsection{Przytycki's generalized formula}

The proof of Przytycki's  formula relies on a simple integration by parts.
It is more delicate to perform such an integration by part in the case of a rational map $f$.
To this purpose we will work in the line bundle ${\cal O}_{\pp}(D)$ for $D:=2(d-1)$ which we endow with two metrics, the flat one $\|[z,x]\|_0$ and the Green metric
 $\|[z,x]\|_{G_F}$ whose potential is the Green function $G_F$ of some lift $F$ of $f$ (see subsection \ref{sslyap}).\\
 
 In this general situation the integration by part yields the following formula. We refer to definition \ref{DefiGreen} for the definitions and notations related to the Green functions.

\begin{prop}\label{propPrzgen}
Let $f$ be a rational map of degree $d\ge 2$ and $F$ be one of its lifts. Let $D:=2(d-1)$
and $Jac_F$ be the holomorphic section of ${\cal O}_{\pp}(D)$ induced  by $\det F'$.
Let $g_F$ be the Green function of $f$ on $\pp$ and $\mu_f$ be the  Green measure of $f$. Then
\begin{center}
$L(f)+\ln d=\int_{\pp} g_F[\Crf]\;-\;2(d-1)\int_{\pp} g_F(\mu_f+\omega)\;+\;\int_{\pp}\ln\|Jac_F\|_{0}\;\omega$.
\end{center}
\end{prop}

\proof
Let us recall that $\Vert\cdot\Vert_G=e^{-Dg_F}\Vert\cdot\Vert_0$.
According to Lemma \ref{lemLyapBund} we have:
\begin{eqnarray*} 
L(f)+\ln d=\displaystyle\int_{\pp}\ln\|Jac_F\|_{G_F}\;\mu_f=\displaystyle\int_{\pp}\ln\|Jac_F\|_{G_F}\;dd^c g_F+\displaystyle\int_{\pp}\ln\|Jac_F\|_{G_F}\;\omega
\end{eqnarray*}
which, after integrating by parts, yields
\begin{eqnarray*} 
L(f)+\ln d=\displaystyle\int_{\pp} g_F\left( dd^c\ln\|Jac_F\|_{G_F}\right)+\displaystyle\int_{\pp}\ln\|Jac_F\|_{G_F}\;\omega
\end{eqnarray*}
and by Poincar\'e-Lelong formula:
\begin{eqnarray} \label{przgen}
L(f)+\ln d=\displaystyle\int_{\pp}g_F\left([\Crf]-D\mu_f\right)+\displaystyle\int_{\pp}\left(\ln\|Jac_F\|_{0}-Dg_F\right)\;\omega.
\end{eqnarray}
\qed\\

When working with a holomorphic family of polynomials $\left(P_\la\right)_M$, Przytycki's formula says that the Lyapunov function $L(P_\la)$ and the 
sum of values of the Green function on critical points differ from a constant. In particular, $L(P_\la)$ is a $p.s.h$ function on $M$ and these two functions induce the same $(1,1)$ current on $M$. It is then rather clear, using Ma\~{n}\'e-Sad-Sullivan Theorem \ref{theoMSS}, that this current is exactly supported by the bifurcation locus.\\

 We aim to generalize this to holomorphic families of rational maps and will therefore 
compute the $dd^c$ of the left part of formula \ref{przgen}. The following formula has been established in \cite{BB1}.

\begin{theo}\label{theoformula}
Let $\left(f_\la\right)_M$ be a holomorphic family of degree $d$ rational maps which admits a holomorphic family of lifts $\left(F_\la\right)_M$.
 Let $p_{\tiny M}$ (resp $p_{\tiny\pp}$) be the canonical projection from $M\times\pp$
onto $M$ (resp.$\pp$). Then
\begin{eqnarray}\label{ddcprzgen}
dd^c L(\la) =(p_M)_{\star}\left( (dd^c_{\la,z} g_\la(z)+\hat\omega)\wedge [C]\right)
\end{eqnarray}
where $C:=\{(\la,z)\in M\times\pp\;/\; z\in \Crla\}$,
$g_\la:=g_{F_\la}$ is the Green function of $f_\la$ on $\pp$, $L(\la)$ the Lyapunov exponent of $f_\la$ and $\hat\omega:=p_{\pp}^{\star} \omega$.\\
Moreover, the function $\la\mapsto \int_{\pp} g_\la(\mu_\la+\omega)$ is pluriharmonic on $M$.
\end{theo}

\proof  
We may work locally and define a holomorphic section  $Jac_\la$ of ${\cal O}_{\pp}(D)$ induced  by $\det F_\la'$.
 Then $\widetilde{Jac}(\la,[z]):=\left(\la,Jac_\la([z])\right)$ is a holomorphic section of the line bundle $M\times {\cal O}_{\pp}(D)$ over $M\times\pp$.\\

Let us rewrite 
Proposition \ref{propPrzgen} on the form $L(\la)+\ln d=H(\la)-D\;B(\la)$
where $H(\la):=\int_{\pp} g_\la[\Crla]\;+\;\int_{\pp}\ln\|Jac_\la\|_{0}\;\omega$, $B(\la):=\int_{\pp} g_\la(\mu_\la+\omega)$ and $D:=2(d-1)$.\\

We  first compute $dd^c H$. Let $\Phi$ denote a $(m-1,m-1)$ test form where $m$ is the dimension of $M$. Then 
$\langle dd^c H,\Phi\rangle=I_1+I_2$ where 
$I_1:=\int_M  dd^c \Phi\int_{\pp} g_\la[\Crla]$ and $I_2:=\int_M  dd^c \Phi\int_{\pp} \ln\|Jac_\la\|_{0}\;\omega$.\\

Slicing and then integrating by parts, we get 
\begin{eqnarray*}
I_1=\int_M dd^c\Phi\int_{\pp} \left(g[C]\right)_\la=\int_{M\times\pp} (p_M)^{\star}\left(dd^c \Phi\right)\wedge \left(g[C]\right)=\\
=\int_{M\times\pp} (p_M)^{\star}\Phi\wedge dd^c g\wedge [C].
\end{eqnarray*}

Performing the same computation and then using the Poincar\'e-Lelong identity $dd^c  \ln\Vert \widetilde{Jac}\Vert_0=[C]-D\hat\omega$ one gets
\begin{eqnarray*}
I_2=\int_M dd^c\Phi\int_{\pp} \left(\ln\Vert\widetilde{Jac}\Vert_0\;\hat\omega\right)_\la=\int_{M\times\pp} (p_M)^{\star}\left(dd^c \Phi\right)\wedge \ln\Vert\widetilde{Jac}\Vert_0\;\hat\omega=\\
=\int_{M\times\pp} (p_M)^{\star}\Phi\wedge dd^c\ln\Vert \widetilde{Jac}\Vert_0 \wedge \hat\omega=
\int_{M\times\pp} (p_M)^{\star}\Phi\wedge \hat\omega\wedge [C].
\end{eqnarray*}

This shows that $dd^c H =(p_M)_{\star}\left( (dd^c_{\la,z} g_\la(z)+\hat\omega)\wedge [C]\right)$.\\

 It remains to show that $dd^c B =0$.
It seems very difficult to prove this by calculus, we will use a trick which exploits the dynamical situation. If one replace the family $\left(f_\la\right)_M$ by the family
 $\left(f_\la^2\right)_M$ then $L$ becomes $2L$ and $dd^c H$ becomes $2dd^c H$ while $B$ is unchanged. Thus, applying Proposition \ref{propPrzgen} to $\left(f_\la\right)_M$ and taking $dd^c$ yields
 $dd^c L=dd^c H-2(d-1)dd^c B$ 
but, with the family $\left(f_\la^2\right)_M$, this yields $2dd^c L=2dd^c H-2(d^2-1)dd^c B$. By comparison one obtains $dd^c B =0$.\qed\\

\begin{rem} The paper \cite{BB1} actually covers  the case of holomorphic families of endomorphisms of $\PP$. In this setting, 
the formula \ref{ddcprzgen} becomes
\begin{eqnarray*}
dd^c L(\la) =(p_M)_{\star}\left( (dd^c_{\la,z} g_\la(z)+\hat\omega)^k\wedge [C]\right)
\end{eqnarray*}
where $L(\la)$ is now the sum of Lyapunov exponents of $f_\la$ with respect to the Green measure $\mu_\la$.
\end{rem}

\subsection{Activity currents and the bifurcation current}

We have seen with Ma\~{n}\'e-Sad-Sullivan Theorem \ref{theoMSS} that the bifurcations, within a holomorphic family of degree $d$ rational maps
$\left(f_{\la}\right)_{\la\in M}$, are due to the activity of the critical points (see in particular Lemma \ref{lembifa}). We will use this and the formula given by Theorem
\ref{theoformula}, to define a closed positive $(1,1)$-current $\Bc$ on $M$ whose support is the bifurcation locus and which admits the Lyapunov function
as global potential.
Besides, we will also introduce  a collection of $2d-2$  closed positive $(1,1)$-currents which detect the activity of each critical point and see that the bifurcation current
$\Bc$ is the sum of these currents.
Here are the formal definitions.\\

\begin{defi}\label{defiBC} Let $\left(f_{\la}\right)_{\la\in M}$ be any holomorphic family of degree $d$ rational maps
with marked critical points $\{c_i(\la);\;1\le i\le 2d-2\}$. 
The \emph{bifurcation current} $\Bc$ is defined by 
\begin{center}
$\Bc:=dd^c L(\la)$
\end{center}
where $L(\la)$ is the Lyapunov exponent of $f_\la$ with respect to its Green measure.\\
The \emph{activity current} $T_i$ of the marked critical point
$c_i$ is defined by
\begin{center}
$T_i:=(p_M)_{\star}\left( (dd^c_{\la,z} g_\la(z)+\hat\omega)\wedge [C_i]\right)$
\end{center}
where $C_i$ is the graph $\{(\la,c_i(\la);\;\la\in M\}$ in $M\times\pp$
and $g_\la:=g_{F_\la}$  is the Green function of $f_\la$ on $\pp$ for some (local) holomorphic family of lifts $\left(F_\la\right)$.\\
\end{defi}

Let us observe that the curent $dd^c_{\la,z} g_\la(z)+\hat\omega$ somehow interpolates the Green measures $\mu_\la$ of $f_\la $ and is actually dynamically 
obtained. Indeed, one may consider the map $\hat f : M\times \pp\to M\times\pp$ defined
by $\hat f(\la,z):=(\la,f_\la(z))$ and then get it by taking the limit of the sequence of currents $d^{-n}(\hat f )^n (\hat \omega)$ in $M\times \pp$.
To establish the convergence one proceeds like in Lemma \ref{lemexg}.
The activity current $T_i$ is the projection on $M$ of the restriction of this current to the hypersurface $C_i$.\\

We now give local potentials for the activity currents.

\begin{lem}\label{lemPotAC}
Let $F_\la$ be a local holomorphic family of lifts of $f_\la$ and $G_\la$ be the Green function of $F_\la$ on $\CC^{2}$.
Let $\hat c_i(\la)$ be a local lift of $c_i(\la)$.
Then $G_\la\left(\hat c_i(\la)\right)$ is a local potential of $T_i$.
\end{lem}

\proof
This is a straightforward computation using the fact that, for any local section $\sigma$ of the canonical projection
$\pi: \CC^{2}\setminus\{0\}\to \pp$, the function $G_\la\left(\sigma(z)\right)$ is a local potential of $dd^c_{\la,z} g_\la(z)+\hat\omega$
(see Proposition \ref{PropG_F}).\qed\\ 

The following result has been originally proved by DeMarco in \cite{DeM1}, \cite{DeM2}. 

\begin{theo}\label{theosuppBC}
Let $\left(f_{\la}\right)_{\la\in M}$ be any holomorphic family of degree $d$ rational maps
with marked critical points $\{c_i(\la);\;1\le i\le 2d-2\}$. The support of the activity current $T_i$ is the activity locus
of the marked critical point $c_i$.\\
The support of the bifurcation current $\Bc$ is the bifurcation locus of the family
$\left(f_{\la}\right)_{\la\in M}$ and $\Bc=\sum_{1}^{2d-2} T_i$.
\end{theo}

\proof
Let us first show that $c_i$ is passive on the complement of $Supp\;T_i$.
If $T_i=0$ on a small ball $B\subset M$ then, by Lemma \ref{lemPotAC}, $G_\la\left(\hat c_i(\la)\right)$ is pluriharmonic 
and therefore equal to $\ln \vert h_i(\la)\vert$ for some non-vanishing holomorphic function $h_i$ on $B$. Replacing 
$\hat c_i(\la)$ by $\frac{\hat c_i(\la)}{h_i(\la)}$ one gets, thanks to the homogeneity property of $G_\la$ (see Proposition
\ref{PropG_F}), $G_\la\left(\hat c_i(\la)\right)=0$. This implies that 
\begin{center}
$\{F_\la^n(\hat c_i(\la))\;/\; n\ge 1, \la\in B\} \subset \cup_{\la\in B} G_\la^{-1}\{0\}$
\end{center}
where, after reducing $B$,  the set $\cup_{\la\in B} G_\la^{-1}\{0\}$ is a relatively compact in $\CC^2$. Montel's theorem then tells us that $\left(F_\la^n(\hat c_i(\la))\right)_n$ and, thus,
$\left(f_\la^n( c_i(\la))\right)_n$ are normal on $B$.\\

Let us now show that $T_i$ vanishes where $c_i$ is passive. 
Assume that a subsequence $\left(f_\la^{n_k} (c_i(\la))\right)_k$ is uniformly converging  on a small ball $B\subset M$.
Then we may find a local section $\sigma$ of $\pi: \CC^{2}\setminus\{0\}\to \pp$ such that
$F_\la^{n_k}(\hat c_i(\la))=h_{n_k}(\la)\cdot \sigma\circ f_\la^{n_k} (c_i(\la))$
where $h_{n_k}$ is a non-vanishing holomorphic function on $B$. As $G_\la\circ F_\la=dG_\la$ (see Proposition \ref{PropG_F}),
this yields
\begin{center}
$G_\la(\hat c_i(\la))=d^{-n_k}\big( \ln \vert h_{n_k}(\la)\vert + G_\la \circ\sigma\circ f_\la^{n_k} (c_i(\la))\big)$
\end{center}
which, after taking $dd^c$ and making $k\to+\infty$, implies that $T_i$ vanishes on $B$.\\

That $dd^c L=\Bc=\sum_{1}^{2d-2} T_i$ follows immediately from Theorem
\ref{theoformula}. Then, Ma\~{n}\'e-Sad-Sullivan Theorem \ref{theoMSS} implies that the support of $\Bc$ is the bifurcation locus. 
\qed\\

\begin{rem}
As proved by Ma\~{n}\'e,  the product of the Lyapunov exponent of a  degree $d$ rational map and the Hausdorff dimension of its Green measure  is equal to $\ln d$ (\cite{Mane}).
This suggest that for bifurcation investigations it could be more natural to study the borelian support of the measure  than its topological support.
\end{rem}

It is important to stress here that the identity $dd^c L=\Bc=\sum_{1}^{2d-2} T_i$ may be seen as a potential-theoretic expression of 
Ma\~{n}\'e-Sad-Sullivan theory. In practice, the expression $\Bc=dd^c L$ will be used to investigate the set of parameters $\Shi$ while the expression
$\Bc=\sum_{1}^{2d-2} T_i$ will be used for the parameters $\Mis$ (see subsection \ref{ssRP}).\\ 
In the next chapters, we will deeply use the fact that the Lyapunov function is a potential of $\Bc$, together with the approximation formula 
\ref{theoapprox}, to analyse how the hypersurfaces $Per_n(w)$ may shape
the bifurcation locus.\\ The continuity of the Lyapunov function $L$ will turn out to be  decisive for this study.
Ma\~{n}\'e was the first to establish the continuity of $L$ (see \cite{Mane}) , using Theorem \ref{theoformula} and Lemma \ref{lemPotAC}
one may see that this function is actually H\"{o}lder continuous.

\begin{theo}\label{theoLH}
The Lyapunov function of any holomorphic family of degree $d$ rational maps
is $p.s.h$ and H\"{o}lder continuous.
\end{theo}

\proof
According to Theorem \ref{theoformula} and Lemma \ref{lemPotAC} the functions $L$ and $\sum G_\la(\hat c_i(\la))$ differ from a
pluriharmonic function. The conclusion follows from the fact that $G_\la(z)$ is H\"{o}lder continuous in $(\la,z)$
(see \cite{BB1} Proposition 1.2.).
\qed\\

\subsection{DeMarco's formula}

Using Theorem \ref{theoformula} we will get an explicit version of the formula given by Proposition \ref{propPrzgen}.
 This result was first obtained by DeMarco who used a completely different
 method. We refer to the paper of Okuyama \cite{Oku} for yet another proof.\\

The key will be to compute the integral $\int_{\pp} g_F\left(\mu_f+\omega\right)$ which appears in the
formula 
given by Proposition \ref{propPrzgen}. To this purpose we shall use the \emph{resultant} of a homogeneous polynomial map $F$ of degree $d$
on $\CC^2$. The space of such maps can be identified with $\CC^{2d+2}$. The  resultant $Res\;F$ of $F$ polynomialy depends on $F$ and vanishes
if and only if $F$ is degenerate. Moreover $Res\;(z_1^d,z_2^d)=1$ and $Res$ is $2d$-homogeneous: $Res\;aF=a^{2d}\;Res\;F$.

\begin{lem}\label{lemB}
$\int_{\pp} g_F\left(\mu_f+\omega\right)=\frac{1}{d(d-1)}\ln \vert Res\;F\vert -\frac{1}{2}$
\end{lem}

\proof
The function $B(F):=\int_{\pp} g_F\left(\mu_f+\omega\right)$ is well defined on $\CC^{2d+2}\setminus \Sigma$ where $\Sigma$ is the hypersurface
where $Res$ vanishes. Moreover, according to Theorem \ref{theoformula}, $B$ is pluriharmonic. 
As $B$ is locally bounded from above, it extends to some $p.s.h$ function through $\Sigma$.
Then, by Siu's theorem, there exists some positive constant $c$ such that $dd^c B(F)=c\;dd^c \ln \vert Res F\vert$
which means that $B-c\ln \vert Res F\vert$ is pluriharmonic on $\CC^{2d+2}$.\\
Let $\varphi$ be a non-vanishing holomorphic function on $\CC^{2d+2}$ such that 
\begin{center}
$B-c\ln \vert Res F\vert=\ln \vert \varphi\vert$.
\end{center}
Using the homogeneity of $Res$ and the fact that $B(aF)=\frac{2}{d-1} \ln \vert a\vert +B(F)$ (one easily checks that $g_{aF}=
\frac{1}{d-1} \ln \vert a\vert +g_F$) one gets:
\begin{center}
$\vert \varphi(aF)\vert = \vert a\vert^{\frac{2}{d-1}-2cd}$
\end{center}
Making $a\to 0$ one sees that $c=\frac{1}{d(d-1)}$ and $\varphi$ is constant. To compute this constant one essentially tests the formula
on $F_0:=(z_1^d,z_2^d)$ (see \cite{BB1}, Proposition 4.10).\qed\\

We are now ready to prove the main result of this subsection.

\begin{theo}\label{theoDeMFor}
Let $f$ be a rational map of degree $d\ge 2$ and $F$ be one of its lifts.
Let $G_F$ be the Green function of $F$ on $\CC^{2}$ and $Res\;F$ be the resultant of $F$. Then
\begin{center}
$\displaystyle L(f)+\ln d=\sum_{j=1}^{2d-2} G_F(\hat c_j) -\frac{2}{d} \ln \vert Res\; F\vert$
\end{center}
where $\hat c_1,\hat c_2,\cdot\cdot\cdot,\hat c_{2d-2}$ are chosen so that
$det F'(z)=\prod_{j=1}^{2d-2} \hat c_j\wedge z$.
\end{theo}

\proof
Taking Lemma \ref{lemB} into account, the formula 
given by Proposition \ref{propPrzgen} becomes
\begin{center}
$L(f)+\ln d=\int_{\pp} g_F[\Crf]\;-\;2(d-1)\left(\frac{1}{d(d-1)}\ln \vert Res\;F\vert -\frac{1}{2}\right)\;+
\;\int_{\pp}\ln\|Jac_F\|_{0}\;\omega.$
\end{center}
Observe that $\omega=\pi_{\star} m$ where $m:=\left(dd^c \ln^+\Vert\cdot\Vert\right)^2$ is the normalized Lebesgue measure on the euclidean unit sphere
of $\CC^2$. Then 
\begin{center}
$\int_{\CC^2} \ln \vert det F'\vert\;m=\int_{\CC^2} \ln\left(  e^{-D\Vert \cdot\Vert}\vert det F'\vert\right)\;m=
\int_{\CC^2}\ln \Vert J_F\circ \pi\Vert_0\; m=\int_{\pp}\ln\Vert Jac_F\Vert_{0}\;\omega.$
\end{center}
Let us pick $U_j$ in the unitary group of $\CC^2$ such that $U_j^{-1}(\hat c_j)=\left(\Vert \hat c_j\Vert,0\right)$. Then $U_j(z)\wedge \hat c_j=-z_2\Vert \hat c_j\Vert$
and, since $\int_{\CC^2} \ln \vert z_2\vert\, m=-\frac{1}{2}$ one gets
\begin{center}
$\int_{\pp}\ln\Vert Jac_F\Vert_{0}\;\omega=\int_{\CC^2} \ln \vert det F'\vert\;m=\sum_j \int_{\CC^2}\ln \vert U_j(z)\wedge \hat c_j \vert =
\sum_j \Vert\hat c_j\Vert -(d-1)$.
\end{center}
On the other hand, $\int_{\pp} g_F[\Crf]=\sum_j g_F\circ \pi (\hat c_i)=\sum_j G_F(\hat c_j)-\sum_j\ln \Vert \hat c_j\Vert$ and the conclusion
follows.\qed

\chapter{Equidistribution towards the bifurcation current}

In this chapter we study how the bifurcation current may be  approximated by various weighted hypersurfaces which are dynamically defined.
As an application, we construct holomorphic motions and describe the laminated structure of the bifurcation locus in certain regions of the moduli space $Mod_2$. 
 
\section{A general method}\label{Sgs}

We will present a general method which may be used to prove that a sequence of currents $T_n:=dd^{c} h_n$, which admit global $p.s.h$ potentials $h_n$
on the parameter space $\CC^k$ of some holomorphic family $\big(f_\la\big)_{\CC^k}$, is converging to the bifurcation current $\Bc$.\\

 The method itself only relies on standard potential-theoretic arguments but dynamical informations are then necessary to  apply it. More specifically, one needs to control the bifurcation current near infinity,  this is why the method works well in polynomial families but seems difficult to handle in more general families.\\
   In section \ref{SecCPP} we shall use it  to study the distribution of critically periodic parameters
 and in subsection \ref{CNC} to study the distribution of the hypersurfaces $Per_n(e^{i\theta})$.\\

We have to prove that $(h_n)_n$ converges in $L_{loc}^1$ to $L$ and, 
for this,
will use a well known \emph{compactness principle} for subharmonic functions:

\begin{theo}\label{compapsh}
Let $(\varphi_j)$ be a sequence of subharmonic functions  which is locally uniformly bounded from above on some domain $\Omega \subset \RR^n$. 
 If $(\varphi_j)$ does not converge
to $-\infty$ then a subsequence $(\varphi_{j_k})$  converges in $L^1_{loc}(\Omega)$ to some subharmonic function $\varphi$. 
In particular, $(\varphi_j)$ converges in $L^1_{loc}(\Omega)$ to some subharmonic function $\varphi$ if it converges pointwise  to $\varphi$.
\end{theo}

We thus first need to check that $(h_n)_n$ is locally bounded from above and does not converge to $-\infty$,
the following result may be used to check this last property.
 
\begin{lem}\label{lemHar}{\bf (Hartogs)}
Let $(\varphi_j)$ be a sequence of subharmonic functions  and $g$ be a continuous function defined on some domain $\Omega \subset \RR^n$.
If
$\limsup_j \varphi_j(x)\le g(x)$ for every $x\in \Omega$ then, for any compact $K\subset \Omega$ and every $\epsilon>0$ one has
$\varphi_j(x)\le g(x)+\epsilon $ on $K$ for $j$ big enough.
\end{lem}

Then one has to show that 
$L$ is the unique limit value of $(h_n)_n$ for the $L_{loc}^1$ topology. To this purpose, we shall apply the following generalized maximum principle to some limit value $\varphi$ of $(h_n)_n$ and to $\psi=L$.

\begin{lem}\label{lemPmax}
Two $p.s.h$ functions $\varphi, \psi$  on 
 $\CC^{k}$ coincide if the following conditions are satisfied:
 \begin{itemize}
 \item[i)] $\psi$ is continuous
 \item[ii)] $\varphi\le \psi$
 \item[iii)]  $Supp\;(dd^{c} \varphi) \subset Supp\;(dd^{c} \psi)$
\item[iv)] $\varphi=\psi$ on $Supp\;(dd^{c} \psi)$
\item[v)] for any $\la_0\in \CC^{k}$ there exists a complex line ${\cal L}$ through $\la_0$ such that $\varphi=\psi$ on 
 the unbounded component of ${\cal L}\setminus\big({\cal L}\cap Supp\;(dd^{c} \psi)\big)$.
 \end{itemize}

\end{lem}

\proof Because of iv), we only have to show that $\varphi(\la_0)=\psi(\la_0)$ when $\la_0$ lies in the complement of
$Supp\;(dd^{c} \psi)$.
According to v), we may find a complex line ${\cal L}$ in $\CC^{k}$ containing $\la_0$ and such that $\psi$ and $\varphi$ coincide on the unbounded
component $\Omega_{\infty}$ of ${\cal L}\setminus\big({\cal L}\cap Supp\;(dd^{c} \psi)\big)$. We may therefore 
assume that
$\la_0\notin \Omega_{\infty}$. By i), iii) and iv) $\varphi\vert_{\cal L}$ coincides with the continuous function 
$\psi\vert_{\cal L}$ on $Supp\;\Delta \varphi\vert_{\cal L}$ which, by the continuity principle, implies that
$\varphi\vert _{\cal L}$ is continuous. Let $\Omega_0$ be the (bounded) component of ${\cal L}\setminus\big({\cal L}\cap Supp\;(dd^{c} \psi)\big)$ containing $\la_0$.
The continuous function $\big(\varphi-\psi\big)\vert_{\cal L}$ vanishes on $b\Omega_0$ (see iv)), is harmonic on $\Omega_0$
(see iii)) and negative (see ii)). The maximum principle now implies that $\varphi(\la_0)=\psi(\la_0)$.\qed\\

It is worth emphasize that quite precise informations about the behaviour of the bifurcation locus at infinity are required to apply the above Lemma to our problems.

\section{Distribution of critically periodic parameters in polynomial families}\label{SecCPP}

The aim of this section is to present a result due to Dujardin and Favre (see \cite{DF}) concerning the asymptotic distribution of degree $d$
polynomials which have a pre-periodic critical point. We will work in the context of
polynomial families, this will allow us to modify the original proof and significantly simplify it. We refer to the paper of Dujardin-Favre for results dealing with general holomorphic families of rational maps.\\

We work here in the family $\big(P_{c,a}\big)_{(c,a)\in \tiny{\CC^{d-1}}}$ of degree $d$ polynomials which has been introduced in the subsection \ref{ssPd}.
Let us recall that  $P_{c,a}$ is the polynomial of degree $d$ whose critical set is $\{0=c_0,c_1,\cdot\cdot\cdot,c_{d-2}\}$
and such that $P_{c,a}(0)=a^d$.\\

For $0\le i\le d-2$  and $0\le k<n$,  we denote by $\Per (i,n,k)$ the hypersurface of  $\CC^{d-1}$ defined by
\begin{center}
$\Per (i,n,k):=\{(c,a) \in \CC^{d-1}\;/\; P_{c,a}^{n}(c_i)= P_{c,a}^{k}(c_i)\}$.
\end{center}

The result we want to establish, which  has been first proved by Dujardin and Favre in \cite{DF}, is the following.

\begin{theo}\label{TheoDF}
In the family of degree $d$ polynomials, for any sequence of integers $(k_n)_n$ such that $0\le k_n<n$ one has
$\sum_{i=0}^{d-2}\lim_n d^{-n}[\Per (i,n,k_n)] =\Bc$.
\end{theo}

To simplify, we shall write $\la$ the parameters $(c,a)\in\CC^{d-1}$.
We follow the strategy described in the previous subsection.\\

The bifurcation current $\Bc$ is given by $\Bc=\sum_{i=0}^{d-2} dd^{c} g_\la(c_i)$ (see Lemma \ref{lemPotAC}) 
where $g_{\la}$ is the Green function of $P_\la$ (see the subsection \ref{condiscon}). It thus suffices to show that for any fixed $0\le i\le d-2$
 the following sequence of potentials
\begin{center}
$h_n(\la):=d^{-n}\ln \vert P_{\la}^{n}(c_i) - P_{\la}^{k_n}(c_i)\vert$
\end{center}
converges in $L_{loc}^1$ to $g_\la(c_i)$.

To this purpose we shall compare these potentials with the functions

\begin{center}
$g_n(\la):=d^{-n}\ln \max\big(1, \vert P_{\la}^{n}(c_i) \vert\big)$
\end{center}

which do  converge locally uniformly to $g_\la(c_i)$.\\

The first point is to check that the sequence $(h_n)_n$ is locally uniformly bounded from above. We shall actually prove a little bit more.
\begin{lem}\label{lemhnlb}
For any compact  $K\subset \CC^{d-1}$ and any $\epsilon>0$ there exists an integer $n_0$ such that $h_n\vert_K \le g_n\vert_K +\epsilon$
for $n\ge n_0$.
\end{lem}

\proof 
It is not difficult to see that there exists $R\ge 1$ such that
\begin{eqnarray}\label{hnlb1}
(1-\epsilon)\vert z\vert ^{d^{n}}\le \vert P_{\la}^{n}(z) \vert \le (1+\epsilon) \vert z\vert ^{d^{n}}
\end{eqnarray}
for every $\la \in K$, every $n\in \NN$ and every $\vert z\vert \ge R$.

We now proceed by contradiction and assume that there exists $\la_p\in K$ and $n_p\to +\infty$ such that 
$h_{n_p}(\la_p) \ge  g_{n_p}(\la_p) +\epsilon$. This means that
\begin{eqnarray}\label{hnlb2}
\vert P_{\la_p}^{n_p}(c_{i}) - P_{\la_p}^{k_{n_p}}(c_i)\vert \ge e^{\epsilon d^{n_p}}
 \max\big(1, \vert P_{\la_p}^{{n_p}}(c_{i}) \vert\big).
\end{eqnarray}
Let us set $B_p:=P_{\la_p}^{k_{n_p}}(c_{i})$. By \ref{hnlb2} we have $\lim_p \vert B_p\vert=+\infty$ and thus $\vert B_p\vert \ge R$ for $p$ big enough. Then, using \ref{hnlb1}, one may write 
\begin{center}
$P_{\la_p}^{n_p}(c_{i})=P_{\la_p}^{n_p-k_{n_p}}(B_p)=(u_pB_p)^{d^{n_p-k_{n_p}}}$
\end{center}
where $(1-\epsilon)\le \vert u_p\vert \le (1+\epsilon)$
and the estimate \ref{hnlb2} becomes
\begin{center}
$\vert (u_pB_p)^{d^{n_p-k_{n_p}}} -B_p\vert \ge e^{\epsilon d^{n_p}}
\vert u_pB_p\vert^{d^{n_p-k_{n_p}}}$.
\end{center}
This is clearly impossible when $p\to +\infty$.\qed\\

We now have to check that $(h_n)_n$ does not converge to $-\infty$.
The following technical Lemma deals with that and will also play an important role in the remaining of the proof.

\begin{lem}\label{cinbassatt}
If $c_i$ belongs to some attracting basin of $P_{\la_0}$ then there exists a neighbourhood $V_0$ of $\la_0$ such that
$\sup_n \sup_{V_0} \big(h_n - g_n\big)\ge 0$.
\end{lem}

\proof If $V_0$ is a sufficently small neighbourhood of $\la_0$ then $P_\la^{n}(c_i)\to a_\la$ where $a_\la$ is an attracting cycle of $P_\la$
for every $\la\in V_0$. We will assume that $P_\la(a_\la)=a_\la$.\\
Let us now proceed by contradiction and suppose that there exists $\epsilon >0$ such that
\begin{eqnarray}\label{cin1}
\vert P_\la^{n} (c_i) - P_\la^{k_n} (c_i)\vert \le e^{-\epsilon d^{n}} \max\big(1, \vert P_{\la}^{n}(c_i) \vert\big),\;\;\;\forall \la\in V_0,\;\forall n.
\end{eqnarray}
Since   $P_\la^{n}(c_i)\to a_\la$, \ref{cin1} would imply
\begin{eqnarray}\label{cin2}
\vert P_\la^{n} (c_i) - P_\la^{k_n} (c_i)\vert \le C e^{-\epsilon d^{n}},\;\;\;\forall \la\in V_0,\;\forall n.
\end{eqnarray}
The estimate \ref{cin2} implies that $a_\la$ is a super-attracting  fixed point for any $\la\in V_0$ which, in turn,
implies that $a_\la =\infty$ for all $\la\in V_0$. But in that case we would have $\vert P_\la^{k_n} (c_i)\vert \le \frac{1}{2}
\vert P_\la^{n} (c_i)\vert$ for $n$ big enough and the estimate \ref{cin1} would be violated.\qed\\

We finally have to show that $g:=g_\la(c_i)=\lim g_n$ is the only limit value of the sequence
$(h_n)_n$ for the $L_{loc}^1$ convergence.
Assume that (after taking a subsequence!) $h_n$ is converging in $L_{loc}^1$ to $h$. To prove that the functions $h$ and $g$ coincide we shall
check that they satisfy the assumptions of Lemma \ref{lemPmax}.\\

Our modification of the original proof essentially stays in the third step. \\

\underline{First step}: $h\le g$.\\

Let $B_0$ be a ball of radius $r$ centered at $\la_0$ and let $\epsilon>0$. By the mean value property we have
\begin{center}
$h(\la_0)\le \frac{1}{\vert B_0\vert}\int_{B_0} h =\lim_n \frac{1}{\vert B_0\vert}\int_{B_0} h_n$
\end{center}
but, according to Lemma \ref{lemhnlb}, $h_n\le g_n +\epsilon$ on $B_0$
for $n$ big enough and thus
\begin{center}
$h(\la_0)\le \epsilon + \lim_n \frac{1}{\vert B_0\vert}\int_{B_0} g_n=\epsilon + \frac{1}{\vert B_0\vert}\int_{B_0} g$.
\end{center}
As $g$ is continuous, the conclusion follows by making $r\to 0$ and then $\epsilon \to 0$. \\

\underline{Second step}: $h=g$ on $Supp\;dd^{c}g$.\\

Combining Lemma \ref{cinbassatt} and the result of step one we will first establish the following 
\begin{center}
{\bf Fact} $(h-g)$ vanishes when $c_i$ is captured by an attracting basin.
\end{center} 
Suppose to the contrary that $c_i$ is captured by an attracting cycle of $P_{\la_0}$ and $(h-g)(\la_0) <0$. As the function $(h-g)$
is upper semi-continuous, we may shrink $V_0$ so that $(h-g)\le -\epsilon<0$ on $V_0$. Now, as $c_i$ is passive on $V_0$,
the function $(h-g)$ is $p.s.h$ on $V_0$ and, after shrinking $V_0$ again, Hartogs Lemma \ref{lemHar} implies that $(h_n-g_n) \le -\frac{\epsilon}{2}$ on $V_0$ for $n$ big enough. This contradicts Lemma \ref{cinbassatt}.\\

Now, if $\la_0\in Supp\;dd^{c}g$ then $\la_0=\lim_k\la_k$
where  $\la_k$ is a parameter for which $c_i$ is captured by some attracting cycle (see Lemma \ref{lemappperno}).
As $(h-g)$ is upper semi-continuous we get $(h-g)(\la_0)\ge\limsup (h-g)(\la_k)=0$ and, by the first step,  $(h-g)(\la_0)=0$.\\

\underline{Third step}: $Supp\;dd^{c}h\subset Supp\;dd^{c}g$.\\

Let $\Omega$ be a connected component of $\CC^{d-1}\setminus Supp\;dd^{c}g$. We have to show that $h$ is pluriharmonic on $\Omega$. We proceed by contradiction. If $dd^{c} h$ does not vanish on $\Omega$ then there exists some $n_0$ for which some irreducible component ${\cal H}$ of
\begin{center}
$\{P_\la^{n_0}(c_i) -P_\la^{k_{n_0}}(c_i)=0\}$
\end{center}
meets $\Omega$. 
When $\la\in {\cal H}$ then $c_i$ is captured by a cycle since $P_\la^{k_{n_0}}(c_i)=:z(\la)$ satisfies
$P_\la^{m_0}(z(\la))= P_\la^{m_0}\circ P_\la^{k_{n_0}}(c_i)= P_\la^{k_{n_0}}(c_i)= z(\la)$ for $m_0:=n_0-k_{n_0}>0$.\\

Let us show that $z(\la)$ is  a neutral periodic point. We first observe that the vanishing of  $dd^{c} g$ on $\Omega$ forces $z(\la)$ to be non-repelling and thus $\vert \big(P_\la^{m_0}\big)'\big(z(\la)\big)\vert \le 1$
on ${\cal H}\cap \Omega$.\\
 Let us now see why  $z(\la)$ cannot be attracting.
If this would be the case then, by the above Fact, we would have $h(\la_0)=g(\la_0)$ for a certain $\la_0\in {\cal H}\cap \Omega$.
As $(h-g)$ is negative and $p.s.h$ on $\Omega$ this implies, via the maximum principle, that $h=g$ on $\Omega$.
This is imposible since $dd^{c}h$ is supposed to be non vanishing on $\Omega$.\\
We thus have $ \big(P_\la^{m_0}\big)'\big(z(\la)\big)=e^{i\nu_0}$
on ${\cal H}\cap \Omega$
and therefore  $z(\la)$ belongs to a neutral cycle whose period $p_0$ divides $m_0$ and whose multiplier is a $q_0$-root of $e^{i\nu_0}$ where $m_0=p_0 q_0$.
In other words, ${\cal H}\cap \Omega$ is contained in a finite union of hypersurfaces of the form $\Per_{n}(e^{i\theta}) $.
This implies that
\begin{center}
${\cal H} \subset \Per_{n_0}(e^{i\theta_0}) $.
\end{center}
for some integer $n_0$ and some
real number $\theta_0$.\\

Finally, using a global argument, we will see that this is impossible.
Let us recall the following dynamical fact.

\begin{lem}
Every polynomial which has a neutral cycle also has a bounded, non-preperiodic, critical orbit.
\end{lem}

Thus, when $\la\in{\cal H}$, the polynomial $P_\la$ has two distinct bounded critical orbits; the orbit of $c_i$ which is preperiodic and the orbit of some other critical point which is given by the above Lemma. This shows that ${\cal H}$ cannot meet the line $\{c_0=c_1=\cdot\cdot\cdot=c_{d-2}=0\}:=M_d$ since the corresponding polynomials (which are given by $\frac{1}{d}z^{d}+a^{d}$ were $a\in\CC$)
have only one critical orbit. 
We will now work in the projective compactification of $\CC^{d-1}$ introduced in subsection \ref{ssPd}. By Theorem \ref{controlinfty},
${\cal H}$ and $M_d$ cannot meet at infinity. This contradicts  Bezout's theorem.\\

\underline{Fourth step}: for any $\la_0\in \CC^{d-1}$ there exists a complex line ${\cal L}$ through $\la_0$ such that $h=g$ on 
 the unbounded component of ${\cal L}\setminus\big({\cal L}\cap Supp (dd^{c} g)\big)$.\\
 
 Here one uses again Theorem \ref{controlinfty} to pick a line ${\cal L}$ through $\la_0$ which meets infinity at some point 
 $\xi_0\notin \overline{Supp\;dd^{c} g}$. Then, for any $\la$ in the unbounded component of 
 ${\cal L}\setminus\big({\cal L}\cap Supp\;(dd^{c} g)\big)$ the critical point $c_i$ belongs to the super-attracting basin of $\infty$ and thus,
 as we saw in second step, $h(\la)=g(\la)$.\qed\\

The general result obtained by Dujardin and Favre may be stated as follows.

\begin{theo}\label{DFGen}
Let $\left(f_\la\right)_M$ be a holomorphic family of degree $d$ rational maps with a marked critical point $c_\la$ which is not stably preperiodic.
Let $Per(n,k)$ be the hypersurface in $M$ defined by:
\begin{center}
$\Per (n,k):=\{\la\in M\;/\; f_\la^{n}(c_\la)= f_\la^{k}(c_\la)\}$.
\end{center}
 Assume that the following assumption is satisfied:
 \begin{center}
 $(H)$ For every $\la_0\in M$ there exists a curve $\Gamma\subset M$ passing through $\la_0$ such that $\{\la\in\Gamma\;/\;c_\la \textrm{ is attracted by a cycle}\}$ has a relatively compact complement in $\Gamma$.
\end{center}
Then, for any sequence of integers $(k_n)_n$ such that $0\le k_n<n$ one has
\begin{center}
$\lim_n d^{-n-(1-e)k_n}[\Per (n,k_n)] =T_c$
\end{center}
where $T_c$ is the activity current of the marked critical point $c_\la$ and $e$ the cardinal of the exceptional set of $f_\la$ for a generic $\la$.
\end{theo}

It would be interesting to remove the assumption $(H)$ which seems to be a technical one. Observe however that the above Theorem covers the case of the moduli space $Mod_2$.

\section{Distribution of rational maps with  cycles of a given multiplier}

Let $f:M\times\pp\to\pp$  be an arbitrary holomorphic family of degree $d\ge 2$ rational maps.\\
 
We want to  investigate the asymptotic distribution of the hypersurfaces $Per_n(w)$ in $M$ when $\vert w\vert <1$. 
Concretely, we will consider the current of integration $[\Per_n(w)]$ or, more precisely, the currents
\begin{center} 
$[\Per_n(w)] :=dd^c\;\ln \vert p_n(\la,w)\vert$
\end{center}
 where $p_n(\cdot,w)$  are the canonical defining functions for  the hypersurfaces $\Per_n(w)$ constructed in section \ref{secPer} by mean of dynatomic polynomials.
We ask if the following convergence occurs:
 \begin{center}
  $\lim_n\frac{1}{d^n}[\Per_n(w)]=\Bc$.
 \end{center}
The question is easy to handle when $\vert w\vert <1$, more delicate when $\vert w\vert =1$ and widely open when
$\vert w \vert >1$.
 
\subsection{The case of attracting cycles}

We aim to prove the following general result (see \cite{BB2}). 

\begin{theo}\label{theodisatt}
For any holomorphic family of degree $d$ rational maps $\big(f_{\la}\big)_{\la \in M}$ one has
$d^{-n}\;[\Per_n(w)]\to \Bc$ when $\vert w\vert <1$.
\end{theo}

Let us have a look to  the case  $w=0$.  Comparing $Per_n(0)$ with the hypersurfaces $Per(n,0)$ considered in the
the last section, one sees that the above result may be derived from the Theorem \ref{DFGen} of Dujardin and Favre but that we
do not need any special assumption. Moreover, for the quadratic polynomial family one obtains the equidistribution of  centers of  hyperbolic components of the Mandelbrot set.
This was first proved by Levin \cite{Lev}.\\
Using arithmetical methods,  Favre and Rivera-Letelier \cite {FRL} have estimated the equidistribution speed of 
$d^{-n}\;[\Per_n(0)]$ for unicritical families $(z^d+c)$.\\

\proof
Let us set 
\begin{center}
$\displaystyle L_n(\la,w):=d^{-n}\ln \vert p_n(\la,w)\vert$.
\end{center}
Since, by definition, $\Bc=dd^c\;L(\la)$ 
where $L(\la)$ be the Lyapunov exponent of $(\pp,f_{\la},\mu_{\la})$ and $\mu_{\la}$ is the Green measure of 
$f_{\la}$, all we have to show is that $L_n$ converges to $L$ in $L^1_{loc}(M)$. 
Here again we shall use the compactness principle for subharmonic functions (see Theorem \ref{compapsh}).\\

The situation is purely local and therefore, taking charts, we may assume that $M={\bf C}^k$.
We write the polynomials $p_n$ as follows :$$p_n(\la,w)=:\prod_{i=1}^{N_d(n)}\big(w-w_{n,j}(\la)\big).$$
Using the fact that
$d^{-n}N_d(n)\sim\frac{1}{n}$ (see Theorem \ref{theopoly})
one sees that the sequence $L_n$ is locally uniformly bounded from above.\\

 According to Theorem \ref{theopoly}, the set $\{w_{n,j}(\la)\; /\;w_{n,j}(\la)\ne 1\}$ coincides with the set of multipliers
of cycles of exact period $n$ (counted with multiplicity) from which the cycles of multiplier $1$ are deleted. We thus have

\begin{eqnarray}\label{A1}
\sum_{j=1}^{N_d(n)} \ln^+ \vert w_{n,j}(\la)\vert =\frac{1}{n} \sum_{p\in R^{*}_n(\la)} \ln \vert (f_{\la}^n)'(p)\vert
\end{eqnarray}

where
$R^{*}_n(\la):=\{p\in \pp \;/\; \;p\; \textrm{has exact period}\; n \;  \textrm{and}\;\vert (f_{\la}^{n })' (p)\vert > 1 \}$.
Since $f_\la$ has a finite number of non-repelling cycles (Fatou's theorem), one sees that there exists $n(\la) \in\NN$ such that
\begin{eqnarray}\label{A2}
n\ge n(\la) \Rightarrow \vert w_{n,j}(\la)\vert > 1,\;\textrm{for any}\;1\le j\le N_d(n).
\end{eqnarray}

By \ref{A1} and \ref{A2}, one gets
\begin{eqnarray*}
L_n(\la,0)=d^{-n}\sum_{j=1}^{N_d(n)} \ln \vert w_{n,j}(\la)\vert =d^{-n}\sum_{j=1}^{N_d(n)} \ln^+ \vert w_{n,j}(\la)\vert =\frac{d^{-n}}{n} \sum_{R^{*}_n(\la)} \ln \vert (f_{\la}^n)'(p)\vert
\end{eqnarray*}

 for $n\ge n(\la)$ which, by Theorem \ref{theoapprox}, yields:

\begin{eqnarray}\label{A3}
\lim_n L_n(\la,0)=L(\la),\; \forall \la \in M.
\end{eqnarray}

If now $\vert w\vert <1$, it follows from \ref{A2} that
 $L_n(\la,w)-L_n(\la,0)=d^{-n}\sum_j\ln \frac{\vert w_{n,j}(\la) - w\vert }{\vert w_{n,j}(\la) \vert}$ 
 and 
$\ln(1-\vert w\vert)\le\ln\frac{\vert w_{n,j}(\la) - w\vert }{\vert w_{n,j}(\la) \vert}\le  \ln(1+\vert w\vert)$ for $1\le j\le N_d(n)$ and $n\ge n(\la)$ .
We thus get 
\begin{eqnarray*}
d^{-n}N_d(n)\ln(1-\vert w\vert)\le\vert L_n(\la,w)-L_n(\la,0) \vert \le d^{-n}N_d(n) \ln(1+\vert w\vert)
\end{eqnarray*}
for $n\ge n(\la)$. Using \ref{A3} and the fact that  $d^{-n}N_d(n)\sim\frac{1}{n}$  we obtain 
$\lim_n L_n(\la,w)=L(\la)$ for any $(\la,w) \in M\times {\Delta}$.
The $L_{loc}^1$ convergence of $L_n(\cdot,w)$ now follows immediately from Theorem \ref{compapsh}.\qed\\

\begin{rem}\label{ptLn}
We have proved that $\displaystyle L_n(\la,w):=d^{-n}\ln \vert p_n(\la,w)\vert$ converges pointwise to $L(\la)$ on $M$ when $\vert w\vert<1$.
\end{rem}

The above discussion shows that
the pointwise convergence of $L_n(\la,w)$ to $L$ (and therefore the convergence $d^{-n}[\Per_n(w)]\to\Bc$)  is quite a straightforward consequence of Theorem \ref{theoapprox} when $\vert w\vert <1$.
However, when $\vert w\vert \ge 1$ and $\la$ is a non-hyperbolic parameter, the control of $L_n(\la,w)=d^{-n}\sum\ln\vert w- w_{n,j}(\la)\vert$ is very delicate because $f_\la$ may have many cycles whose multipliers are close to $w$.
This is why we introduce the $p.s.h$ functions $ L_n^+$ which both coincide with $L_n$ on the hyperbolic components and are rather easily seen to converge nicely. These functions will be extremely helpful later.\\

\begin{defi}
 The $p.s.h$ functions $ L_n^+$ are defined by:
\begin{center}
$\displaystyle L_n^+(\la,w):=d^{-n}\sum_{j=1}^{N_d(n)} \ln^+\vert w- w_{n,j} (\la)\vert$
\end{center}
where  $p_n(\la,w)=:\prod_{j=1}^{N_d(n)}(w-w_{n,j}(\la))$ are the polynomials associated to the family $\left(f_\la\right)_{\la\in M}$
by Theorem \ref{theopoly}.
\end{defi}

The interest of considering these functions stays in the next Lemma.

\begin{lem}\label{lemLn+}
 The sequence $L^+_n$ converges pointwise and in $L_{loc}^1$ to $L$ on $M\times{\bf C}$. For every $w\in {\bf C}$ the sequence $L^+_n(\cdot,w)$ converges in  $L_{loc}^1$ to $L$ on $M$. 
\end{lem}

\proof We will show that $L_n^+(\cdot,w)$ converges pointwise to $L$ on $M$ for every $w\in{\bf C}$. As
$(L_n^+)_n$ is locally uniformly bounded, this implies the convergence of $L_n^+(\cdot,w)$ in $L_{loc}^1(M)$ (Theorem \ref{compapsh}) and the convergence of $L_n^+$ in 
$L_{loc}^1(M\times{\bf C})$ then follows by Lebesgue's theorem.\\

As $L_n(\la,0)\to L(\la)$ (see Remark \ref{ptLn}),
we have to estimate $L^+_n(\la,w)-L_n(\la,0)=:\epsilon_n(\la,w)$ on $M$. Let us fix $\la \in M$, $w\in {\bf C}$ and pick $R>\vert w\vert $. Since $f_\la$ has a finite number of non-repelling cycles (Fatou's theorem), one sees that there exists $n(\la) \in\NN$ such that
\begin{eqnarray*}
n\ge n(\la) \Rightarrow \vert w_{n,j}(\la)\vert > 1,\;\textrm{for any}\;1\le j\le N_d(n).
\end{eqnarray*}
We may then decompose $\epsilon_n(\la,w)$ in the following way:

\begin{eqnarray*}
\epsilon_n(\la,w)=d^{-n}\sum_{1\le\vert w_{n,j}(\la)\vert<R+1} \ln^+\vert w_{n,j}(\la) - w\vert 
+ d^{-n}\sum_{\vert w_{n,j}(\la)\vert\ge R+1} \ln\frac{\vert w_{n,j}(\la) - w\vert }{\vert  w_{n,j}(\la)\vert }\\
-d^{-n}\sum_{1\le\vert w_{n,j}(\la)\vert<R+1} \ln \vert w_{n,j}(\la) \vert .
\end{eqnarray*}

We may write this decomposition as
$
\epsilon_n(\la,w)=:\epsilon_{n,1}(\la,w) + \epsilon_{n,2}(\la,w) - \epsilon_{n,1}(\la,0).
$\\

Clearly,
$
0\le \epsilon_{n,1}(\la,w)\le d^{-n} N_d(n) \ln\big(2R+1\big)
$
and thus
$\lim_n  \epsilon_{n,1}(\la,w)=0$.
Similarly, 
$
\lim_n \epsilon_{n,2}(\la,w)=0
$ follows from the fact that one has:

\begin{eqnarray*}
\ln(1-\frac{R}{R+1})\le\ln\frac{\vert w_{n,j}(\la) \vert -R}{\vert w_{n,j}(\la) \vert}\le\ln\frac{\vert w_{n,j}(\la) - w\vert }{\vert w_{n,j}(\la) \vert}\\
\le\ln\frac{\vert w_{n,j}(\la) \vert +R}{\vert w_{n,j}(\la) \vert}\le \ln(1+\frac{R}{R+1}).
\end{eqnarray*}
for $\vert w_{n,j}(\la)\vert >R+1>\vert w\vert +1$.
\qed\\

As the functions $L_n^+$ and $L_n$ coincide on hyperbolic components, the above Lemma would easily yield the convergence of 
$d^{-n}[\Per_n(w)]$ towards $\Bc$ \emph{for any} $w\in {\bf C}$ if the density of hyperbolic parameters in $M$ was known.
The remaining of this section is, in some sense, devoted to overcome this difficulty.
We shall first do this in a general setting by averaging the multipliers. Then we will restrict ourself to polynomial families and, using the nice distribution of hyperbolic parameters near infinity, will show that $d^{-n}[\Per_n(e^{i\theta})]$ converges towards $\Bc$.

\subsection{Averaging the multipliers}

Although the convergence of $\lim_n\frac{1}{d^n}[\Per_n(w)]$ to $\Bc$ is not clear when $\vert w\vert \ge 1$, one easily obtains
the convergence by averaging over the argument of the multiplier $w$. The following result is due to Bassanelli and the author
\cite{BB2}.

\begin{theo}\label{theodistaver}
For any holomorphic family of degree $d$ rational maps $\big(f_{\la}\big)_{\la \in M}$ one has
$\frac{d^{-n}}{2\pi}\int_0^{2\pi} [\Per_n(re^{i\theta})]\;d\theta \to \Bc,\;\textrm{when}\;
r\ge 0$.
\end{theo}

\proof

One essentially has to investigate the following sequences of $p.s.h$ functions
\begin{center} 
$\displaystyle L_n^r(\la):=\frac{d^{-n}}{2\pi}\int_0^{2\pi} \ln \vert p_n(\la,re^{i\theta})\vert\;d\theta$.
\end{center}
We will see that $L_n^r(\la)\ge C\frac{\ln r}{n}$ where $C$ only depends on the family and that $L_n^r$ converges to $L$ in $L^1_{loc} (M)$.\\

For that, we essentially will compare  $ L_n^r $ with $L_n(\la,0)=L_n^0 $ by 
using the formula $\ln \max(\vert a \vert,r)=\frac{1}{2\pi}\int_0^{2\pi}\ln \vert a-re^{i\theta}\vert d\theta$.
Indeed, writting $$p_n(\la,w)=:\prod_{i=1}^{N_d(n)}\big(w-w_{n,j}(\la)\big)$$ this formula yields

\begin{eqnarray}\label{Comp}
L_n^r(\la)=\frac{1}{2\pi d^n}\int_0^{2\pi} \ln \prod_j\vert re^{i\theta}-w_{n,j}(\la)\vert d\theta=
d^{-n}\sum_j \ln \max(\vert w_{n,j}(\la)\vert,r).
\end{eqnarray}

 According to Theorem \ref{theopoly}, the set $\{w_{n,j}(\la)\; /\;w_{n,j}(\la)\ne 1\}$ coincides with the set of multipliers
of cycles of exact period $n$ (counted with multiplicity) from which the cycles of multiplier $1$ are deleted.
Using the fact that
$d^{-n}N_d(n)\sim\frac{1}{n}$ (see Theorem \ref{theopoly}) one sees that the sequence $L_n^r(\la)$ is locally bounded from above
and is uniformly bounded from below by  $C\frac{\ln r}{n}$.
Since $f_\la$ has a finite number of non-repelling cycles (Fatou's theorem),  there exists $n(\la) \in\NN$ such that
\begin{eqnarray*}
n\ge n(\la) \Rightarrow \vert w_{n,j}(\la)\vert > 1,\;\textrm{for any}\;1\le j\le N_d(n).
\end{eqnarray*}
Now we deduce from \ref{Comp} that for $n\ge n(\la)$:

\begin{eqnarray*}
L_n^r(\la)=d^{-n}\sum_j \ln \vert w_{n,j}(\la)\vert +d^{-n}\sum_{1\le\vert w_{n,j}(\la)\vert < r} \ln \frac{r}{\vert w_{n,j}(\la)\vert }=\\
L_n(\la,0) + d^{-n}\sum_{1\le\vert w_{n,j}(\la)\vert < r}  \ln \frac{r}{\vert w_{n,j}(\la)\vert }
\end{eqnarray*}

and thus

\begin{eqnarray*}
0\le L_n^r(\la) -
L_n(\la,0)= d^{-n}\sum_{1\le\vert w_{n,j}(\la)\vert < r}  \ln \frac{r}{\vert w_{n,j}(\la)\vert }\le
d^{-n}N_d(n)\ln^+ r.
\end{eqnarray*}

Recalling that $d^{-n}N_d(n)\sim\frac{1}{n}$ and $L_n(\la, 0)\to L(\la)$ (see Remark \ref{ptLn}), this implies that $L_n^r$ converges pointwise  to $L$ and, by Theorem \ref{compapsh}, 
 that $(L_n^r)_n$ converges to $L$ in $L^1_{loc} (M)$.\\
 
Now, to get the conclusion, one has to justify the following identity:
\begin{center}
$dd^c L_n^r=\frac{d^{-n}}{2\pi}\int_{0}^{2\pi}[Per_n(re^{i\theta})]d\theta.$
\end{center}

Going back to definitions and taking $dd^c \ln\vert p_n(\la,re^{i\theta})\vert = [Per_n(re^{i\theta})]$ into account, one sees that this is a conseqence of Fubini's theorem if
one checks that $\ln\vert p_n(\la,re^{i\theta})\vert$ is locally integrable.
Let $K$ be a compact subset of $M$ and $c_n$ be an upper bound for $\ln\vert p_n(\la,r e^{i\theta})\vert$ on $K\times[0,2\pi]$.
Then, the negative function   
$\ln\vert p_n(\la,e^{i\theta})\vert -c_n$
is indeed integrable on $K\times[0,2\pi]$
as it follows from the fact that $L_n^r(\la)\ge C\frac{\ln r}{n}$:
\begin{eqnarray*}
\int_{K}\big(\int_{0}^{2\pi}\big( \ln\vert p_n(\la,r e^{i\theta})\vert -c_n\big) d\theta\big)dV&=& 2\pi d^n\int_{K}
L_n^r dV - 2\pi c_n \int_{K}dV\\
&\ge& \left(d^n C\frac{\ln r}{n} -  c_n\right) 2\pi \int_{K}dV.\\
\end{eqnarray*}
\qed\\

\begin{rem}\label{ptLnr}
We have proved that $\displaystyle L_n^r(\la):=\frac{d^{-n}}{2\pi}\int_0^{2\pi} \ln \vert p_n(\la,re^{i\theta})\vert\;d\theta$ is pointwise converging to $L(\la)$ on $M$.
\end{rem}

The following result is essentially a potential-theoretic consequence of the former one. It implicitely contains some information about the convergence of $d^{-n}\;[\Per_n(w)]$ for arbitrary choices of $w$ but seems hard to improve without furtherly use dynamical properties (see \cite{BB3}). 

\begin{theo}\label{theogene}
For any family of degree $d$ rational maps $\big(f_{\la}\big)_{\la \in M}$ one has 
\begin{center}
$d^{-n}\;dd^c_{(\la,w)} \ln \vert p_n(\la,w)\vert\to dd^c L(\la)$
\end{center}
where $p_n(\cdot,w)$  are the canonical defining functions for  the hypersurfaces $\Per_n(w)$ given by Theorem \ref{theopoly}.
\end{theo}

\proof
Let us set 
\begin{center}
$ L_n(\la,w):=d^{-n}\ln \vert p_n(\la,w)\vert$.
\end{center}

As we have seen in the two last subsections (see remarks \ref{ptLn} and \ref{ptLnr})
\begin{center}
$L_n(\la,0)\to L(\la)$\\
$ L_n^r(\la):=\frac{d^{-n}}{2\pi}\int_0^{2\pi} \ln \vert p_n(\la,re^{i\theta})\vert\;d\theta \to L(\la)$ for any $r\ge 0$.
\end{center}

Let us also recall that the function $L$ is continuous on $M$ (see Theorem \ref{theoLH}).\\

As
 the functions $L_n$ are $p.s.h$ and the sequence $(L_n)_n$ is locally uniformly bounded from above, we shall again
use the compacity properties of $p.s.h$ functions given by Theorem \ref{compapsh}. Since $L_n(\la,0)$ converges to $L(\la)$, the sequence
$(L_n)_n$ does not converge to $-\infty$ and it therefore suffices to show that,  among $p.s.h$ functions on $M\times {\bf C}$, the function $L$ is the only possible
limit for  $(L_n)_n$ in $L^1_{loc} (M\times {\bf C})$.

Let $\varphi$ be a $p.s.h$ function on $M\times {\bf C}$ and $(L_{n_j})_j$ a subsequence of $(L_n)_n$  which converges to $\varphi$
in $L^1_{loc} (M\times {\bf C})$. Pick $(\la_0,w_0)\in M\times{\bf C }$. We have to prove that $\varphi(\la_0,w_0)= L(\la_0)$.\\

Let us first observe that $\varphi(\la_0,w_0)\le L(\la_0)$.
Take a ball $B_{\epsilon}$ of radius $\epsilon$ and centered at $(\la_0,w_0)\in M\times{\bf C }$. By the submean value property and the $L_{loc}^1$-
convergence
of $L_n^+$ (see Lemma \ref{lemLn+}) we have:

\begin{eqnarray*}
\varphi(\la_0,w_0)\le \frac{1}{\vert B_{\epsilon}\vert}\int_{B_{\epsilon}}\varphi\; dm=\lim_j \frac{1}{\vert B_{\epsilon}\vert}\int_{B_{\epsilon}} L_{n_j}\;dm\\
\le \lim_j \frac{1}{\vert B_{\epsilon}\vert}\int_{B_{\epsilon}} L_{n_j}^+\;dm = \frac{1}{\vert B_{\epsilon}\vert}\int_{B_{\epsilon}} L\;dm
\end{eqnarray*}

making then $\epsilon\to 0$, one obtains $\varphi(\la_0,w_0)\le L(\la_0)$ since $L$ is continuous.\\
Let us now check that $\limsup_j L_{n_j}(\la_0,w_0e^{i\theta})=L(\la_0)$ for almost all $\theta\in [0,2\pi]$. Let $r_0:=\vert w_0\vert$. 
By Lemma \ref{lemLn+}, the sequence $L_n^+$ converges pointwise  to $L$ and therefore:
$$\limsup_j L_{n_j}(\la_0,w_0e^{i\theta})\le \limsup_j L_{n_j}^+(\la_0,w_0e^{i\theta})=L(\la_0).$$

On the other hand, by pointwise convergence of $L_{n}^{r_0}$ to $L$ and Fatou's lemma we have:

\begin{eqnarray*}
L(\la_0)=\lim_nL_{n}^{r_0} (\la_0)=\limsup_j\frac{1}{2\pi}\int_0^{2\pi} L_{n_j}(\la_0,r_0 e^{i\theta}) d\theta \le\\
\frac{1}{2\pi}\int_0^{2\pi} \limsup_j L_{n_j}(\la_0,r_0 e^{i\theta}) d\theta
\end{eqnarray*}

and the desired property follows immediately.\\
To end the proof we argue by contradiction and assume that $\varphi(\la_0,w_0) < L(\la_0)$.
As $\varphi$ is upper semi-continuous and $L$ continuous, there exists a  neighbourhood $V_0$ of 
$(\la_0,w_0)$ and $\epsilon>0$ such that 
$$\varphi - L \le -\epsilon\;\;\textrm{on}\;V_0.$$
Pick a small ball $B_{\la_0}$ centered at $\la_0$ and a small disc $\Delta_{w_0}$ centered at $w_0$ such that  
$B_0:=B_{\la_0}\times\Delta_{w_0}$ is relatively compact in $V_0$.  Then, according to Hartogs Lemma \ref{lemHar},
we have:

$$\limsup_j \big(\sup_{B_0} (L_{n_j}-L)\big)\le \sup_{B_0} (\varphi -L)\le -\epsilon.$$
This is impossible since, as we have seen before, we may find $(\la_0,r_0 e^{i\theta_0})\in B_0$ such that
$\limsup_j\big(L_{n_j}(\la_0,r_0e^{i\theta_0})-L(\la_0)\big)=0$.\qed\\

\begin{rem} Using standard techniques, one may deduce from the above Theorem that the set of multipliers $w$ for which
the bifurcation current $\Bc$ is not a limit of the sequence $d^{-n} [\Per_n(w)]$ is contained in a polar subset of the complex plane.
\end{rem}

\subsection{The case of neutral cycles in polynomial families}\label{CNC}

We return to the family $\big(P_{c,a}\big)_{(c,a)\in \tiny{\CC^{d-1}}}$ of degree $d$ polynomials.
Recall that  $P_{c,a}$ is the polynomial of degree $d$ whose critical set is $\{0=c_0,c_1,\cdot\cdot\cdot,c_{d-2}\}$
and such that $P_{c,a}(0)=a^d$ (see subsection \ref{ssPd}).\\

We want to prove that, in this family,
$\lim_n d^{-n}[\Per_n(w)]=\Bc$ for $\vert w\vert\le 1$.
Taking the results of the previous subsection into account (see Theorem \ref{theodisatt}), it remains to treat the case $\vert w\vert=1$
and prove the following result due to Bassanelli and the author (see \cite{BB3}).

\begin{theo}\label{theoequipol} In the family of degree $d$ polynomials 
 $\lim_n d^{-n}[\Per_n(e^{i\theta})]=\Bc$ for any $\theta\in [0,2\pi]$.\\
\end{theo}

We will follow the strategy described in subsection \ref{Sgs}.
As we shall see, the proof would be rather simple if we would know that the bifurcation locus is accumulated by hyperbolic parameters.
This is however unknown when $d\ge 3$ and is a source of technical difficulties (see the fourth step).\\

\proof
We denote by $\la$ the parameter in ${\bf C}^{d-1}$ (i.e. $\la:=(c,a)$) and set 
$$L_n(\la):=d^{-n}\ln \vert p_n(\la,e^{i\theta})\vert$$ where the polynomials $p_n(\la,w)$ are those given by Theorem \ref{theopoly}. 
We have to show that the sequence $(L_n)_n$ converges to $L$ in $L_{loc}^1$.\\ 

We have already seen that $(L_n)_n$ is a  uniformly locally bounded sequence of $p.s.h$ functions on ${\bf C}^{d-1}$. 
 Since the family  $\{P_{c,a}\}_{ (c,a)\in{\bf C}^{d-1}}$ contains hyperbolic parameters, on which the $L_n(\la)=L_n^+(\la,e^{i\theta})$,
 it follows from Lemma \ref{lemLn+} that
 the sequence $(L_n)_n$ does not converge to $-\infty$.
Thus, according to Theorem \ref{compapsh}, we have to show that $L$ is the only limit value of the sequence
$(L_n)_n$ for the $L_{loc}^1$ convergence.\\

Assume that (after taking a subsequence!) $(L_n)_n$ is converging in $L_{loc}^1$ to $\varphi$. To prove that the $p.s.h$ functions $\varphi$ and $L$ coincide we shall
check that they satisfy the assumptions of Lemma \ref{lemPmax}.\\

\underline{First step:} $\varphi\le L$.\\

Since $L^+_n(\la,e^{i\theta})$ converges to $L$ in $L^1_{loc}$ (see Lemma \ref{lemLn+}) and $L_{n}(\la)\le \L^+_{n}(\la,e^{i\theta})$ we get
$$\varphi(\la_0)\le \frac{1}{\vert B_{\epsilon}\vert}\int_{B_{\epsilon}}\varphi\; dm \le \frac{1}{\vert B_{\epsilon}\vert}\int_{B_{\epsilon}} L\;dm$$ for any small ball $B_{\epsilon}$ centered at 
$\la_0$. The desired inequality then follows by making $\epsilon\to 0$ since the function $L$ is continuous (see Theorem \ref{theoLH}).\\

\underline{Second step}: $Supp\;dd^c \varphi\subset Supp\;dd^c L$.\\

Since there are no persistent neutral cycles in the family $\big(P_{c,a}\big)_{(c,a)\in \tiny{\CC^{d-1}}}$, the hypersurfaces
$Per_n(e^{i\theta})$ are contained in the bifurcation locus. This means that the functions $L_n$ are pluriharmonic on
$\CC^{d-1}\setminus Supp\;dd^c L$. The same is thus true for the limit $\varphi$.\\

\underline{Third step}: for any $\la_0\in \CC^{d-1}$ there exists a complex line ${\cal L}$ through $\la_0$ such that $\varphi=L$ on 
 the unbounded component of ${\cal L}\setminus\big({\cal L}\cap Supp (dd^{c} L)\big)$.\\
 
By Theorem \ref{controlinfty} we may pick a line ${\cal L}$ through $\la_0$ which meets infinity far from
the cluster set of $\cup_i{\cal B}_i$ in ${\bf P}_{\infty}$. This means that 
 for any $\la$ in the unbounded component of 
 ${\cal L}\setminus\big({\cal L}\cap Supp\;(dd^{c} L)\big)$ all critical points $c_i$ belong to the super-attracting basin of $\infty$ and thus,
$\la$ is a hyperbolic parameter.
This implies that $L_n(\la)=L_n^+(\la,e^{i\theta})$
and that, by Lemma \ref{lemLn+}, $\varphi(\la)=L(\la)$.\\

\underline{Fourth step}: $\varphi=L$ on $Supp\;dd^c L$.\\

This is the most delicate part of the proof, it somehow proceeds by induction on $d$. 
To simplify the exposition, we will only treat the cases  $d=2$ and $d=3$.\\

When $d=2$ the parameter space is  ${\bf C}$ and the bifurcation locus is the boundary $b{\cal M}$ of the Mandelbrot set.
The unbounded stable component $\big({\cal M}\big)^c$ is hyperbolic and thus, as we saw in the last step, $\left(\varphi-L\right)=0$ there.
Since $\left(\varphi-L\right)$ is negative and $u.s.c$, this implies that $\varphi=L$
on $b{\cal M}=Supp\;dd^c L$.\\
Let us stress that this ends the proof when $d=2$;
the complex line ${\cal L}$ of the third step is the parameter space itself in that case !\\

We now assume that $d=3$, the parameter space is then ${\bf C}^2$. Let us consider the sets 
$U_k$ of parameters which do admit an attracting $k$-cycle:

$$U_k:=\bigcup_{\vert w\vert <1} Per_k(w).$$

We have to show that $\left(\varphi-L\right)$ vanishes on the bifurcation locus.
Since the bifurcation locus is accumulated by curves of the form $Per_k(0)$ (by Theorem \ref{theodisatt} $\lim_k d^{-k}[\Per_k(0)]=\Bc$),
and the function $\left(\varphi-L\right)$ is negative and upper semi-continuous, it suffices to prove that
$$\left(\varphi-L\right)=0\;\textrm{on all sets}\;U_k.$$

 Let us first treat the problem on a curve ${\cal C}:=\Per_{k}(\eta)$ for $\vert \eta\vert <1$ and show that 
 $$(\star)\;\;\textrm{the sequence}\;L_n\vert_{\cal C}\;\textrm{ converges uniformly  to }\;L\vert_{\cal C} \;\textrm {on the stable components}.$$
 
We may assume that ${\cal C}$ is irreducible and 
desingularize it. This gives a one-dimensional holomorphic family
$(P_{\pi(u)})_{u\in M}$. Keeping in mind that the elements of this family are degree 3 polynomials which do admit an attracting basin of period $k$ and using the fact that the 
connectedness locus in ${\bf C}^2$ is compact (see Theorem \ref{controlinfty}), one sees that the family $(P_{\pi(u)})_{u\in M}$ enjoys the same properties than the quadratic polynomial family:
\begin{itemize}
\item[1-] the bifurcation locus is contained in the closure of  hyperbolic parameters
\item[2-] the set of non-hyperbolic parameters is compact in $M$.
\end{itemize}

Exactly as for the quadratic polynomial family this implies that the sequence $L_n\vert_{\cal C}$
 converges  in $L^1_{loc}$ to $L\vert_{\cal C}$ and the convergence is locally uniform on stable components since,
 as we already observed,  the functions $L_n$ are pluriharmonic there.\\

We now want to show that $\varphi= L$ on any open subset $U_k$.  Again, as the stable parameters are dense and 
$(\varphi-L)$  is $u.s.c$ and negative, it suffices to show that
$\varphi=L$ on any stable component  of $U_k$. On such a component the functions $L_n$ are pluriharmonic and thus actually converge
locally uniformly to $\varphi$. Then, ($\star$) clearly implies that $\varphi=L$ on $\Omega$.
 
\qed\\

\begin{rem}
The relatively simple behaviour of the bifurcation locus near infinity within polynomials families is crucial in the above proof. It is an open problem to show that
$\lim_n d^{-n}[\Per_n(e^{i\theta})]=\Bc$ in general families of rational maps. The first reasonable case to study would be  that of the moduli space $Mod_2$ for which
precise informations concerning the bifurcation locus at infinity have been obtained  by Epstein (\cite{epstein}).
\end{rem}

\section{Laminated structures in bifurcation loci}\label{secLam}
\subsection{Holomorphic motion of the Mandelbrot set in $Mod_2$}

We work here in the moduli space $Mod_2$ of degree two rational maps which, as we saw in section \ref{secMod}, can be identified
to $\CC^2$. 
Our aim is to show that the bifurcation locus in the region
$$U_1:=
\{\la\in{\bf C}^2/ f_{\la} \;\textrm{has an attracting fixed point}\}$$
can be obtained by holomorphically moving the boundary $b{\cal M}$ of the Mandelbrot set. We remind that $b{\cal M}$ is the bifurcation locus of  $Per_1(0)$
which is a complex line contained in $U_1$ and can be identified to the family of quadratic polynomials.\\
We will see simultaneously that the bifurcation current is uniformly laminar in the region $U_1$. Let us first recall some basic facts about holomorphic motions.

\begin{defi} 
Let $M$ be a complex manifold and $E\subset M$ be any subset. A \emph{holomorphic motion} of $E$ in $M$ is 
a map 
\begin{center}
$\sigma : E\times \Delta \ni(z,u)\mapsto \sigma(z,u)=:\sigma_{u}(z)\in M$
\end{center}
which satisfies the following properties:
\begin{itemize}
\item[i)]$\sigma_0=Id\vert_E$
\item[ii)]$E\ni z\mapsto \sigma_{u}(z)\in M$ is one-to-one for every $u\in \Delta$
\item[iii)]$\Delta \ni u \mapsto \sigma_{u}(z)\in M$ is holomorphic for every $z\in E$. 
\end{itemize} 
\end{defi}

When the family of holomorphic discs in $M$ enjoys good compactness properties, any holomorphic motion extends to the
closure. In particular, when $M$ is the Riemann sphere  ${\widehat \CC}$, the Picard-Montel theorem combined with Hurwitz lemma
easily leads to some famous extension statement which is usually called $\la$-lemma since the "time" parameter is 
denoted $\la$ rather than $u$ (see Lemma \ref{lamlem}).\\

The main result of this subsection is the following. It was first proved by Goldberg and Keene (see \cite {GK}). The proof we present here exploits the formalism of bifurcation currents and is due to Bassanelli and the author (see \cite{BB2}), this approach turns out to be much simpler and also provides some information on the laminarity of the bifurcation currents.

\begin{theo}\label{theomouvMand} Let $\Omega_{hyp}$ be the union of all hyperbolic
components of the Mandelbrot set ${\cal M}$ and $\car$ the main cardioid. 
Let $\Bif_1$ be the bifurcation locus in $U_1$ and $\Bc{\arrowvert_{U_1}}$ be 
the associated bifurcation current. Let $ \mu_{1}$ be the harmonic measure of ${\cal M}$.\\
There exists a continuous holomorphic motion 
$$\sigma : \big(\big(\Omega_{hyp}\setminus \car\big)\cup b{\cal M}\big)\times \Delta \to U_1$$
such that
$$\sigma\left(b{\cal M}\times\Delta\right)=\Bif_1\;\;\textrm{and}\;\;
\Bc{\arrowvert_{U_1}} = \int_{Per_1(0)}[\sigma (z,\Delta)]\;\mu_{1}.$$

In particular, $ \Bif_1$  is a lamination with $ \mu_{1}$ as transverse measure. 
Moreover, the map $\sigma$ is holomorphic on $\big(\Omega_{hyp}\setminus \car\big)\times \Delta$ and preserves the curves 
$Per_n(w)$ for $n\ge 2$ and $\vert w\vert\le 1$.
\end{theo}

\proof

\underline{First step:} Holomorphic motion of $\big(\Omega_{hyp}\setminus \car\big)$.\\

The curve $Per_1(0)$ is actually the complex line $\la_1=2$. We will write $(2,\la_2)=:z$ the points of this line.\\

Let us consider $U_n:=
\{\la\in{\bf C}^2/ f_{\la} \;\textrm{has an attracting cycle of period}\; n \}$ and $\Omega_n:=U_n\cap Per_1(0)$.
We recall that 
$\Omega_{hyp}:= \car\cup\bigcup_{n\ge 2}\Omega_n$.\\

Let us also set $U_{n,1}:=U_n\cap U_1$. By the Fatou-Shishikura inequality, a quadratic rational map has at most two non-repelling cycles. Thus $U_{n,1}\cap U_{m,1}=\emptyset$ when $n\ne m$ and there exists a well defined holomorphic map
$$\psi_n:U_{n,1}\to\Delta\times\Delta$$
which associates to every $\la\in U_{n,1}$ the pair $(w_n(\la),w_1(\la))$ where $w_n(\la)$ is the multiplier of the attracting $n$-cycle of $\la$
and $w_1(\la)$ the multiplier of its attracting fixed point.
 
The cornerstone is the following transversality statement due to Douady and Hubbard (see also \cite{BB1}):
\begin{lem}\label{lemTrDH}
The map $\psi$ induces a biholomorphism
$$\psi_{n,j} : {U_{n,1,j}}\to \Delta\times\Delta$$ on each connected component ${U_{n,1,j}}$ of
$U_{n,1}$.
\end{lem}

By the above Lemma, the connected components $\Omega_{n,j}$ of $\Omega_n$ coincides with $U_{n,1,j}\cap Per_1(0)$ and one clearly obtains a holomorphic motion 
$\sigma: \big(\Omega_{hyp}\setminus \car\big)\times \Delta\to U_1$ by setting:
$$\sigma(z,t):=({\psi_{n,j}})^{-1}\big( w_n (z),t\big)$$
for any $z\in \Omega_{n,j}$.\\

\underline{Second step:} extension of $\sigma$ to $b{\cal M}$.\\

The key point here is that $\sigma(z,t)=:(\alpha(z,t),\beta(z,t))$ belongs to the complex line
$Per_1(t)$ which, according to Proposition \ref{Per/L}, is given by the equation 
$$(t^2+1)\la_1 - t\la_2 -(t^3+2)=0.$$

Thus $\sigma(z,t)$ is completely determined by $\beta(z,t)$:
\begin{eqnarray}\label{albet}
\alpha (z,t)= \frac{1}{1+t^2}\left(t\beta (z,t)+t^3+2\right),\;\;\forall t\in \Delta.
\end{eqnarray}

We will now identify $Per_1(0)$ with 
 the deleted Riemann sphere ${\widehat \CC}\setminus\{\infty\}$ and set $\beta(\infty,t)=\infty$ for all $t\in \Delta$. Then, the map 
$\beta: \big(\{\infty\}\cup\big(\Omega_{hyp}\setminus \car\big)\big)\times \Delta\to {\widehat \CC}$ is clearly a holomorphic motion
which, by Lemma \ref{lamlem}, extends to the closure of $\big(\Omega_{hyp}\setminus \car\big)$.
We thus obtain a continuous holomorphic motion 
$$\beta: \big(\{\infty\}\cup\big(\Omega_{hyp}\setminus \car\big)\cup b{\cal M}\big)\times \Delta\to {\widehat \CC}.$$
As, by construction, $\beta(z,t)\ne \infty$ when $z\ne \infty$, the identity \ref{albet} shows that $\sigma(z,t)=(\alpha(z,t),\beta(z,t))$
extends to a continuous holomorphic motion of $\big(\big(\Omega_{hyp}\setminus \car\big)\cup b{\cal M}\big)$.\\

\underline{Third step:} laminarity properties.\\

Let us show that $\Bc{\arrowvert_{U_1}} = \int_{Per_1(0)}[\sigma (z,\Delta)]\;\mu_{1}.$
According to the approximation formula given by Theorem \ref{theodisatt} and applied on $Per_1(0)$ we have:
\begin{eqnarray}\label{z}
\mu_{1}=\lim_m 2^{-m}\sum_{z\in Per_1(0)\cap Per_m(0)}\delta_{\sigma(z,0)}.
\end{eqnarray}
Let us set $T:=\int_{Per_1(0)}[\sigma(z,\Delta)]\;\mu_{1}$. We have to check that $T=\Bc{\arrowvert_{U_1}}$.
Let $\phi$ be a $(1,1)$-test form in $U_1$. As the holomorphic motion $\sigma$ is continuous, the function
$z\mapsto \langle[\sigma(z,\Delta)],\phi\rangle$ is continuous as well. Then, using \ref{z} one gets
\begin{eqnarray}\label{zz}
\langle T,\phi\rangle =\lim_m 2^{-m}\sum_{z\in Per_1(0)\cap Per_m(0)}\langle[\sigma(z,\Delta)],\phi\rangle =
\lim_m2^{-m}\langle [Per_m(0)],\phi\rangle
\end{eqnarray}
where the last equality uses the fact that, according to Proposition \ref{multi}, 
 the curves $Per_m(0)$ have no multiplicity in $U_1$. 
Now the conclusion follows by using \ref{zz} and the approximation formula of Theorem \ref{theodisatt} in $U_1$.\\

By construction,
the map $\sigma$ is holomorphic on $\big(\Omega_{hyp}\setminus \car\big)\times \Delta$ and preserves the curves 
$Per_n(w)$ for $n\ge 2$ and $\vert w\vert<1$. This extends to $\vert w\vert =1$ by continuity.\\ 

Using the continuity of $\sigma$, one easily sees that $\sigma\big(b{\cal M},\Delta\big)$ is closed in $U_1$ and therefore contains the support of $\int_{Per_1(0)}[\sigma(z,\Delta)]\;\mu_{1}$. By the above formula we thus have 
$$\Bif_1=Supp\big(\Bc{\arrowvert_{U_1}}\big)\subset
\sigma\big(b{\cal M},\Delta\big).$$ The opposite inclusion easily follows from the construction of $\sigma$:
any point in $\sigma\big(b{\cal M},\Delta\big)$ is a limit of $z_m\in Per_1(0)\cap Per_m(0)$ where $m\to+\infty$.\qed\\ 

%COROLLARY: Closure of Bif in P2???.\\

Instead of using the basic $\la$-Lemma we could have use its far advanced generalization due to  Slodkowski 
and get a motion on the full line $Per_1(0)$. 

\begin{theo}{\bf (Slodkowski $\la$-lemma)}\label{lamslo}
Let $E\subset {\widehat \CC}$ be a subset of the Riemann sphere and
$\sigma : E\times \Delta \ni(z,t)\mapsto \sigma(z,t)\in {\widehat \CC}$ be a holomorphic motion.
Then $\sigma$ extends to a holomorphic motion $\tilde \sigma$ of ${\widehat \CC}$. Moreover
$\tilde\sigma$ is continuous on $\overline{E}\times \Delta$ and
$z\mapsto {\tilde\sigma}(z,t)$ is a $K$-quasi-conformal
homeomorphism for $K:=\frac{1+\vert t\vert}{1-\vert t\vert}$.
\end{theo}

Our reference for quasi-conformal maps is the book \cite{Hub} where one can also find 
a nice proof of Slodkowski theorem due to Chirka and Rosay.\\

Using Slodkowski Theorem one may obtain further informations on the motion given by Theorem \ref{theomouvMand}.
We refer to our paper \cite{BB2} for a proof.

\begin{theo}\label{farf}
Let 
$\sigma : \big(\big(\Omega_{hyp}\setminus \car\big)\cup b{\cal M}\big)\times \Delta \to U_1$
be the holomorphic motion given by Theorem \ref{theomouvMand}. Then $\sigma$ extends to
 a continuous holomorphic motion
${\tilde\sigma}:Per_1(0)\times\Delta \longrightarrow U_1$ which is onto.  All stable components 
in $U_1$ are of the form 
${\tilde\sigma}\left(\omega\times\Delta\right)$ for some stable component $\omega$ in 
$Per_1(0)$.
Moreover, the map $z\mapsto{\tilde\sigma}(z,t)$ is a
quasi-conformal homeomorphism for each $t$ and
$\tilde\sigma$ is one-to-one on $\left( Per_1(0)\setminus{\overline\car}\right)\times \Delta$ where $\car$ is the 
 main cardioid.  
\end{theo}

Theorem \ref{farf} shows that non-hyperbolic components exist in $U_1$ if and only if such components exist within the quadratic 
polynomial family $Per_1(0)$.
Let us underline that,  in relation with 
Fatou's problem on the density of hyperbolic rational maps, it is conjectured that such components do not exist.\\

It might be useful to note that holomorphic motions enjoy good H\"{o}lder regularity properties. We will end this subsection by 
giving a basic result in this direction.

\begin{lem}\label{HoHoMo}
Let $h: B(0,r)\times E\longrightarrow\pp$ be a holomorphic motion of some $E\subset\pp$ parametrized by a ball in $\CC^k$. Then, for any $z_0\in E$ we may find  $0<r_1\leq r$ and $\eta>0$ such that the following estimates hold for $0<r'<r_1$,  $\lambda\in B(0,r')$ and $z,z'\in D(z_0,\eta)\cap E$:
\begin{eqnarray*}
C'(r')|z-z'|^{\frac{r_1+\|\lambda\|}{r_1-\|\lambda\|}}\leq|h(\lambda,z)-h(\lambda,z')|\leq C(r')|z-z'|^{\frac{r_1-\|\lambda\|}{r_1+\|\lambda\|}}
\label{biholder}
\end{eqnarray*}
where  $C(r'),C'(r')$ are strictly positive constants.
\end{lem}

\proof   Let $z_0\in E$. As $h$ is continuous, one finds $0<\eta<1$ and $0<r_1\leq r$ such that
\begin{center}
$|h(\lambda,z)-h(\lambda,z_0)|<1$
\end{center}
for any $\lambda\in B(0,r_1)$ and any $z\in D(z_0,\eta)\cap E$. Let us now pick two distinct points  $z,z'$ in $ D(z_0,\eta)\cap E$, a parameter $\lambda$ in $B(0,r_1)\setminus\{0\}$ and set
\begin{center}
$g_{z,z'}^\lambda(t):= -\log\displaystyle\frac{|h(r_1t\frac{\lambda}{\|\lambda\|},z)-h(r_1t\frac{\lambda}{\|\lambda\|},z')|}{2},\ $ $t\in D$.
\end{center}
Since $z\neq z'$ and $z,z'\in D(z_0,\eta)\cap E$ one has
\begin{center}
$0<|h(r_1t\frac{\lambda}{\|\lambda\|},z)-h(r_1t\frac{\lambda}{\|\lambda\|},z')|<2$
\end{center}
for any $t\in \Delta$ and therefore $g_{z,z'}^\lambda$ is a positive harmonic function on the unit disc $ \Delta$. Harnack inequalities yield
\begin{eqnarray*}
\frac{1-|t|}{1+|t|}g_{z,z'}^\lambda(0)\leq g_{z,z'}^\lambda(t)\leq\frac{1+|t|}{1-|t|}g_{z,z'}^\lambda(0)
\end{eqnarray*}
for any $t\in D$, which means that
\begin{eqnarray*}
\left(\frac{|z-z'|}{2}\right)^{\frac{1+|t|}{1-|t|}}\leq \frac{|h(rt\frac{\lambda}{\|\lambda\|},z)-h(rt\frac{\lambda}{\|\lambda\|},z')|}{2}\leq\left(\frac{|z-z'|}{2}\right)^{\frac{1-|t|}{1+|t|}}
\end{eqnarray*}
for $t\in D$. Taking $t=\|\lambda\|/r_1$ we get
\begin{eqnarray*}
\left(\frac{|z-z'|}{2}\right)^{\frac{r_1+\|\lambda\|}{r_1-\|\lambda\|}}\leq \frac{|h(\lambda,z)-h(\lambda,z')|}{2}\leq{\left(\frac{|z-z'|}{2}\right)}^{\frac{r_1-\|\lambda\|}{r_1+\|\lambda\|}}.
\end{eqnarray*}
The conclusion follows by setting $C(r'):=2^{\frac{2r'}{r_1+r'}}$ and $C'(r'):=2^{\frac{-2r'}{r_1-r'}}$ for $0<r'<r_1$.\qed

\subsection{Further laminarity statements for $\Bc$}

The following result is an analogue of Theorem \ref{theomouvMand} in the regions
$$U_n:=
\{\la\in Mod_2 / f_{\la} \;\textrm{has an attracting cycle of period}\; n \}.$$

It shows, in particular, that the bifurcation current in $Mod_2$ is uniformly laminar in the regions $U_n$.
It has been established by Bassanelli and the author in \cite{BB2}.

\begin{theo}\label{thMouvBif}
Let $\Bif_n$ be the bifurcation locus in $U_n$ and $\Bc{\arrowvert_{U_n}}$ be 
the associated bifurcation current.
Let $\Bif_{n}^{\;c}$ be the bifurcation locus in the central curve $Per_n(0)$
and $\mu_{n}^c$ be the associated bifurcation measure.
Then, there exists a map
\begin{eqnarray*}
\sigma:\Bif_{n}^{\;c}\times\Delta & \longrightarrow &\Bif_n \\
(\la,t)&\longmapsto&\sigma(\la,t)\\
\end{eqnarray*}
such that:

\begin{itemize}
\item[1)] $\sigma\left(\Bif_{n}^{\;c}\times\Delta\right)=\Bif_n$
\item[2)] $\sigma$\; \textrm{\it is continuous,} $\sigma(\la,\cdot)\; \textrm{\it is one-to-one and holomorphic for each }\la \in \Bif_{n}^{\;c}$  
\item[3)] $p_n\left(\sigma(\lambda,t),t\right)=0;\;\forall \la\in \Bif_{n}^{\;c},\forall t\in \Delta$
\item[4)] \textrm{\it the discs} 
$\left(\sigma(\la, \Delta)\right)_{\la\in\Bif_{n}^{\;c}}$ 
\textrm{\it are mutually disjoint}.
\end{itemize}

Moreover the bifurcation current in $U_n$ is given by
$$\Bc{\arrowvert_{U_n}} = \int_{\Bif_{n}^{\;c}}[\sigma\left(\la,\Delta\right)]\;\mu_{n}^c$$

and, in particular, $ \Bif_n$  is a lamination with $ \mu_{n}^c$ as transverse measure. 
\end{theo}

The proof is similar to that of Theorem \ref{thMouvBif} but requires a special treatment for the extension problems since there
is no $\la$-lemma available. The key is to use the fact that, by construction, the starting motion
$ \sigma:\left(\cup_{m\ne n} (U_m\cap U_n)\right)\times\Delta \longrightarrow U_n$ satisfies the following
property:
$$p_n\left(\sigma(\lambda,t),t\right)=0;\;\;\forall \la\in \Bif_{n}^{\;c},\;\;\forall t\in \Delta.$$ 

Such a motion is what we call a 
$p_n$-guided holomorphic motion. Using Zalcman rescaling lemma, one proves the following compactness property for
guided holomorphic motions. We stress that here, an holomorphic motion ${\cal G}$ is seen as a family of disjoints holomorphic discs
$\sigma$ and
${\cal G}_{t_0}$ is the set of points $\sigma(t_0)$.

\begin{theo}\label{thext}
Let $p(\la,w)$ be a polynomial on ${\bf C}^2\times{\bf C}$ such that the degree of $p(\cdot,w)$ does not depend on $w\in \Delta$.
Let ${\cal G}$ be a $p$-guided holomorphic motion in ${\bf C}^2$ such that any component of the algebraic curve $\{p(\cdot,t)=0\}$
contains at least three points of ${\cal G}_t$ for every $t\in \Delta$.
Then, for any ${\cal F}\subset{\cal G}$ such that ${\cal F}_{t_0}$ is relatively compact in ${\bf C}^2$ for some $t_0\in \Delta$,
there exists a continuous $p$-guided holomorphic motion $\widehat{{\cal F}}$ in ${\bf C}^2$
such that ${\cal F}\subset \widehat{{\cal F}}$
and $\widehat{{\cal F}}_{t_0}=\overline{{\cal F}_{t_0}}.$
\end{theo}

The above Theorem plays the role of the $\la$-lemma in the proof of Theorem \ref{thMouvBif}. We refer to the paper \cite{BB2} for details.\\

We will now end this section by presenting some more precise laminarity result for $\Bc$ which is due to Dujardin (see \cite{dujardin2})

\begin{theo}\label{theoDuL}
Within the polynomial family of degree $3$,
the bifurcation current is laminar on every open set where one critical point is passive and, in particular, outside the connectedness locus.
\end{theo}

We refer to section \ref{conloc} for the discussion of the connectedness locus within polynomial families.\\

A few comments on the concept of laminar currents are necessary here. We restrict ourself to the case of a positive $(1,1)$-current $T$ on a complex manifold $M$. One says that $T$ is
\emph{locally uniformly laminar} on $M$ if any point of the support of $T$ admits a neighbourhood on which $T$ is of the form
\begin{center}
$T=\int_\tau [\Delta_t]\; \mu(t)$
\end{center}
where $\Delta_t$ is a lamination by holomorphic discs with a transverse measure $\mu$.\\

With this terminology, Theorem \ref{thMouvBif} says that $\Bc$ is uniformly laminar on the open sets of $Mod_2$ where an attracting basin of given period exists. In particular,
$\Bc$ is locally uniformly laminar on any open set of $Mod_2$ where a critical point is attracted by a cycle.\\

One says that $T$ is \emph{laminar} on $M$ if there exists a sequence of open sets $\Omega_i \subset M$ and a sequence of  currents $T_i$ which are locally uniformly laminar
on $\Omega_i$ and such that the sequence $\left(T_i\right)_i$ increasingly converges to $T$. This is actually equivalent to say that 
\begin{center}
$T=\int_A [\Delta_a]\;\mu(a)$
\end{center}
 where $\left( \Delta_a\right)_{a\in A}$ is a family of compatible (in the sense of analytic continuation) holomorphic discs parametrized by an abstract set $A$ and $\mu$ is a measure on $A$.\\

The methods used in $Mod_2$ for proving Theorem \ref{thMouvBif}  work also in the cubic polynomial family (the situation is actually technically simpler) and thus, for this family too, the
bifurcation current is locally uniformly laminar on the open subsets where a critical point is attracted by a cycle.\\
As it is not hard to see, this would imply Theorem \ref{theoDuL} if the hyperbolicity conjecture was known to hold in the cubic polynomial family. To overcome this difficulty, one uses more sophisticated tools on laminar currents. Precisely, the proof of theorem \ref{theoDuL}  is based on a laminarity criterion due to De Thelin (\cite{DeT}) which says that a 
current $T:=\lim_n \frac{1}{d_n} [C_n]$, where $[C_n]$ is a sequence of integration currents on curves in $\CC^2$, is laminar if $\textrm{\it genus}\;(C_n)=O(d_n)$.\\

\chapter{The bifurcation measure}
 
The powers $(\Bc)^k$ of the bifurcation current detect stronger bifurcations and allow to define a very interesting stratification of the bifurcation locus.  Among them, the highest power is a measure whose support should be the seat of the strongest bifurcations.  In this chapter, we will survey the basic properties of this measure
 and then describe its support as the closure of some sets of remarkable parameters. Although some results about  the intermediate powers of $\Bc$ can be easily deduced from our exposition, we will not discuss explicitely these currents here and refer the reader to the survey of Dujardin \cite{DuS} for more details on this topic. 

\section{A Monge-Amp\`ere mass related with strong bifurcations}
\subsection{Basic properties}

Since the bifurcation current $\Bc$  of any holomorphic family $\left(f_{\la}\right)_{\la\in M}$ of degree $d$ rational maps has a continuous potential $L$ (see Definition \ref{defiBC} and Theorem \ref{theoLH}), one may define the powers $(\Bc)^k:=\Bc\wedge\Bc \wedge \cdot\cdot\cdot \wedge\Bc$ for any $k\le m:=dim M$. 
We recall that for any closed positive current $T$, the product $dd^c L\wedge T$ is defined by $dd^c L\wedge T:=dd^c (LT)$. In particular,
$(\Bc)^m$ is a positive measure on $M$ which is equal to the Monge-Amp\`ere mass of the Lyapunov function $L$.

\begin{defi}
Let $\left(f_{\la}\right)_{\la\in M}$ be a
holomorphic family  of degree $d$ rational maps parametrized by a complex manifold $M$ of dimension $m$. The \emph{bifurcation measure} $\Bm$
of the family is the positive measure on $M$ defined by
\begin{center}
$\Bm=\frac{1}{m!} (\Bc)^m=\frac{1}{m!}\left(dd^c L\right)^m$
\end{center} 
where $\Bc$ is the bifurcation current and $L$ the Lyapunov function of the family.
\end{defi}

The following proposition is a direct consequence of the definition and the fact that $\Bm$ has locally bounded potentials.

\begin{prop}\label{propCLN}
The support of $\Bm$ is contained in the bifurcation locus and $\Bm$ does not charge pluripolar sets.  
\end{prop}

It is actually possible to define the bifurcation measure in the moduli space $Mod_d$ of degree $d$ rational maps and show that this measure has
strictly positive and finite mass (see \cite{BB1} Proposition 6.6). Although all the results we will present here are true in $Mod_d$, we will restrict ourself to the
technically simpler situation of holomorphic families. The example we have in mind are the polynomial families
and the moduli space $Mod_2$ which, in some sense, can be treated as a holomorphic family (see Theorem \ref{theoMiMod2}).\\

In arbitrary holomorphic families, the measure $\Bm$ can identically vanish. Moreover, when $\Bm >0$, it is usually quite involved to prove it.
Note however that this will follow from standard arguments in polynomial families. The following simple observation
already shows that $\Bm>0$ in $Mod_2$ (and more generally $Mod_d$), it has also its own interest.

\begin{prop}\label{propLatMb}
In any holomorphic family, all rigid Latt\`es examples belong to the support of $\Bm$.
\end{prop}

\proof
The parameters corresponding to rigid Latt\`es examples are isolated.
Thus, as Theorem \ref{theoCarLat} shows,  the Lyapunov function $L$ takes the value $\frac{\ln d}{2}$  and has a strict minimum at any such parameter.
Applying the comparison principle to $L$ and some constant function $\frac{\ln d}{2}+\epsilon $, one sees that the Monge-Amp\`ere measure $\left(dd^c L\right)^m$ cannot vanish around a strict minimum of $L$. This yields the conclusion.
\qed\\

Using pluripotential theory, it is possible to show that any Latt\`es example in $Mod_d$ lies in the support of $(\Bc)^{2d-3}$. Buff and Gauthier
have recently shown that such maps actually
belong to the support of $\Bm$ (see \cite{BG}). Their proof requires transversality statements and uses the quadratic differentials techniques.

\begin{theo}
In $Mod_d$,
 Latt\`es examples belong to the support of the bifurcation measure.
\end{theo}

We end this subsection by showing that the activity currents have no self-intersection. This is a useful geometric information which,
in particular,  shows that the activity of all critical points is a necessary condition for a parameter
to be in the support of the bifurcation measure. It was first proved by 
Dujardin-Favre in the context of polynomial families (see \cite{DF} Proposition 6.9), 
we present here a general argument due to Gauthier (\cite{Gau}).

\begin{theo}\label{theoNIAC}
Let $\left(f_{\la}\right)_{\la\in M}$ be any holomorphic family of degree $d$ rational maps
with marked critical points. The activity currents $T_i$ satisfy
$T_i\wedge T_i=0$. In particular, when  $m:=dim\;M= 2d-2$ (or $m= d-1$
for polynomial families)
then
\begin{center}
$\Bm= T_1\wedge T_2\wedge\cdot\cdot\cdot\wedge T_{m}$
\end{center}
and the support of $\Bm$ is contained in the intersection of the activity loci of the critical points.

\end{theo}

The proof is very close to that of a density statement which will be presented in the next section,
it combines the following potential-theoretic Lemma with a dynamical observation.

\begin{lem}\label{lemBrDu}
Let $u$ be a continuous $p.s.h$ function on some open subset $\Omega$ in $\CC^2$.
Let $\Gamma$ be the union of all analytic subsets of $\Omega$ on which $u$ is harmonic.
If the support of $dd^c u$ is contained in $\overline{\Gamma}$ then $dd^c u\wedge dd^c u$
vanishes on $\Omega$.
\end{lem}

\proof
Let us set $\mu:=dd^c u\wedge dd^c u$. 
Let $B_r$ be an open ball of radius $r$ whose closure is contained in $\Omega$, we have to show that $\mu\left(B_{\frac{r}{2}}\right)=0$.\\

Denote by $h$ the solution of the Dirichlet-Monge-Amp\`ere problem with data $u$ on $bB_r$:

\begin{center}
$h=u\;\textrm{on the boundary of}\; B_r$\\
$dd^c h\wedge dd^c h=0\;\textrm{on}\;B_r\;(\textrm{i.e.}\;h\;\textrm{is maximal on}\;B_r)$.
\end{center}

The function $h$ is $p.s.h$ and continuous on $\overline{B_r}$ (see \cite{BT}).
As $h$ is $p.s.h$ maximal and coincides with the $p.s.h$ function $u$ on $b\overline{B_r}$, we have $u\le h$ on $\overline{B_r}$. For any $\epsilon>0$ we define 
\begin{center}
$D_\epsilon:=\{\la\in B_{\frac{r}{2}}\;/\; 0\le h(\la)-u(\la) \le \epsilon\}$.
\end{center}
We will see that our assumption implies that
\begin{eqnarray}\label{D_e}
Supp\;\mu\cap B_{\frac{r}{2}} \subset D_\epsilon\;\textrm{for all }\;\epsilon>0.
\end{eqnarray}

Indeed, if $\gamma$ is a complex curve in $\Omega$ on which $u$ is harmonic then, the maximum modulus principle, applied to
$(h-u)$ on $\gamma\cap B_r$ implies that $h=u$ on $\gamma$.
Then, as $(h-u)$ is continuous on $B_r$ and 
$Supp\; \mu\subset Supp\;dd^c u \subset \overline{\Gamma}$, we get $Supp\;\mu\cap B_r\subset \{h=u\}$.\\

Now, a result due to Briend-Duval (see \cite{BrDu} or \cite{DS2} Th\'eor\`eme A.10.2) says that 
\begin{eqnarray}\label{BDu}
\mu\left(D_\epsilon\right) \le C\epsilon
\end{eqnarray}
where $C$ only depends on $u$ and $B_r$. From \ref{D_e} and \ref{BDu} we deduce that $\mu\left(B_{\frac{r}{2}}\right)=0$.\qed\\

We may now end the proof of Theorem \ref{theoNIAC}.\\

\proof
We only  have to show that $T_i\wedge T_i=0$, the remaining then follows from the identity
$\Bc=\sum_i T_i$ (see Theorem \ref{theosuppBC}).\\
The statement is local and we may therefore assume that $M=\CC^k$. Moreover, an elementary slicing argument 
allows to reduce the dimension to $k=2$.
We apply the above Lemma with $u=G_\la\left(\hat c_i(\la)\right)$ (see Lemma \ref{lemPotAC}).
We have to show that the support of $dd^c u=T_i$ is accumulated by curves on which the critical point $c_i$
is passive. These curves are of the form $\{f_\la^n(c_i(\la)=c_i(\la)\}$ and their existence follows from  Lemma \ref{lemappperno}.
\qed.\\

\begin{rem}\label{rkdou}
An example, due to A.Douady, shows that the activity of all critical points is not sufficent for a parameter to be in the support of the bifurcation measure.
We will present this example in the next subsection (see Example \ref{exaDou}).
\end{rem}

\subsection{Some concrete families}

We first discuss the case of the polynomial families introduced in subsection \ref{ssPd}.
We follow here the paper \cite{DF} by Dujardin and Favre.

\begin{prop}
The bifurcation measure $\Bm$ of the degree $d$ polynomial family  
$\big(P_{c,a}\big)_{(c,a)\in {\tiny {\bf C}^{d-1}}}$
is a probability measure supported on the connectedness locus ${\cal C}$.
It coincides with the pluricomplex equilibrium measure of the compact set ${\cal C}$
and its support is the Shilov boundary of ${\cal C}$.
\end{prop}

\proof
Let us recall that the Green function of the polynomial $P_{c,a}$ is denoted $g_{c,a}$.
The connectedness locus ${\cal C}$ is a compact subset of $ {\bf C}^{d-1}$ which coincides
with the intersection 
 $\cap_{0\le i\le d-2} {\cal B}_i$  where ${\cal B}_i$ is the set of parameters for which the orbit of the critical point $c_i$ is bounded (see Theorem \ref{controlinfty}).
 As the support of the activity current $T_i$ is contained in $b{\cal B}_i$ we deduce that $Supp\;\Bm\subset {\cal C}$ from Theorem \ref{theoNIAC}.
 All the remaining follows from the fact that 
 \begin{center}
 $\Bm=\left( dd^c  {\cal G}\right)^{d-1}$
 \end{center}
  where 
 \begin{center}
 ${\cal G}:=sup\{u\;p.s.h\;/\; u-\ln^+ max\{\vert a\vert,\vert c_k\vert\}\le O(1), u\le 0\;\textrm{on}\;{\cal C}\}$.
 \end{center}
 is the pluricomplex Green function of $\cal C$ with pole at infinity.
An identity which we shall now prove.\\

 Let us first establish that $\Bm=\left( dd^c  G\right)^{d-1}$ where $G:=max\{g_0,g_1,\cdot\cdot\cdot,g_{d-2}\}$ and $g_i:=g_{c,a}(c_i)$.
 We show by induction that $T_0\wedge T_1\wedge\cdot\cdot\cdot\wedge T_l=\left( dd^c  G_l\right)^{l+1}$ for $0\le l\le d-2$
 where $G_l:=max\{g_0,g_1,\cdot\cdot\cdot,g_{l}\}$. This 
comes from the following computation:
 \begin{eqnarray*}
 T_0\wedge T_1\wedge\cdot\cdot\cdot\wedge T_{l-1}\wedge T_l=dd^c \left( g_l (dd^c G_{l-1})^{l}\right)=dd^c \left(G_l (dd^c G_{l-1})^{l}\right)=\\
 dd^c \left(G_{l-1}(dd^c G_{l-1})^{l-1}\wedge dd^c G_l\right)= dd^c \left(G_{l}(dd^c G_{l-1})^{l-1}\wedge dd^c G_l\right)=\left(dd^c G_l\right)^{l+1}
\end{eqnarray*}

 the second equality follows from $g_l=G_l$ on the support of $ (dd^c G_{l-1})^{l}$ and the fourth from $G_l=G_{l-1}$ on the support of $dd^c G_{l}$,
 the last equality is obtained by repeating the same arguments $l-1$ times.\\
 
 It remains to show that $G={\cal G}$. The proof is standard and relies on the estimate given by Proposition \ref{estimgreen} 
 and the fact that $G$ is maximal outside ${\cal C}$. We refer to the paper \cite{DF}, Proposition 6.14 for more details.\qed\\
 
 \begin{rem}\label{remShi}
 The above result shows that the support of the bifurcation measure is topologically much smaller than the bifurcation locus. Indeed, the Shilov boundary is usually a tiny part of the
 full boundary. For instance, the boundary of the bidisc $\Delta\times\Delta$ is $\left(\Delta\times S^1\right)\cup \left(S^1\times \Delta\right)$ while its Shilov boundary is the real torus $S^1\times S^1$.
 \end{rem}
 
 We will now present the example mentionned at the end of last subsection.
 
 \begin{exa}\label{exaDou}
 In the holomorphic family of degree $3$ polynomials
 \begin{center}
 $\left((1+\alpha_1) z+ (\frac{1}{2}+\alpha_2) z^2 + z^3\right)_{\alpha\in V_0}$
 \end{center}
  where $V_0$ is a neighbourhood of the origin in $\CC^2$
 the critical points are both active at the origin $(0,0)$ but $(0,0)\notin Supp\;\Bm$.
 \end{exa}
 
This family is a deformation of the polynomial
 $P_0:=z+\frac{1}{2} z^2 + z^3$. If $V_0$ is small enough we have two marked critical points $c_1(\alpha)$ and $c_2(\alpha)$.
The origin $0$ is a parabolic fixed point for $P_0$ and we may assume that $P_{\alpha}$  has two fixed points counted with multiplicity near $0$
for all $\alpha \in V_0$. As $P$ is real, its critical points are complex conjugate and both of them are attracted by the 
parabolic fixed point at $0$, moreover their orbits are not stationnary.\\

 We first show that $(0,0)\notin Supp\;\Bm$.
 When the fixed points of $P_{\alpha}$ are distinct,
 we denote by $m_1(\alpha)$ and $m_2(\alpha)$ their multipliers.
When this is the case, it turns out that either $\vert m_1(\alpha)\vert <1$ or $\vert m_2(\alpha)\vert <1$ and thus one of the fixed points attracts a critical point.
This can be seen by using the holomorphic fixed point formula (see \cite{Milnor4}, Lecture 12).\\
By Theorem \ref{theoNIAC}, this implies that $\alpha\notin Supp\; \Bm$. We thus see that $\Bm$ is supported on the subvariety of parameters $\alpha$
for which the fixed point is double. By Proposition \ref{propCLN} this implies that $\Bm$ vanishes near $(0,0)$.\\

Let us now see that both critical points are active.  We may assume that the family is parametrized by a disc $D$ in $\CC$
such that $P_\alpha$ has two distinct fixed points when $\alpha\ne 0$. Assume to the contrary that a critical point $c(\alpha)$ is passive.  Then, after taking a subsequence, the sequence  $u_n(\alpha):=P_{\alpha}^n (c(\alpha))$ is uniformly converging  to $u(\alpha)$. Since the polynomial $P_0$ is real, its critical points are complex conjugate and must therefore
both be attracted by the parabolic fixed point $0$. Thus $u_n(0)$ converges to $0$ and never belongs to the analytic set $Z:=\{(\alpha,z)\in D\times \pp\;/\; P_\alpha(z)=z\}$.
If $u_n(\alpha)\in Z$ for $\alpha\ne 0$ then the orbit of $c(\alpha)$ is captured by a fixed point which, since $c(\alpha)$ is passive, must be attracting. In that case 
the curve $(\alpha,u(\alpha))$ is contained in $Z$. If $u_n(\alpha)$ never belongs to $Z$ then, by Hurwitz lemma, the curve $(\alpha,u(\alpha))$ is also contained in $Z$.
This is impossible since, for some $\alpha$ close to $0$, the critical orbit should be attracted by a repelling fixed point.\qed\\

The situation in the moduli space $Mod_2$ is more complicated. We recall that $Mod_2$ can be identified to $\CC^2$. Using the results which will be obtained in the last section of this chapter and the holomorphic motions
constructed in section \ref{secLam}, it is possible to show that the support of the bifurcation locus is not bounded.\\

\section{Density statements}

Our aim here is to explain why the support of the bifurcation measure may be considered as the locus of the strongest bifurcations. To this purpose we will  show that the remarkable parameters introduced in subsection \ref{ssRP} accumulate the support of the bifurcation measure.
These informations will be obtained through equidistribution arguments for the bifurcation current and its powers.\\
Let us mention that, using further techniques, some of the above mentioned remarkable parameters will be shown to belong to the support of the bifurcation measure in the next section.

\subsection{ Strongly Misiurewicz parameters}

The results given in this subsection are essentially due to Dujardin and Favre.
We present them in the setting of polynomial families and refer the reader to the original paper (\cite{DF}) for a greater generality.

\begin{theo}

In the degree $d$ polynomial family  
$\big(P_{c,a}\big)_{(c,a)\in {\tiny {\bf C}^{d-1}}}$ let us  define a sequence of analytic sets
by:
\begin{center}
$W_{n_0,\cdot\cdot\cdot,n_l}:=\cap_{j=0}^{l}\{P_{c,a}^{n_j}(c_j)=P_{c,a}^{k(n_j)}(c_j)\}$
\end{center}
where $l\le d-2$ and $k(n_j)<n_j$.
Then 
\begin{center}
$\lim_{n_{d-2}\to\infty}\cdot\cdot\cdot\lim_{n_{0}\to\infty}\frac{1}{d^{n_{d-2}}+\cdot\cdot\cdot+d^{n_{0}}}[W_{n_0,\cdot\cdot\cdot,n_{d-2}}]=\Bm$
\end{center}
and $W_{n_0,\cdot\cdot\cdot,n_{d-2}}$ is finite.
\end{theo}

\proof
We treat the case $d=3$ which is actually not very different from the general case.\\

Let us first observe that $W_{n_0,n_1}$ has codimension at least two and is contained in the connectedness locus which is compact (see Theorem \ref{controlinfty}).
Thus $W_{n_0,n_1}$ is a finite set.\\
 
Applying a version of Theorem \ref{TheoDF} suitably adapted 
 to the family $W_{n_0}$ yields 
 \begin{eqnarray*}
 \lim_{n_1\to\infty} d^{n_1} [W_{n_0,n_1}]=T_1\wedge [W_{n_0}]
 \end{eqnarray*}
 where $T_1$ is the activity current of the critical point $c_1$.
 By the same Theorem one has $ \lim_{n_0\to\infty} d^{n_0} [W_{n_0}] =T_0$ where $T_0$ is the activity current of $c_0$
 and this, since $T_1$ has continuous potentials, gives
  \begin{eqnarray*}
 \lim_{n_0\to\infty} T_1\wedge d^{n_0} [W_{n_0}]=T_1\wedge T_0.
 \end{eqnarray*}
 The conclusion follows immediately since, according to Theorem \ref{theoNIAC}, $\Bm=T_0\wedge T_1$.\qed\\
 
 An important consequence of the above result is that the support of the bifurcation measure is accumulated by strongly Misiurewicz polynomials.
 An alternative proof of that fact will be given in the next subsection for arbitrary families.
 We refer to \ref{defiMis} for a definition of Misiurewicz parameters.
 
 \begin{cor} In polynomial families,
the support of the bifurcation measure  is contained in the closure of strongly Misiurewicz parameters:  $Supp\;\Bm\subset\overline{\SMis}$.
 \end{cor}
 
 \proof
 By the above Theorem
 \begin{eqnarray}\label{DFF}
\lim_{n_{d-2}\to\infty}\cdot\cdot\cdot\lim_{n_{0}\to\infty}\frac{1}{d^{n_{d-2}}+\cdot\cdot\cdot+d^{n_{0}}}
[\cap_0^{d-2}\{P_{c,a}^{n_j}(c_j)=P_{c,a}^{n_j -1}(c_j)\}]=\Bm.
\end{eqnarray}

Let us observe that
\begin{eqnarray*}
 H_j:=\{P_{c,a}^{n_j}(c_j)=P_{c,a}^{n_j -1}(c_j)\}=Preper_{n_j} \cup Fix_j
 \end{eqnarray*}
 where, for parameters in  $Preper_{n_j}$ the critical point $c_j$ is strictly preperiodic to a (necessarily) repelling fixed point while $c_j$ is fixed
 for parameters in $Fix_j$.

  Now Theorem \ref{TheoDF} may be rewritten as
  
 \begin{eqnarray}
\lim_{n_{j}\to\infty}\frac{1}{d^{n_{j}}}
[Preper_{n_j}] + \frac{\alpha_{n_j}}{d^{n_{j}}} [Fix_j]=T_j
 \end{eqnarray} 

but, as  $T_j$ cannot charge the hypersurface $Fix_j$, we must have $\frac{\alpha_{n_j}}{d^{n_{j}}}\to 0$.
Thus $\lim_{n_{j}\to\infty}\frac{1}{d^{n_{j}}}[Preper_j] =T_j$ and \ref{DFF} yields

\begin{eqnarray*}
\lim_{n_{d-2}\to\infty}\cdot\cdot\cdot\lim_{n_{0}\to\infty}\frac{1}{d^{n_{d-2}}+\cdot\cdot\cdot+d^{n_{0}}}
[\cap_{j=0}^{d-2} Preper_{n_j}]=\Bm.
\end{eqnarray*}

The conclusion follows immediately since $\cap_{j=0}^{d-2} Preper_{n_j} \subset \SMis$.\qed\\
 
\subsection{Shishikura or hyperbolic parameters}\label{densShiHyp}

We aim here to show that the support of the bifurcation measure in $Mod_d$ is simultaneously accumulated by Shishikura and hyperbolic parameters
(see subsection \ref{ssRP} for definitions): 
\begin{center}
$Supp\;\Bm\subset \overline{\Shi} \cap\overline{\Hyp}$.
\end{center}

It is worth emphasize that both statements will be  deduced in the same way from the following generalized version of Theorem \ref{theodistaver}.
The results discussed in this subsection and the next one are due to Bassanelli and the author (\cite{BB1},\cite{BB2}).

\begin{theo}\label{AppMu} Let $\Bm$ be the bifurcation measure 
of a holomorphic family $(f_{\la})_{\la\in M}$ of rational maps. Let $m$ denote the complex dimension of $M$.
Let $0<r\le 1$. Then
there exists increasing sequences of integers $k_2(n),...,k_m(n)$ such that:
$$\Bm=\lim_n \frac{d^{-(n+k_2(n)+\cdot\cdot\cdot+k_m(n))}}{m! (2\pi)^m}\int_{[0,2\pi]^m}
[Per_{n}(r e^{i\theta_1})]\wedge\bigwedge_{j=2}^m[Per_{k_j(n)}(r e^{i\theta_j})]\;d\theta_1\cdot\cdot\cdot d\theta_m.$$
Moreover, we may assume that $k_j(n)\ne k_i(n)$ when $i\ne j$. 
\end{theo}

We will derive that result from Theorem \ref{theodistaver} by simple calculus arguments with currents.\\

\proof
For any fixed variety $Per_p(r e^{i\theta_p})$,
the set of $\theta\in [0,2\pi]$ for which $Per_p(r e^{i\theta_p})$
shares a non trivial component with $Per_m(r e^{i\theta})$ for some $m\in \NN^*$
is at most countable. This follows from Fatou's theorem on the
finiteness of the set of non-repelling cycles. Thus, the wedge 
products $[Per_{n}(r e^{i\theta_1})]\wedge [Per_{k_2(n)}(r e^{i\theta_2})]\cdot\cdot\cdot\wedge [Per_{k_m(n)}(r e^{i\theta_m})]$
 make sense for almost every $(\theta_1,\cdot\cdot\cdot,
\theta_m)\in[0,2\pi]^m$ and the integrals
\begin{center}
$\int_{[0,2\pi]^m}
[Per_{n}(r e^{i\theta_1})]\wedge\bigwedge_{j=2}^m[Per_{k_j(n)}(r e^{i\theta_j})]\;d\theta_1\cdot\cdot\cdot d\theta_m$
\end{center}
 are well defined.\\

Next, we need the following formula which 
has been justified for $q=1$ at the end of the proof of Theorem \ref{theodistaver}.
The proof is similar for $q>1$ and we shall omit it.
Recall that $L_n^r(\la):=\frac{d^{-n}}{2\pi}\int_0^{2\pi} \ln \vert p_n(\la,re^{i\theta})\vert\;d\theta$.

\begin{center}
 $dd^c L_{n_1}^r\wedge\cdot\cdot\cdot\wedge dd^c L_{n_q}^r=
\frac{d^{-(n_1+\cdot\cdot\cdot+n_q)}}{(2\pi)^q}
\int_{[0,2\pi]^q}\bigwedge_{k=1}^{q}[Per_{n_k}(r e^{i\theta_k})]d\theta_1\cdot\cdot\cdot
d\theta_q.$
\end{center}

To prove the convergence,
 we may replace $M$
by ${\bf C}^m$ since the problem is local.
The conclusion is obtained by using Theorem \ref{theodistaver}, the above formula
and the next Lemma inductively.

\begin{lem}\label{LemEgorov}
If $S_n\to (dd^c L)^p$ for some sequence $(S_n)_n$ of closed, positive $(p,p)$-currents on $M$ then  
 $dd^c L_{k(n)}^r\wedge S_n\to(dd^c L)^{p+1}$ for some increasing sequence of integers $k(n)$. 
\end{lem} 

 Let us briefly justify Lemma \ref{LemEgorov}.
Let us denote by $s_n$ the trace measure of $S_n$, as $M$ has been identified with ${\bf C}^m$ this measure is given by 
$s_n:=S_n\wedge (dd^c \vert z\vert^2)^{m-p}$. Since $S_n$ is positive, $s_n$ is positive as well.
Let us consider the sequence $(u_k)_k$ defined by $u_k:=L_k^r-L$.
We know that $(u_k)_k$ converges pointwise to $0$ (see remark \ref{ptLnr}) and is locally uniformly bounded (the function $L$ is continuous).
The positive current $S_n$ may be considered as a
$(p,p)$ form whose coefficients are measures which are dominated by the trace measure $s_n$.
Thus, by the dominated convergence theorem, $(L_k^r-L)S_n=u_k S_n$ tends to $0$ as $k\to\infty$ and $n$ is fixed.
On the other hand, $LS_n$ converges to $LS$ because $L$ is continuous. It follows that some subsequence $L_{k(n)}^rS_n$ converges to $LS$.
\qed\\

\begin{cor}\label{corShiHyp}
In the moduli space $Mod_d$ the support of the bifurcation measure $\Bm$ is contained in 
$\overline{\Shi} \cap\overline{\Hyp}$.
\end{cor}

\proof 
Use Remark \ref{helpwhen} to work with families and then apply Theorem \ref{AppMu}.
For $0<r<1$ one gets $Supp\;\Bm \subset \overline{\Hyp}$ and for $r=1$ $Supp\;\Bm \subset \overline{\Shi}$\qed\\

We also stress that the above result yields a rather simple proof of the existence of Shishikura maps, the original one was based on
quasi-conformal surgery.

\begin{cor}
In the moduli space $Mod_d$ one may find maps having $2d-2$ neutral cycles.

\end{cor}

\proof Combine Corollary \ref{corShiHyp} with the fact that $\Bm>0$ (see Proposition \ref{propLatMb}).\qed\\

Let us end this subsection with a final remark.\\

Theorem \ref{AppMu} remains true, with the same proof, if one replace the integrals by 
$$\int_{[0,2\pi]^m}
[Per_{n}(r_1 e^{i\theta_1})]\wedge\bigwedge_{j=2}^m[Per_{k_j(n)}(r_j e^{i\theta_j})]\;d\theta_1\cdot\cdot\cdot d\theta_m$$
where $0<r_j\le 1$ for $1\le j\le m$. As a consequence, if $\alpha+\nu=m$ and 
 ${\cal P}_{\alpha,\nu}$ is  the set of parameters $\la$ such that $f_{\la}$ has $\alpha$ distinct attracting cycles and $\nu$ distinct neutral cycles
then $Supp\;\Bm$ is contained in the closure of ${\cal P}_{\alpha,\nu}$.

\subsection{Shishikura or hyperbolic parameters with chosen multipliers}

As Theorem \ref{theoequipol} shows, polynomial with a neutral cycle of a given multiplier are dense in  the support of the bifurcation current.
We
 believe that  a similar property is still true for the bifurcation measure which means that
 Shishikura parameters with arbitrarily fixed multipliers should be  dense in the support $\Bm$.
 The following result goes in this direction.

\begin{theo}\label{theoBB1}
Denote by $p(f)$ (resp. $s(f)$, $c(f)$) the number of distinct parabolic (resp. Siegel, Cremer) cycles of $f\in Mod_d$.
Then
\begin{center} 
$Supp\;\Bm\subset \overline{\{f\in Mod_d\;/\; p(f)=p, s(f)=s\;\textrm{and}\;c(f)=c\}}$
\end{center}
for any triple of integers $p$, $s$ and $c$ such that $p+s+c=2d-2$. 
\end{theo}

The proof is essentially based on Lemma \ref{lemBrDu} and Ma\~{n}\'e-Sad-Sullivan theorem. More precisely we will use the following

\begin{lem}\label{lemE} Let $E$ be a dense subset of $[0,2\pi]$. Then for any holomorphic family of degree $d$ rational map $\left( f_\la\right)_M$
the set
\begin{center}
$\cup_{n}\cup_{\theta\in E} Per_n(e^{i\theta})$
\end{center}
is dense in the bifurcation locus.
\end{lem}

\proof Use Ma\~{n}\'e-Sad-Sullivan theorem or Theorem \ref{theodistaver} with $r=1$.\qed\\

Let us now prove Theorem \ref{theoBB1}. We restrict ourself to $Mod_2$. The general case requires to use a slicing argument, we refer to \cite{BB1} for
details.\\

\proof
Let $E_1$ and $E_2$ be two dense subsets of $[0,2\pi]$. Let $\la_0$ be a point in the support of $\Bm$ and $U_0$ be an arbitrarily small neighbourhood of $\la_0$.
By Lemma \ref{lemE}, the support of the bifurcation current $dd^c L$ is accumulated by holomorphic discs contained in $\cup_{n}\cup_{\theta\in E_1} Per_n(e^{i\theta})$.
Among such discs, let us consider those which go through $U_0$ and pick one disc $\Gamma_1$  on which the Lyapunov function $L$ is not harmonic. Such a disc exists
since otherwise, according to Lemma \ref{lemBrDu}, the measure
$\Bm$ would vanish on $U_0$. The bifurcation locus of $\left(f_\la\right)_{\Gamma_1}$ is not empty and thus,  to get a Shishikura parameter in $U_0$ with multipliers $e^{i\theta_1}$ and $e^{i\theta_2}$ where $\theta_j\in E_j$, it suffices to apply again 
Lemma \ref{lemE}  with the dense set $E_2$ to the family $\left(f_\la\right)_{\Gamma_1}$.\qed\\ 

Using Theorem \ref{theodisatt} instead Theorem \ref{theodistaver} one may prove, with exactly the same arguments as above, that the support of $\Bm$ is accumulated by
hyperbolic parameters with attracting cycles of given multipliers.

\begin{theo}
Let $w_1, w_2,\cdot\cdot\cdot,w_{2d-2}$ be complex numbers belonging to the open unit disc.
Then any $\la_0 \in Supp\;\Bm$ is accumulated by maps $f\in Mod_d$ having 
$2d-2$ attracting cycles whose respective multipliers are the $w_i$.
\end{theo}

\section{The support of the bifurcation measure}

We  will establish that the inclusions obtained in the former section are actually equalities. This will give a precise meaning to our interpretation  of the support of the bifurcation measure in $Mod_d$ as a
 strong-bifurcation locus.

\subsection{A transversality result}\label{ssTran}

Transversality statements play a very important role for understanding the structure of parameter spaces. We have already encounter such results
like for instance Lemma \ref{lemTrDH} or the Fatou-Shishikura inequality. We refer to the fundamental work of Epstein \cite{epstein2} for a general and synthetic treatement of transversality
problems in holomorphic dynamics.
\\

All the results presented here are true in $Rat_d$ or in the moduli spaces $Mod_d$. For simplicity we shall restrict ourself to the moduli space
$Mod_2$ which will be treated as a holomorphic family $\left(f_\la\right)_{\la \in \CC^{2}}$ (see Theorem \ref{theoMiMod2}). We shall also assume, to simplify,
to have two marked critical points $c_j(\la)\;j=1,2$. Our  exposition also covers the case of the family of cubic polynomials.\\

Assume that $f_0\in Mod_2$ is strongly Misiurewicz. This means that there exists two repelling cycles
\begin{center}
${\cal C}_j(0):=\{z_j(0),\cdot\cdot\cdot,f_0^{n_j-1} (z_j(0))\}$ 
\end{center}
and an integer $k_0\ge 1$ such that
\begin{center}
$f_0^{k_0} (c_j(0)) =z_j(0)$ but $c_j(0)\notin {\cal C}_j(0)$
\end{center}
 for $j=1,2$.\\

By the implicit function theorem, we may follow the cycles ${\cal C}_j$ on a small ball $B(0,r)$ centered at the origin. Writting 
${\cal C}_j(\la):=\{z_j(\la),\cdot\cdot\cdot,f_\la^{n_j-1} (z_j(\la))\}$ the cycles corresponding to the parameter $\la\in B(0,r)$,
we may define an important tool for studying the parameter space near $f_0$.

\begin{defi}\label{defiAct}
The map $\chi:B(0,r)\to\CC^2$ 
 defined by 
 \begin{center}
$ \la\mapsto \left(f_\la^{k_0}(c_j(\la))-z_j(\la)\right)_{j=1,2}$
 \end{center}
 is called \emph{activity map} near the strongly Misiurewicz parameter $f_0$. 
\end{defi}

The activity map $\chi$ measures the difference between two natural holomorphic motions of the point $z_j(0)$. One as a repelling periodic point and the other as a post-critical point.
The definition of $\chi$ implicitely uses local charts of $\pp$ near $z_j(0)$.
As we shall see later, the activity map will allow to transfer informations from the dynamical space of $f_0$ to the parameter space.
For this, the next result will be essential.

\begin{theo}\label{theoTBE}
The activity map $\chi$ near a strongly Misiurewicz parameter $f_0$ is locally invertible.
\end{theo}

This Theorem was proved by Buff and Epstein (see \cite{BuEp}) in the general setting of $Rat_d$. In that case one has to assume that $f_0$ is
not a flexible Latt\`es map (such maps do not exist in degree two). The proof of Buff and Epstein uses quadratic differentials thechniques.
We shall prove here a weaker statement which is due to Gauthier (\cite{Gau}) and is sufficent for the applications we have in mind.

 \begin{theo}\label{theoTG}
The activity map $\chi$ near a strongly Misiurewicz parameter $f_0$ is locally proper.
\end{theo}

The proof of that result is based on more classical arguments going back to Sullivan
(see also \cite{vanstrien} or \cite{Aspenberg2}). The key point relies on the following Lemma.

\begin{lem}\label{lemSul}
Any holomorphic curve contained in $\SMis$ must consist of flexible Latt\`es maps. 
\end{lem}
\proof
 Assume that $\left(f_\la\right)_{\la\in D}$ is a holomorphic family parametrized by some one-dimensional disc $D$ and that all $f_\la$
are strongly Misiurewicz parameters when $\la \in D$.  It actually suffices to show that $f_{\la_1}$ is Latt\`es for some $\la_1\in D$.\\

Since all $f_\la$ are strongly Misiurewicz, the Julia set of $f_\la$ coincides with $\pp$ for all $\la\in D$. Moreover,  by Lemma \ref{lemappperno}, the  family $\left(f_\la\right)_{\la\in D}$ is stable. According to a Theorem of Ma\~{n}\'e-Sad-Sullivan (see \cite{MSS} Theorem B), 
the stability of $\left(f_\la\right)_{\la\in D}$ implies the existence of  
a quasiconformal holomorphic motion $\Phi: D\times \pp\to \pp$ which conjugates $f_\la$ to $f_0$ on $\pp$. Let us denote by $\eta^\la$ the Beltrami form satisfying
\begin{center}
$\frac{\partial \Phi_\la}{\partial \bar z}=\eta^\la\frac{\partial \Phi_\la}{\partial z}$.
\end{center}
There exists $\la_1\in D\setminus\{0\}$ for which the support of $\eta^{\la_1}$ has strictly positive Lebesgue measure. Indeed, if this would not be the case, 
$f_\la$ would be holomorphically conjugated to $f_0$ for all $\la\in D$. Then the Julia set of $f_{\la_1}$ carries an invariant line field and thus $f_{\la_1}$ is a
flexible Latt\`es map. This last argument uses the fact that the conical set of a strongly Misiurewicz map coincides with its Julia set  (see \cite{BM} Theorem VII. 22 and \cite{McMullen} corollary 3.18). \qed\\

We may now easily prove Theorem \ref{theoTG}\\

\proof
Let us first establish that at least one critical point must be active. Assume to the contrary that both critical points are passive around $f_0$. Then, according to 
Lemma \ref{lemAcMi}, $f_\la$ is strongly Misiurewicz for all $\la\in B(0,r)$ after maybe reducing $r$. Cutting $B(0,r)$ by a disc $D$ passing through the origin we 
obtain, by Lemma \ref{lemSul}, a disc of flexible Latt\`es maps.  Since this is impossible in $Mod_2$ (and by assumption in other cases)
we have reached a contradiction and proved that at least one critical point, say $c_1$, is active at $f_0$.\\

The activity of $c_1$ means that $\chi_1^{-1}(0)$ has codimension one. The conclusion is obtained by repeating the argument on the hypersurface
$\chi_1^{-1}(0)$.\qed\\

\subsection{The bifurcation measure and strong-bifurcation loci}

We want to establish that the inclusion $Supp\;\Bm\subset \overline{\Shi}\cap\overline{\SMis}$ obtained in subsection 
\ref{densShiHyp} is actually an equality.  This is the reason for which we shall consider the support of the bifurcation measure as a
strong-bifurcation locus. This is essentially
a consequence of the following result due to Buff and Epstein \cite{BuEp}.

\begin{theo}\label{theoSBE}
In the moduli space $Mod_d$ the set of strongly Misiurewicz parameters is contained in the support of the bifurcation measure:
$\SMis\subset Supp\;\Bm$.
\end{theo}

\begin{cor}
In the moduli space $Mod_d$ one has $Supp\;\Bm=\overline{\Shi}=\overline{\SMis}$.
\end{cor}

\proof
To simplify the presentation we will work in the degree $3$ polynomial family
$\big(P_{c,a}\big)_{(c,a)\in {\tiny {\bf C}^{2}}}$.
As usual we write $\la$ the parameter $(c,a)$ and $c_1(\la)$, $c_2(\la)$ the marked critical points, the fact that in this setting $c_1=0$ does not play any role here.\\

Assume that $p_0$ is a strongly Misiurewicz polynomial. By definition, there exists an integer $k_0$ such that $p_0^{k_0}(c_j(0))=:z_j(0)$ is a repelling periodic point for $j=1,2$. To get lighter notations we whall assume that the $z_j(0)$ are fixed repelling points.\\

We denote by $z_j(\la)$ the repelling fixed points which are obtained by holomorphically moving $z_j(0)$ on some neighbourhood of the origin
and by $w_j(\la)$ the corresponding multipliers. Observe that $\vert w_j(\la)\vert \ge a>1$ on a sufficently small neighbourhood of $0$.\\

The activity map $\chi$ (see definition \ref{defiAct}) may be written: $\chi=\left(\chi_1,\chi_2\right)$ where
\begin{center}
$\chi_j(\la)=p_\la^{k_0}\left( c_j(\la)\right) - z_j(\la)$.
\end{center}

We will use here Theorem \ref{theoTBE} and assume that $\chi$ is locally invertible at the origin. 
It is possible to adapt the proof for using the weaker transversality statement given by Theorem \ref{theoTG}. For this one uses the fact that the sets
obtained by rescaling the ramification locus of $\chi$ are not charged by the measure $\Bm$ and,  thanks to  some Besicovitch covering argument,
reduces the problem to some estimate similar to those which we will now perform in the invertible case. We refer to the papers \cite{BuEp} and \cite{Gau} for details.\\

Let us denote by $D^2(0,\epsilon)$ the bidisc centered at the origin and of multiradius $(\epsilon,\epsilon)$ in ${\bf C}^{2}$. For $\epsilon$ small enough we may define a sequence of rescaling 
\begin{center}
$\delta_n: D^2(0,\epsilon) \to \Omega_n$
\end{center}

by setting $\delta_n(x):=\chi^{-1} \left(\frac{x_1}{\left(w_1(0)\right)^n},\frac{x_2}{\left(w_2(0)\right)^n}\right)$.
To prove the Theorem, it suffices to show that $\Bm(\Omega_n)>0$ for all $n$. The crucial point of the proof is revealed by the following computations.

\begin{eqnarray*}\label{venezia}
2\Bm(\Omega_n)=\int_{\Omega_n} \Bc^2=\int_{D^2(0,\epsilon)} \delta_n^{\star} \left(T_1+T_2\right)^2
\ge \int_{D^2(0,\epsilon)} \delta_n^{\star} \left( T_1\wedge T_2\right)=\\
=\int_{D^2(0,\epsilon)} \delta_n^{\star} \left[dd^c g_\la(c_1(\la))\wedge dd^c g_\la(c_2(\la))\right]
\end{eqnarray*}

using the homogeneity property of the Green function 
\begin{eqnarray*}
g_\la\left(c_j(\la)\right)=3^{-(k_0+n)} g_\la\left( p_\la^{k_0+n}(c_j(\la))\right)
\end{eqnarray*}
one thus gets
 
 \begin{eqnarray*}
 2\cdot 3^{(k_0+n)} \Bm(\Omega_n)\ge
\int_{D^2(0,\epsilon)} \delta_n^{\star} \left[ dd^c g_\la\circ p_\la^{k_0+n}(c_1(\la))\wedge dd^c g_\la\circ p_\la^{k_0+n}(c_2(\la))\right]=\\
= \int_{D^2(0,\epsilon)} dd^c g_{\delta_n(x)}\circ p_{\delta_n(x)}^{k_0+n}(c_1({\delta_n(x)}))\wedge dd^c g_{\delta_n(x)}\circ p_{\delta_n(x)}^{k_0+n}(c_2({\delta_n(x)})).
 \end{eqnarray*}
 
Let us express the quantities $p_{\delta_n(x)}^{k_0+n}(c_j({\delta_n(x)}))$ by using the activity map $\chi$. By definition we have 
$p_\la^{k_0+n}(c_j(\la))=p_\la^{n}\left(z_j(\la)+\chi_j(\la)\right)$ and thus
\begin{eqnarray*}
p_{\delta_n(x)}^{k_0+n}(c_j({\delta_n(x)}))=p_{\delta_n(x)}^{n}\left(z_j({\delta_n(x)})+\frac{x_j}{w_j(0)^n}\right).
\end{eqnarray*}

To conclude, we momentarily admit the following\\

{\bf Claim}: 
$p_{\delta_n(x)}^{n}\left(z_j({\delta_n(x)})+\frac{x_j}{w_j(0)^n}\right)$ is uniformly converging to some local biholomorphism
$\psi_j:\CC,_{0}\to\pp,_{z_j(0)}$.\\

As the Green function $g_\la(z)$ is continuous in $(\la,z)$, the Claim clearly implies that $g_{\delta_n(x)}\circ p_{\delta_n(x)}^{k_0+n}(c_2({\delta_n(x)}))$
uniformly converges towards $g_0\left(\psi_0(x_j)\right)$ and our estimate yields
\begin{eqnarray*}
\liminf_n  \Bm (\Omega_n) \ge \int_{D^2(0,\epsilon)} dd^c g_{0}\left(\psi_1(x_1)\right)\wedge dd^c g_{0}\left(\psi_2(x_2)\right)=\\
=\left( \int_{\psi_1 \left(D(0,\epsilon)\right)}  dd^c g_{0}\right) \left( \int_{\psi_2 \left(D(0,\epsilon)\right)}  dd^c g_{0}\right) >0
\end{eqnarray*} 
 where the positivity of the last term follows from the fact that the repelling fixed points $z_j(0)=\psi_j(0)$ belong to the Julia set of $p_0$. \\
 
 It remains to justify the Claim. Let us write $w_j(\la)$ on the form $w_j(0)\left(1+\epsilon_j(\la)\right)$. As $\Vert \delta_n(x)\Vert \le C \frac{1}{a^n}$ one sees that
 $\left( 1+\epsilon_j(\delta_n(x))\right)^n$ is uniformly converging to $1$. Then
 \begin{eqnarray*}
 p_{\delta_n(x)}^{n}\left(z_j({\delta_n(x)})+\frac{x_j}{w_j(0)^n}\right)=
 p_{\delta_n(x)}^{n}\left[z_j({\delta_n(x)})+\frac{x_j}{w_j(\delta_n(x))^n}\left( 1+\epsilon_j(\delta_n(x))\right)^n\right]
 \end{eqnarray*}
 behaves like 
  $p_{\delta_n(x)}^{n}\left[z_j({\delta_n(x)})+\frac{x_j}{w_j(\delta_n(x))^n}\right]$.\\

  Now let us linearize $p_\la$ near the repelling fixed points $z_j(\la)$. The linearization holomorphically depends on the parameter $\la$ and one gets
  local biholomorphisms $\psi_{j,\la}$ such that $\psi_{j,\la}(0)=z_j(\la)$ and 
  \begin{center}
  $p_\la \circ \psi_{j,\la} (z)=\psi_{j,\la} (w(\la) z)$ on $B(0,\frac{\epsilon}{\vert w(\la)\vert})$.
  \end{center}
 Then the local biholomorphism $\psi_j$ of the Claim is simply $\psi_{j,0}$.\qed\\
 
 \begin{rem}
By Theorem \ref{theoNIAC}, $(T_1+T_2)^2=2T_1\wedge T_2$ and therefore the estimate in the proof of Theorem \ref{theoSBE}
becomes an equality:

 \begin{eqnarray*}
 2\cdot 3^{(k_0+n)} \Bm(\Omega_n)
= \int_{D^2(0,\epsilon)} dd^c g_{\delta_n(x)}\circ p_{\delta_n(x)}^{k_0+n}(c_1({\delta_n(x)}))\wedge dd^c g_{\delta_n(x)}\circ p_{\delta_n(x)}^{k_0+n}(c_2({\delta_n(x)})).
 \end{eqnarray*}
This allows to estimate the pointwise Hausdorff dimension of $\Bm$ at strongly Misiurewicz parameters.
\end{rem}

Let us recall  that, for polynomial families, Theorem \ref{theoSBE} was first proved by Dujardin and Favre in \cite{DF}. Their proof is based on a totally different approach, which consists in comparing the measure $\Bm$ with a landing measure for a family of external rays and uses some deep result due to Bielefeld, Fisher and Hubbard \cite {BFH}
and Kiwi \cite{Kiwi}.\\
This approach has the advantage to yield to a generalization in any degree of some results due to Graczyk and Swiatek \cite{GS} and Smirnov \cite{Sm}. They prove the following.

\begin{theo}
A $\Bm$ generic polynomial has the topological Collet-Eckman property. In particular the following  properties are satisfied for $\Bm$-almost all parameters:
\begin{itemize}
\item[i)] all cycles are repelling
\item[ii)] the Julia set and the filled-in Julia set coincide and their Hausdorff  dimension is  strictly less than $2$
\item[iii)] each critical point has a dense forward orbit in the Julia set.
\end{itemize}
\end{theo}

Let us end this subsection with a generalization of Buff-Epstein's theorem.
It is actually true that all Misiurewicz parameters and not only strongly ones (see \ref{defiMis} for definitions) belong to the support of the bifurcation measure.
This has been proved by Gauthier in \cite{Gau}, his result can be stated as follows:

\begin{theo}\label{theoGauM}
In the moduli space $Mod_d$ one has $\overline{\Mis}=Supp\;\Bm$.
\end{theo}

The proof is conceptually similar to that of Theorem \ref{theoSBE} but is technically more involved. 
Let us stress the main new difficulties. To define the activity map $\chi$ one has to construct a natural holomorphic motion for hyperbolic invariant set, this has 
been done by de Melo-van Strien (\cite{dMvS}, Theorem 2.3 page 225). The transversality statement concerning $\chi$ is a weaker one but, as we explained in the last subsection, one may overcome this by using the flexibility of positive currents. The final argument requires to linearize the maps along repelling orbits contained in hyperbolic sets.\\

Using the transversality map $\chi$ associated to
Misiurewicz parameters  it is possible to construct a "transfer map" which copies hyperbolic sets in the parameter space with a precise control of distortion. This allows to estimate the dimension of  $Supp\;\Bm$ from below and is particularly useful when dealing with hyperbolic sets of big Hausdorf dimension. This is what we shall discuss in the last subsection.

\subsection{Hausdorff dimension estimates}

Our aim here is to present the following result obtained by T. Gauthier in his thesis (see \cite{Gau}).\\

\begin{theo}\label{theoGauD}
The strong-bifurcation locus in $Mod_d$ is homogeneous and has full Hausdorff dimension.
\end{theo}

Let us recall that a subset $E$ of any metric space is said to be \emph{homogeneous} if its Hausdorff dimension coincides with that of any non-empty intersection
$U\cap E$ with an open set $U$.\\

For the quadratic polynomial family, the strong-bifurcation locus coincides with the usual one.
In that case, Gauthier's theorem restates a celebrated result of Shishikura \cite{Shishikura2} which says that  the boundary of the Mandelbrot set is of
Hausdorff dimension $2$. More generally, Tan Lei \cite{TanLei} has proved that the boundary of the connectedness locus of the degree $d$ polynomial family has maximal Hausdorff dimension.
As Remark \ref{remShi} shows, Gauthiers's result asserts that, for polynomial families, the full dimension is actually concentrated on a tiny part of the  boundary of the connectedness locus.\\

The approach of Gauthier exploits his former result (Theorem \ref{theoGauM}) and the transfer techniques elaborated by Shishikura and Tan-Lei
which he slightly simplifies. We shall explain this below.\\

Combining Theorem \ref{theoGauD} and Theorem \ref{theoBB1}  yields the

\begin{cor}
Denote by $p(f)$ (resp. $s(f)$, $c(f)$) the number of distinct parabolic (resp. Siegel, Cremer) cycles of $f\in Rat_d$.
Let $p$, $s$ and $c$ be three integers such that $p+s+c=2d-2$. 
Then the set 
\begin{center}
$\overline{\{f\in Rat_d\;/\; p(f)=p, s(f)=s\;\textrm{and}\;c(f)=c\}}$
\end{center}
is homogeneous and has maximal Hausdorff dimension $2(2d+1)$.
\end{cor}

\par Let us now sketch the proof of Theorem \ref{theoGauD} in the moduli space $Mod_2$ of quadratic rational maps which we again will consider as a holomorphic family parametrized by $\CC^2$ with marked critical points.\\

We recall that the hyperbolic dimension $\dim_{\textup{hyp}}(f)$ of a map $f$ is the supremum of the dimensions of all compact $f$-invariant hyperbolic homogeneous sets.
We shall actually prove the following basic dimension estimate:

\begin{prop}\label{inegdim}
For  $f\in \Mis$  one has
$\dim_H(\Mis\cap \Omega)\geq 2 \dim_{\textup{hyp}}(f)$ for any neighborhood $\Omega$ of $[f]$ in $Mod_2$.
\end{prop}

It is absolutely not trivial to deduce Theorem \ref{theoGauD} from the above Proposition. Adapting some ideas of McMullen
 \cite{McMullen2}, one may show that the hyperbolic dimension of a Misiurewicz map equals the Hausdorff dimension of its Julia set which is
 the Riemann sphere. Then Proposition \ref{inegdim} actually gives $\dim_H(M_{is})=4$ and the Theorem \ref{theoGauD} follows
since, by Theorem \ref{theoGauM}, $\Mis\subset Supp \ \Bm$.\\

Let us now prove Proposition \ref{inegdim}. Let  $f$ be Misiurewicz (and not a flexible Latt\`es map) and  $\Omega$ be a neighbourhood of $f$ in $Mod_2$. Let $E_f$ be a compact $f$-hyperbolic set such that $$\dim_H(E_f\cap U)\geq\dim_{\textup{hyp}}(f)-\epsilon$$ for any open set $U$ intersecting $E_f$.
Note that the critical orbits of $f$ are captured by some hyperbolic set which, a priori, is different from that set $E_f$. The proof uses two technical tools which are certainly of independant interest.
 The first one allows to make all cititical orbits fall into some \emph{big} hyperbolic set, which is a deformation of $E_f$, by an arbitrarily small perturbation of $f$. More specifically one shows the following :\\

{\bf Fact 1} {\it   One may find a Misiurewicz map $g\in \Omega$ whose critical points eventually fall into $E_g$ (i.e. $g^N(c_i)\in E_g$ for $i=1,2$) where 
 $E_g$ is a compact $g$-hyperbolic set such that 
 \begin{center}
 $\dim_H(E_g\cap U_i)\geq\dim_{\textup{hyp}}(f)-3\epsilon$
 \end{center}
  for any open neighbourhood $U_i$ of $g^N(c_i)$.}\\
 
 The second is a construction of a transfer map with good regularity properties:\\
 
{\bf Fact 2} {\it Let $g$ be a Misiurewicz map with associated $g$-hyperbolic set $E_g$ (i.e. $g^N(c_i)\in E_g$ for $i=1,2$).
For any neighbourhood $\Omega$ of $g$ in $Mod_2$ there exists a neighbourhood $U_1\times U_2$ of $\left(g^N(c_1),g^N(c_2)\right)$ in $\pp\times\pp$ and a transfer map ${\cal T}$
\begin{center}
${\cal T} : \left((E_g\cap U_1)\times (E_g\cap U_2)\right) \to  \Mis\cap \Omega$
\end{center}
such that:
\begin{center}
$\|{\cal T}(z_1,z_2)-{\cal T}(z_1',z_2')\|\geq C_\epsilon\|(z_1,z_2)-(z_1',z_2')\|^{\frac{1+\epsilon}{1-\epsilon}}$
\end{center}
 for $0<\epsilon<\epsilon_0$ and $(z_1,z_2),(z_1',z_2')\in U_1\times U_2$.}\\

Let us stress that the proofs of both Facts require transversality (Theorem \ref{theoTG}) and uses the H\"older continuity properties of holomorphic motions (Lemma \ref{HoHoMo}).
We shall prove them later and right now
 see how they easily lead to the desired dimension estimate.\\

Let $g \in \Omega$ be the Misiurewicz map given by Fact 1. Consider the transfer map 
${\cal T}$ associated to $g$ and given by Fact 2. Then, by construction, we have:
\begin{center}
$\begin{array}{rcl}
\dim_H(\Mis\cap \Omega) 
& \geq & \dim_H\left[{\cal T}\big((E_g\cap U_1)\times (E_g\cap U_2)\big)\right].
\end{array}$
\end{center}
The estimate of Fact 2 then gives:
\begin{center}
$\begin{array}{rcl}
\dim_H(M_{is}\cap \Omega) & \geq & \frac{1-\epsilon}{1+\epsilon}\dim_H\left[(E_g\cap U_1)\times (E_g\cap U_2)\right]\\
& \geq & \frac{1-\epsilon}{1+\epsilon}\big(\dim_H (E_g\cap U_1)+\dim_H(E_g\cap U_2)\big).
\end{array}$
\end{center}
But the choice of the map $g$ (see fact 1) gives
\begin{center}
$\begin{array}{rcl}
\dim_H(M_{is}\cap W) & \geq & \frac{1-\epsilon}{1+\epsilon}2(\dim_{\textup{hyp}}(f)-3\epsilon).
\end{array}$
\end{center}
and the conclusion follows  making $\epsilon \to 0$.\\

We now prove the two Facts.\\

\textbf{Proof of Fact 1.}
 One may construct a holomorphic motion
\begin{center}
$h:B(0,\rho)\times E_f\cup\overline{\bigcup_{n\geq k_0}f^n(C_f)}\longrightarrow\mathbf{P}^1$
\end{center}
 such that $f_\lambda\circ h_\lambda=h_\lambda\circ f_0$ on $E_f$. Here $B(0,\rho)$ is a ball in $\mathbf{C}^2$ centered at $f_0:=f$. 
 As in Definition \ref{defiAct}, one may associate to this motion an activity map
 $\chi$ defined on a neighbourhood of $0$ 
 by $$ \chi : \la\mapsto \left(f_\la^{k_0}(c_j(\la))-h_\la(f_{0}^{k_0}(c_j(0)))\right)_{j=1,2}.$$

 The transversality Theorem \ref{theoTG} being also valid in this context, the hypersurfaces $\chi_1^{-1}\{0\}$ and  $\chi_2^{-1}\{0\}$ intersect properly and, therefore, the critical point $c_2$ is active in  $\chi_1^{-1}\{0\}$. This will allow us to drastically change the orbit of $c_2$ by a small perturbation of $\la$ in  $\chi_1^{-1}\{0\}$. \\
 
  Let $z_1$, $z_2$ and $z_3$ be three distinct points of $E_f$. By Montel's Theorem, there exists $1\leq i\leq 3$, a parameter $\la_1\in \chi_1^{-1}\{0\}$ which is arbitrarily close to $0$ and an integer $N_1\geq1$ such that $f_{\lambda_1}^{N_1}(c_2(\lambda_1))=h_{\lambda_1}(z_i)$. When $f_{\la_1}$ is sufficently close to $f$ (i.e. $\la_1$ close to $0$) one has
  \begin{center}
    $\dim_H(h_{\lambda_1}(E_f)\cap U)\geq\dim_{\textup{hyp}}(f)-2\epsilon$
    \end{center}
     for any small enough neighborhood $U$ of $z_i$. This follows from the H\"older regularity of holomorphic motions.\\
     
 It turns out, since $\la_1$ has been taken in $\chi_1^{-1}\{0\}$, that the map $f_{\la_1}$ is again Misiurewicz with critical orbits captured by $E:=h_{\lambda_1}(E_f)$. We may therefore define a new activity map (still denoted $\chi$) near $\la_1$ and repeat the above procedure working this time in the family $\chi_2^{-1}\{0\}$.
This yields a 
 map $g:= f_{\lambda_2}$ which is arbitrarily close to $f$ and an integer $N\geq1$ such that
\begin{center}
$f_{\lambda_2}^N(C_{f_{\lambda_2}})\subset h_{\lambda_2}(E)$ and $\dim_H(h_{\lambda_2}(E)\cap U)\geq \dim_{\textup{hyp}}(f)-3\epsilon$
\end{center}
for any neighborhood $U\subset \mathbf{P}^1$ of any $c_i(\lambda_2)$.

~

\textbf{Proof of Fact 2.} 
As in the proof of Fact 1,
one may construct a holomorphic motion
\begin{center}
$h:B(0,\rho)\times E_g\longrightarrow\mathbf{P}^1$
\end{center}
 such that $g_\lambda\circ h_\lambda=h_\lambda\circ g_0$ on $E_g$. Here $B(0,\rho)$ is a ball in $\mathbf{C}^2$ centered at $g_0:=g$.\\
  
By a Theorem of Bers and Royden (see \cite{BersRoyden}), this motion extends to a holomorphic motion $$h:B(0,\rho/3)\times \mathbf{P}^1\longrightarrow\mathbf{P}^1.$$

 Set

\begin{eqnarray*}
X:B(0,\rho/3)\times (\mathbf{P}^1)^2 & \longrightarrow & (\mathbf{P}^1)^2\\
(\lambda,(z_1,z_2)) & \longmapsto & \big(g_\lambda^N(c_1(\lambda))-h_\lambda(z_1),g_\lambda^N(c_2(\lambda))-h_\lambda(z_2)\big).
\end{eqnarray*}
For the sake of simplicity, we assume that the activity map of $g$, which coincides with $X\big(\cdot,(g^N(c_1),g^N(c_2))\big)$, is locally invertible at $g$. Then, since $X$ is continuous,  a Rouch\'e-like Theorem shows that  the map $X\big(\cdot,(z_1,z_2)\big)$ is locally invertible at $g$ for $(z_1,z_2)$ close enough to $(g^N(c_1),g^N(c_2))$.\\

 This allows to define the transfer map $\cal T$ by
\begin{center}
${\cal T}(z_1,z_2):=X(\cdot,(z_1,z_2))^{-1}\{(0,0)\}$,
\end{center}
for $(z_1,z_2)\in U_1\times U_2$, where $U_i$ are  small enough neighborhood of $g^N(c_i)$. By construction ${\cal T}(z_1,z_2)$ is a Misiurewicz map
when $(z_1,z_2)\in  (E_g\cap U_1)\times (E_g\cap U_2)$ and thus,  shrinking $U_i$ one  has
\begin{eqnarray*}
{\cal T}\left((E_g\cap U_1)\times (E_g\cap U_2)\right)\subset \Mis\cap \Omega.
\end{eqnarray*}

The estimate is obtained by exploiting the fact that the holomorphic motion $h(\la,z)$ is
 H\"older continuous in $z$ and  holomorphic in $\la$. \qed

\bibliographystyle{plain}
\bibliography{biblio}
\end{document}